\documentclass[11pt]{amsart}

\usepackage{amsfonts,amssymb,amsmath,mathrsfs,graphicx}
\usepackage{float}
\usepackage{epsfig}

\usepackage{CJK,color}

\allowdisplaybreaks

\definecolor{red}{rgb}{1.00,0.00,0.00}
\definecolor{blue}{rgb}{0.00,0.00,0.63}
\definecolor{black}{rgb}{0.00,0.00,0.00}

\newcommand{\nm}{\nonumber}

\newtheorem{theorem}{Theorem}[section]
\newtheorem{lemma}{Lemma}[section]

\newtheorem{proposition}{Proposition}[section]
\newtheorem{remark}{Remark}[section]

\newcommand{\ba}{\begin{aligned}}
\newcommand{\ea}{\end{aligned}}
\newcommand{\be}{\begin{equation}}
\newcommand{\ee}{\end{equation}}
\renewcommand{\div}{ {\rm div }  }
\def\na{\nabla}

\def\charf {\mbox{{\text 1}\kern-.30em {\text l}}}

\def\mb{\mathbf}

\def\f{\frac}

\def\d{\delta}
\def\a{\alpha}

\def\z{\zeta}

\def\di{\displaystyle}
\def\r{\rho}

\def\k{\kappa}

\def\charf {\mbox{{\text 1}\kern-.30em {\text l}}}

\def\mb{\mathbf}
\def\f{\frac}

\def\d{\delta}
\def\s{\sigma}
\def\a{\alpha}

\def\di{\displaystyle}
\def\r{\rho}

\def\k{\kappa}
\def\pa{\partial}

\begin{document}

\title[Stability of planar shock wave]
{Nonlinear stability of planar viscous shock wave to three-dimensional compressible Navier-Stokes equations}


\author[Wang]{Teng Wang}
\address[Teng Wang]{\newline College of Mathematics, Faculty of Science, Beijing University of Technology, Beijing 100124, P. R. China}
\email{tengwang@amss.ac.cn}

\author[Wang]{Yi Wang}
\address[Yi Wang]{\newline  CEMS, HCMS, NCMIS, Academy of Mathematics and Systems Science, Chinese Academy of Sciences, Beijing 100190, P. R. China
\newline and School of Mathematical Sciences, University of Chinese Academy of Sciences, Beijing 100049, P. R. China
}
\email{wangyi@amss.ac.cn}




\maketitle

\begin{abstract}

We prove the nonlinear stability of the planar viscous shock up to a time-dependent shift for the three-dimensional (3D) compressible Navier-Stokes equations under the generic perturbations, in particular, without zero mass conditions. Moreover, the time-dependent shift function keeps the shock profile shape time-asymptotically. Our stability result is unconditional for the weak planar Navier-Stokes shock. Our proof is motivated by the $a$-contraction method (a kind of weighted $L^2$-relative entropy method) with time-dependent shift introduced in \cite{KV, KV-1, KVW} for the stability of viscous shock in one-dimensional (1D) case. Instead of the classical anti-derivative techniques, we perform the stability analysis of planar Navier-Stokes shock in original $H^2$-perturbation framework and therefore zero mass conditions are not necessarily needed, which, in turn, brings out the essential difficulties due to the compressibility of viscous shock. Furthermore, compared with 1D case, there are additional difficulties coming from the wave propagation along the multi-dimensional transverse directions and their interactions with the viscous shock. To overcome these difficulties, a multi-dimensional version sharp weighted Poincar\'{e} inequality (see Lemma \ref{le-poin}), $a$-contraction techniques with time-dependent shift, and some essential physical structures of the multi-dimensional Navier-Stokes system are fully used.

\end{abstract}

\maketitle \centerline{\date}

\tableofcontents

%
%
\section{Introduction}
\setcounter{equation}{0}
We are concerned with the time-asymptotic stability of planar viscous shock wave to 3D compressible Navier-Stokes equations
\begin{equation}\label{ns}
\left\{
\begin{array}{lr}
\di \pa_t\r+\div_x(\r u)=0,&(t,x)\in\mathbb{R}^+\times \Omega,\\[1mm]
\di \pa_t(\r u)+\div_x(\r u\otimes u)+\na_x {p}(\r)=\mu \Delta_x u+(\mu+\lambda)\na_x \div_x u.&
\end{array}
\right.
\end{equation}
Here $\r=\r(t,x):\mathbb{R}^+\times \Omega\rightarrow \mathbb{R}^+,$ $u=u(t,x)=(u_1,u_2,u_3)^t(t,x):\mathbb{R}^+\times \Omega\rightarrow \mathbb{R}^3$ represent the mass density and the velocity
of a fluid in $\Omega\subset\mathbb{R}^3$ respectively, and ${p}(\r)=b\r^{\gamma}(b>0,\gamma>1)$
stands for the classical $\gamma$-law pressure, and both constants $\mu$ and $\lambda$ are viscosity coefficients satisfying the physical constraints
$$
\mu>0,\quad 2\mu+3\lambda\geq0.
$$
We are concerned with Cauchy problem of 3D Navier-Stokes system \eqref{ns} in $x=(x_1, x_2, x_3)^t\in\Omega:={\mathbb R}\times {\mathbb T}^2$ with ${\mathbb T}^2:=(\mathbb {R}/\mathbb {Z})^2$. The initial data
\be\label{initial}
(\r, u)|_{t=0}=(\r_0,u_0)\to (\r_\pm,u_\pm), \quad {\rm as}\quad x_1\rightarrow \pm\infty,
\ee
is prescribed with the far-fields conditions $\r_{\pm}>0$ and $u_\pm=(u_{1\pm},0,0)^t$ as $ x_1\rightarrow \pm\infty$, and the periodic boundary conditions are imposed on $(x_2,x_3)\in\mathbb{T}^2$  for the solution $(\r, u)$.


\

The large-time asymptotic behavior of the solutions to 3D compressible Navier-Stokes system \eqref{ns}-\eqref{initial} with different end states $(\r_\pm,u_\pm)$ without shear is conjectured to be determined by the planar Riemann problem of corresponding 3D Euler system
\begin{equation}\label{3d-ER}
\left\{
\begin{array}{l}
\di \pa_t\r+\div_x(\r u)=0,\\[1mm]
\di \pa_t(\r u)+\div_x(\r u\otimes u)+\na_x {p}(\r)=0,\\[2mm]
(\r,u)(0,x)=\left\{
\begin{array}{ll}
\di (\r_-, u_{-}), &\di x_1<0, \\[1mm]
\di (\r_+, u_{+}), &\di x_1>0.
\end{array}
 \right.
 \end{array}
 \right.
\end{equation}
The solution to Riemann problem \eqref{3d-ER} in general contains two nonlinear waves, i.e., shock and rarefaction waves and the above conjecture towards the time-asymptotic stability
of Riemann solutions is well-established in 1D case. In 1960, Illin-Oleinik \cite{IO} first proved the stability of shock and rarefaction wave to 1D scalar Burgers equation. Then Matsumura-Nishihara \cite{MN-85} proved the stability of viscous shock wave to 1D compressible Navier-Stokes system with physical viscosity under the zero mass condition. Independently, Goodman \cite{G} proved the same result to a general system with "artificial" diffusions. Then Liu \cite{L}, Szepessy-Xin \cite{SX} and Liu-Zeng \cite{LZ} removed the crucial zero mass condition in \cite{MN-85, G} by introducing the constant shift on the viscous shock and the diffusion waves and the coupled diffusion waves in the transverse characteristic fields. Masica-Zumbrun \cite{MZ} proved the spectral stability of viscous shock to 1D compressible Navier-Stokes system under a spectral condition, which is slightly weaker than the zero mass condition. Huang-Matsumura \cite{Huang-matsumura} proved the stability of a composite wave
consisting of two viscous shocks for 1D full compressible Navier-Stokes equations
with non-zero initial mass and the condition that the strengths of two viscous shocks are suitably small with same order. Note that all the above results for the stability of shocks are based on the classical anti-derivative techniques, which is essentially suitable to the 1D case and seems not applicable to the multi-dimensional system \eqref{ns} directly.

\

On the other hand, the time-asymptotic stability of rarefaction wave to the 1D compressible Navier-Stokes system was proved by Matsumura-Nishihara \cite{MN-86, MN-92} by using the direct $L^2$-energy methods due to the expanding property of rarefaction wave.
Very recently, Kang-Vasseur-Wang \cite{KVW} proved the stability of the combination wave of viscous shock and rarefaction to 1D compressible Navier-Stokes system by using $a-$contraction methods with the time-dependent shifts to overcome the difficulties caused by the incompatibility of viscous shock and rarefaction.

\

In multi-dimensions, Goodman \cite{G1} first proved the time-asymptotic stability of weak planar viscous shock for the scalar viscous equation by the anti-derivative techniques with the shift function depending on both time and spatially transverse directions, and then Hoff-Zumbrun \cite{HZ, HZ1} extended Goodman's result to the large amplitude shock case. Recently, Kang-Vasseur-Wang \cite{KVW1} proved $L^2$-contraction of large planar viscous shocks up to a  shift function depending on both time and spatial variables.

\

Comparatively speaking, there are very few results on the nonlinearly time-asymptotic stability of planar viscous shocks to the multi-dimensional Navier-Stokes system \eqref{ns} due to the substantial difficulties in the high-dimensional propagation of shocks and the nonlinearities of the  system. In 2017, Humpherys-Lyng-Zumbrun \cite{HLZ} proved the spectral stability of planar viscous Navier-Stokes shocks under the spectral assumptions by the numerical Evans-function methods and one can refer to the survey paper by Zumbrun \cite{Z} for the related results and the references therein.

\

The aim of this paper is to prove the nonlinearly time-asymptotic stability of planar viscous shock wave up to a time-dependent shift for 3D compressible Navier-Stokes system \eqref{ns} by using the weighted energy method under the generic $H^2$-perturbations without the zero mass conditions.

\

The compressibility of viscous shock, which substantially causes the ``bad" sign terms in $L^2$ elementary entropy estimates, is the main difficulty for proving its time-asymptotic stability by energy methods. In 1D case, the classical anti-derivative technique was developed to make full use of the compressibility of viscous shock, and then the zero mass conditions, or generic perturbations with non-zero mass distribution but the constant shift on the viscous shock and the diffusion waves on transverse characteristic fields, are necessarily needed to well define the anti-derivative variables for the perturbation around the viscous shock (\cite{MN-85, G, L, SX, LZ}). However, the above anti-derivative techniques can not be applied directly to prove the stability of planar viscous shocks for the multi-dimensional Navier-Stokes system \eqref{ns}. Alternatively, Kang-Vasseur \cite{KV, KV-1} developed the $a$-contraction method (a kind of weighted $L^2$-relative entropy method) with time-dependent shift to obtain $L^2$-contraction of shock wave to the viscous conservation laws.  One of the advantages of $a$-contraction method is not necessary to introduce the anti-derivative variables for the perturbation and fully use the time-dependent shift in the original perturbation to control the compressibility of shock. The idea can also be applied to prove the stability of planar viscous shock to the multi-dimensional scalar conservation laws (\cite{KVW1}) and the stability of the combination wave of viscous shock and rarefaction to 1D compressible Navier-Stokes system (\cite{KVW}).

\

 Our proof of the time-asymptotic stability for multi-dimensional Navier-Stokes shock is motivated by the $a$-contraction method. However, compared with 1D stability, there are several new difficulties:
 \begin{itemize}
\item [i).] We need to establish a 3D version sharp weighted Poincar${\rm \acute{e}}$ inequality (see Lemma \ref{le-poin}) together with the time-dependent shift $\mathbf{X}(t)$ defined in \eqref{X} to control the compressibility of planar shock.



\

\item [ii).] For the stability analysis of 1D Navier-Stokes system in \cite{KV, KV-1}, Lagrangian structure of the system is fully utilized. However, this structure can not be kept in Eulerian coordinates, especially in multi-dimensional case. Therefore, we need to find a new effective velocity $h:=u-(2\mu+\lambda)\na_{_\xi}v$ (see also \eqref{h}) in 3D Eulerian coordinates and the rewritten system (see \eqref{new-ns}) has the similar stability structure as Lagrangian coordinates. 

\

\item [iii).] Some physical underlying structures of the multi-dimensional Navier-Stokes system \eqref{ns} are used. We use the Hodge decomposition to decompose the diffusive term $\Delta u$ into the irrotational part $\na\div u$ and the rotational part $\na\times\na\times u$ and borrow some ideas from the stability of planar rarefaction wave in \cite{LWW,LW, LWW-CMP, LWW-3} to overcome the wave propagations along the transverse directions and their interactions with the planar viscous shock.
\end{itemize}

\

The rest part of the paper is organized as follows. In section \ref{Preliminaries}, we first list the property of the planar viscous shock and then state our main result. In section \ref{pf-thm}, we first present some useful functional inequalities, and then construct the shift function and give the proof of our main theorem based on the local existence in Proposition \ref{local} and uniform-in-time a-priori estimates in Proposition \ref{priori}.
In section \ref{estimate-(v,h)}, we reformulate the problem in new variable function $(v,h)$ first, and then prove the uniform-in-time $H^2$-estimates in Proposition \ref{priori}.

\section{Planar viscous shock and Main result}\label{Preliminaries}
\setcounter{equation}{0}
In this section, we first describe the planar viscous shock and then state our main result on the time-asymptotic stability of planar viscous shock to 3D compressible Navier-Stokes equations \eqref{ns} under generic $H^2$-perturbations without zero mass conditions.
\subsection{Viscous shock wave}
First we depict planar viscous shock. For definiteness, we consider 2-shock and 1-shock case can be treated similarly. It is well-known that the Riemann problem of 1D
compressible Euler system
\begin{equation}\label{euler}
\left\{
\begin{array}{l}
\di\pa_t\rho+\pa_{x_1}(\rho u_1)=0,\\[1mm]
\di\pa_t(\rho u_1)+\pa_{x_1}(\rho u_1^2+{p}(\r))=0,
\end{array}
\right.
\end{equation}
with Riemann initial data
\begin{equation}\label{R-in}
(\r,u_1)(0,x_1)=(\r_0,u_{10})(x_1)=\left\{
\begin{array}{ll}
\di (\r_-, u_{1-}), &\di x_1<0, \\[1mm]
\di (\r_+, u_{1+}), &\di x_1>0,
\end{array}
 \right.
\end{equation}
determined by the far-field states \eqref{initial}, admits a 2-shock wave solution with the speed $\sigma$
\begin{equation*}
(\r,u_1)(t,x_1)=\left\{
\begin{array}{ll}
\di (\r_-, u_{1-}), &\di x_1<\sigma t, \\[1mm]
\di (\r_+, u_{1+}), &\di x_1>\sigma t,
\end{array}
 \right.
\end{equation*}
provided that the following Rankine-Hugoniot condition
\be\label{RH}
\left\{
\ba
&-\s(\r_+-\r_-)+(\r_+ u_{1+}-\r_-u_{1-})=0,\\
&-\s(\r_+u_{1+}-\r_-u_{1-})+(\r_+u_{1+}^2-\r_-u_{1-}^2)+({p}(\r_+)-{p}(\r_-))=0,
\ea
\right.
\ee
and the Lax entropy condition
$$
\Lambda_2(\r_+, u_{1+})<\sigma<\Lambda_2(\r_-, u_{1-})
$$
with $\Lambda_2(\r, u_{1})=u_1+\sqrt{p'(\r)}$ being the second eigenvalue of Jacobi matrix of Euler system \eqref{euler}, hold true. Denote that
$$
 \xi=(\xi_1,\xi_2, \xi_3), \ \ {\rm with} \ \ \xi_1=x_1-\s t \ \  {\rm and} \ \   \xi_i=x_i,\  i=2,3.
$$
Correspondingly, planar $2-$viscous shock wave $(\r^s,u^s)(\xi_1)$ with $u^s(\xi_1):=(u_1^s(\xi_1),0,0)^t$, connecting $(\r_-,u_{-})$ and $(\r_+,u_{+})$, to the 3D compressible Navier-Stokes equations satisfies the ODE
\be\label{stationary}
\left\{
\begin{array}{l}
\di -\s(\r^s)'+(\r^s u^s_1)'=0, \qquad\qquad\qquad\qquad\qquad\quad ':=\f{d}{d\xi_1}, \\[2mm]
\di -\s(\r^s u_1^s)'+(\r^s (u_1^s)^2)'+p(\r^s)'=(2\mu+\lambda)(u^s_1)'',
\end{array}
\right.
\ee
for the far-field conditions:
\be\label{s-b}
  (\r^s, u^s_1)(-\infty)=(\r_-, u_{1-}),\quad (\r^s, u^s_1)(+\infty)=(\r_+, u_{1+}).
\ee

In the new variable $(t, \xi)$, we can rewrite the system \eqref{ns} as
\be\label{ns-re}
\left\{
\begin{array}{l}
\di \pa_{t}\r-\s\pa_{_{\xi_1}}\r+\div_{_\xi}(\r u)=0,\\[1mm]
\di \pa_{t}(\r u)-\s\pa_{_{\xi_1}}(\r u)+\div_{_\xi}(\r u\otimes u)+\na_{_\xi} {p}(\r)
=\mu \Delta_{_\xi} u+(\mu+\lambda)\na_{_\xi} \div_{_\xi} u.
\end{array}
\right.
\ee
If we introduce the volume function $v=1/\r$, then we can further rewrite the system \eqref{ns-re} as
\begin{equation}\label{ns-xi}
\left\{
\begin{array}{l}
\di \r(\pa_t v-\s\pa_{_{\xi_1}}v+u\cdot\na_{_\xi} v)=\div_{_\xi} u,\\[1mm]
\di \r(\pa_t u-\s\pa_{_{\xi_1}}u+u\cdot\na_{_\xi} u)+\na_{_\xi}{p}(v)
=(2\mu+\lambda)\na_{_\xi} \div_{_\xi} u-\mu\na_{_\xi}\times\na_{_\xi}\times u,
\end{array}
\right.
\end{equation}
where the pressure is defined by $p(v)=bv^{-\gamma}$ and we have used the identity
$$
\Delta_{_\xi}u=\na_{_\xi} \div_{_\xi} u-\na_{_\xi}\times\na_{_\xi}\times u
$$ for the viscosity term. Without loss of generality, we can normalize the constant $b=1$ in pressure in the sequel.

Similarly, by using the volume function $v^s:=1/\r^s$, the ODE system \eqref{stationary} is transformed  into
\be\label{stationary-1}
\left\{
\begin{array}{l}
\di \r^s(-\s(v^s)'+u_1^s(v^s)')=(u_1^s)', \\[1mm]
\di \r^s(-\s(u_1^s)'+u_1^s(u_1^s)')+p(v^s)'=(2\mu+\lambda)(u^s_1)'',
\end{array}
\right.
\ee
where $p(v^s)= (v^s)^{-\gamma}$. Integrating \eqref{stationary}$_1$ from $(-\infty,\xi_1)$, it holds
\be\label{bar-sigma}
-\s\r^s+\r^s u_1^s=-\s\r_-+\r_-u_{1-}=:-\s_* .
\ee
Therefore, the system \eqref{stationary-1} and far-field conditions \eqref{s-b} can be rewritten as
\be\label{stationary-2}
\left\{
\begin{array}{l}
\di -\s_* (v^s)'=(u_1^s)',  \\[1mm]
\di -\s_* (u_1^s)'+p(v^s)'=(2\mu+\lambda)(u^s_1)'',
\end{array}
\right.
\ee
and
\be\label{far-v}
  (v^s, u^s_1)(-\infty)=(v_-, u_{1-}),\quad (v^s, u^s_1)(+\infty)=(v_+, u_{1+}), \quad v_\pm=1/\r_{\pm}.
\ee
By \eqref{stationary-2} (or \eqref{RH}) and \eqref{bar-sigma}, it holds that
$$
\left\{
\begin{array}{l}
\di -\s_* (v_+-v_-)=u_{1+}-u_{1-}, \\[1mm]
\di -\s_* (u_{1+}-u_{1-})+p(v_+)-p(v_-)=0.
\end{array}
\right.
$$
Therefore, we have $\s_* =\sqrt{-\f{p(v_+)-p(v_-)}{v_+-v_-}}>0$ for $2$-shock.

The existence and properties of the 2$-$viscous shock wave $(v^s, u_1^s)(\xi_1)$ can be summarized in the following lemma, while its proof can be found in \cite{MN-85}.

\begin{lemma}\label{le-shock}
For any right state $(v_+,u_{1+})$, there exists a constant $C>0$ such that the following holds. For any left end state $(v_-,u_{1-})\in S_2(v_+,u_{1+})$ and let $\d$ denotes the shock wave strength $\d:=|p(v_-)-p(v_+)|\sim |v_+-v_-|\sim |u_{1-}-u_{1+}|$, there exists a unique (up to a constant shift) solution $(v^s,u^s_1)(\xi_1)$ to ODE system \eqref{stationary-2}, \eqref{far-v} and moreover, it holds
$$
v_{_{\xi_1}}^s>0,\quad u^s_{\scriptscriptstyle 1\xi_1}=-\s_*  v_{_{\xi_1}}^s<0,
$$
and
$$
\ba
&|v^s(\xi_1)-v_-|\leq C\d e^{-C\d|\xi_1|},\qquad\qquad\quad\,\forall\xi_1<0,\\[1mm]
&|v^s(\xi_1)-v_+|\leq C\d e^{-C\d|\xi_1|},\qquad\qquad\quad\,\forall \xi_1>0,\\[1mm]
&|(v^s_{_{\xi_1}},u^s_{_{1\xi_1}})|\leq C\d^2 e^{-C\d|\xi_1|}, \qquad\qquad\quad\,\,\,
\forall\xi_1\in\mathbb{R},\\[1mm]
&|(v^s_{_{\xi_1\xi_1}},u^s_{_{1\xi_1\xi_1}})|\leq C\d|(v^s_{_{\xi_1}},u^s_{_{1\xi_1}})|, \qquad\quad \forall\xi_1\in\mathbb{R},\\[1mm]
&|(v^s_{_{\xi_1\xi_1\xi_1}},u^s_{_{1\xi_1\xi_1\xi_1}})|\leq C\d^2|(v^s_{_{\xi_1}},u^s_{_{1\xi_1}})|, \quad\,\forall\xi_1\in\mathbb{R}.
\ea
$$
\end{lemma}

\subsection{Main result}

%

Now we can state the main result as follows.

\begin{theorem}\label{thm}
Let $(v^s, u^s )(x_1-\s t)$ be the planar $2-$viscous shock wave with $u^s(x_1-\s t):=(u_1^s(x_1-\s t),0,0)^t$, where $(v^s,u_1^s)(x_1-\s t)$ is defined in Lemma \ref{le-shock} with end states $(v_\pm,u_{1\pm})$.  Then there exist positive constants $\d_0$,
$\varepsilon_0$  such that if the shock wave strength  $\d\leq \d_0$, and the initial data $(v_0,u_0)$ satisfies
\be\label{initial pertub}
\|(v_0(x)-v^s(x_1),u_0(x)-u^s(x_1))\|_{H^2(\mathbb{R}\times\mathbb{T}^2)}\leq \varepsilon_0,
\ee
then 3D compressible Navier-Stokes equations \eqref{ns} or \eqref{ns-xi} admits a unique global-in-time solution
$(v,u)(t,x)$ with $\rho=\frac 1v$, and there exists an absolutely continuous shift $\mathbf{X}(t)$ such that
\be\label{space}
\ba
&v(t,x)-v^s(x_1-\s t-\mathbf{X}(t))\in C([0,+\infty);H^2(\mathbb{R}\times\mathbb{T}^2)),\\
&u(t,x)-u^s(x_1-\s t-\mathbf{X}(t))\in C([0,+\infty);H^2(\mathbb{R}\times\mathbb{T}^2)),\\
&\na_{x}\big(v(t,x)-v^s(x_1-\s t-\mathbf{X}(t))\big)\in L^2(0,+\infty;H^1(\mathbb{R}\times\mathbb{T}^2)),\\
&\na_{x}\big(u(t,x)-u^s(x_1-\s t-\mathbf{X}(t))\big)\in L^2(0,+\infty;H^2(\mathbb{R}\times\mathbb{T}^2)).
\ea
\ee
Furthermore, the planar $2-$viscous shock wave $(v^s, u^s)(x_1-\s t)$ is time-asymptotically stable with the time-dependent shift $\mathbf{X}(t)$:
\begin{equation}\label{large-time}
\lim_{t\rightarrow\infty}\sup_{x\in\mathbb{R}\times\mathbb{T}^2}|(v,u)(t,x)-(v^s,u^s)(x_1-\s t-\mathbf{X}(t))|=0,
\end{equation}
and the shift function $\mathbf{X}(t)$ satisfies the time-asymptotic behavior
\be\label{large-X}
\lim_{t\rightarrow\infty} |\dot{\mathbf{X}}(t)|=0.
\ee
\end{theorem}

\begin{remark}
Theorem \ref{thm} states that if the two far-fields states $(\rho_\pm, u_\pm)$ in \eqref{initial} are connected by the shock wave, then the solution to 3D compressible Navier-Stokes equations \eqref{ns}, or equivalently \eqref{ns-xi}, tends to the corresponding planar viscous shock with the time-dependent shift $\mathbf{X}(t)$ under the generic $H^2$-perturbations, in particular, without zero mass conditions.
\end{remark}

\begin{remark}
The shift function $\mathbf{X}(t)$ is proved to satisfy the time-asymptotic behavior \eqref{large-X}, which implies
$$
\lim_{t\rightarrow+\infty}\frac{\mathbf{X}(t)}{t}=0,
$$
that is, the shift function $\mathbf{X}(t)$ grows at most sub-linearly with respect to the time $t$, and therefore, the shifted planar viscous shock wave $(v^s,u^s)(x_1-\s t-\mathbf{X}(t))$ keeps the original traveling wave profile time-asymptotically.
\end{remark}

\begin{remark}
Theorem \ref{thm} is the first analytical result on the time-asymptotic stability of planar viscous shock wave to the multi-dimensional system \eqref{ns} with physical viscosities as far as we know. Moreover, our stability result is unconditional for the weak planar Navier-Stokes shock in 3D case.
\end{remark}


{\textbf{Notation.}} Throughout this paper, several positive generic constants are denoted by $C$ if without confusions. We define
$$
x':=(x_2,x_3),\quad dx':=dx_2dx_3, \quad{\rm and }\quad \xi':=(\xi_2,\xi_3),\quad d\xi:=d\xi_2d\xi_3.
$$
For $1\leq r\leq \infty$, we denote $L^r:=L^r(\Omega)=L^r(\mathbb{R}\times\mathbb{T}^2)$, and use the notation $\|\cdot\|:=\|\cdot\|_{L^2}$.  For a non-negative integer $m$, the space $H^m(\Omega)$ denotes the $m-$th order Sobolev space over $\Omega$ in the $L^2$ sense with the norm
$$
\|f\|_{H^m}:=\Big(\sum_{l=0}^m\|\na^l f\|^2\Big)^{1/2},\quad
\|f\|:=\Big(\int_{\Omega}|f|^2d\xi\Big)^{1/2}=\Big(\int_{\mathbb{T}^2}\int_{\mathbb{R}}|f|^2d\xi_1d\xi'\Big)^{1/2}.
$$
Also, we denote
$$
\|(f,g)\|_{H^m}=\|f\|_{H^m}+\|g\|_{H^m}.
$$

%
%
%
%


%
%
%
%

\section{Proof of main result}\label{pf-thm}
\setcounter{equation}{0}

\subsection{Some functional inequalities}\label{functional-inequality}
We first present a 3D weighted sharp Poincar${\rm \acute{e}}$ type inequality, which is a 3D version of 1D weighted Poincar${\rm \acute{e}}$ inequality in \cite{KV-1} and plays a very important role in our stability analysis.
\begin{lemma}\label{le-poin}
For any $f:[0,1]\times\mathbb{T}^2\rightarrow\mathbb{R}$ satisfying $$\int_{\mathbb{T}^2}\int_0^1
\big[y_1(1-y_1)|\pa_{y_1}f|^2+\f{|\na_{y'}f|^2}{y_1(1-y_1)}\big]dy_1dy'<\infty,$$ it holds
\be\label{poin-2}
\int_{\mathbb{T}^2}\int_0^1|f-\bar f|^2dy_1dy'\leq \f12\int_{\mathbb{T}^2}\int_0^1
y_1(1-y_1)|\pa_{y_1}f|^2dy_1dy'
+\f{1}{16\pi^2}\int_{\mathbb{T}^2}\int_0^1\f{|\na_{y'}f|^2}{y_1(1-y_1)}dy_1dy'.
\ee
where $\bar f=\int_{\mathbb{T}^2}\int_0^1 f dy_1dy'$ and $y'=(y_2,y_3)$.
\end{lemma}

\textbf{\emph{Proof}}: The proof is motivated by that of 1D weighted Poincar${\rm \acute{e}}$ inequality  in \cite{KV-1} and here we need to concern the transverse directions $(y_2, y_3)\in \mathbb{T}^2$ additionally. Let $\{P_n:[-1,1]\rightarrow\mathbb{R}\}_{n\in \mathbb{Z}, n\geq0}$ be an orthogonal basis of Legendre's polynomials, that
are the solutions to Legendre's differential equations:
\be\label{Le-equation}
\f{d}{dx_1}\Big((1-x_1^2)\f{d}{dx_1}P_n(x_1)\Big)=-n(n+1)P_n(x_1),
\ee
and satisfy the orthogonality in $L^2[-1,1]$, i.e., $\int_{-1}^1 P_i P_j=\d_{ij}$. Then for fixed $x'=(x_2,x_3)\in \mathbb{T}^2$ and any $w(\cdot,x')\in L^2[-1,1]$,
we have
$$
w(x_1,x')=\sum_{i=0}^{\infty}c_i(x')P_i(x_1),\quad
c_i(x')=\int_{-1}^1w(x_1,x')P_i(x_1)dx_1.
$$
In particular, we see that $P_0(x_1)=\f{1}{\sqrt{2}}=:P_0$, thus
$$
\int_{\mathbb{T}^2}c_0(x')P_0(x_1)dx'=P_0\int_{\mathbb{T}^2}c_0(x')dx'=\f12\int_{\mathbb{T}^2}\int_{-1}^1w(x_1,x')dx_1dx'
=:\bar w.
$$
Set $\bar c_0:=\int_{\mathbb{T}^2}c_0(x')dx'$, it holds
$$
\ba
w(x_1,x')-\bar w&=\sum_{i=0}^{\infty}c_i(x')P_i(x_1)-\int_{\mathbb{T}^2}c_0(x')P_0dx'\\
&=\sum_{i=1}^{\infty}c_i(x')P_i(x_1)+(c_0(x')-\bar c_0)P_0.
\ea
$$
Then we have
$$
\ba
\int_{\mathbb{T}^2}\int_{-1}^1|w-\bar w|^2dx_1dx'=\sum_{i=1}^{\infty}\int_{\mathbb{T}^2}\int_{-1}^1c_i^2(x')P_i^2(x_1)dx_1dx'
+\int_{\mathbb{T}^2}|c_0-\bar c_0|^2dx'=:A_1+A_2.
\ea
$$
By Legendre's differential equations \eqref{Le-equation}, we have
$$
\ba
&\int_{\mathbb{T}^2}\int_{-1}^1(1-x_1)^2|\pa_{x_1}w|^2dx_1dx'\\
=&-\int_{\mathbb{T}^2}\int_{-1}^1\pa_{x_1}\Big((1-x_1^2)\pa_{x_1}w\Big)wdx_1dx'
=-\int_{\mathbb{T}^2}\int_{-1}^1\pa_{x_1}\Big((1-x_1^2)\pa_{x_1}w\Big)(w-\bar w)dx_1dx'\\
=&-\int_{\mathbb{T}^2}\int_{-1}^1\pa_{x_1}\Big(\sum_{i=1}^{\infty}c_i(x')(1-x_1^2)\f{d}{dx_1}P_i(x_1)\Big)\sum_{j=1}^{\infty}c_j(x')P_j(x_1)dx_1dx'\\
&-\int_{\mathbb{T}^2}\int_{-1}^1\pa_{x_1}\Big(\sum_{i=1}^{\infty}c_i(x')(1-x_1^2)\f{d}{dx_1}P_i(x_1)\Big)(c_0(x')-\bar c_0)P_0dx_1dx'\\
=&\sum_{i=1}^{\infty}\int_{\mathbb{T}^2}\int_{-1}^1i(i+1)c_i^2(x')P_i^2(x_1)dx_1dx'
\geq2\sum_{i=1}^{\infty}\int_{\mathbb{T}^2}\int_{-1}^1c_i^2(x')P_i^2(x_1)dx_1dx'
=2A_1,
\ea
$$
which implies that
\begin{equation*}
A_1\leq \frac12 \int_{\mathbb{T}^2}\int_{-1}^1(1-x_1)^2|\pa_{x_1}w|^2dx_1dx'.
\end{equation*}

For $A_2,$ by using Poincar${\rm \acute{e}}$ inequality, it holds
$$
\ba
A_2&\leq \f{1}{(2\pi)^2}\int_{\mathbb{T}^2}|\na_{x'}c_0|^2dx'
=\f{1}{4\pi^2}\int_{\mathbb{T}^2}\Big(\big(\int_{-1}^1\pa_{x_2}wP_0dx_1\big)^2+\big(\int_{-1}^1\pa_{x_3}wP_0dx_1\big)^2\Big)dx'\\
&=\f{1}{8\pi^2}\int_{\mathbb{T}^2}\Big(\big(\int_{-1}^1\pa_{x_2}wdx_1\big)^2+\big(\int_{-1}^1\pa_{x_3}wdx_1\big)^2\Big)dx'\\
&\leq\f{1}{4\pi^2}\int_{\mathbb{T}^2}\int_{-1}^1(|\pa_{x_2}w|^2+|\pa_{x_3}w|^2)dx_1dx'
=\f{1}{4\pi^2}\int_{\mathbb{T}^2}\int_{-1}^1|\na_{x'}w|^2dx_1dx'.
\ea
$$
Combining the above two estimates on $A_1$ and $A_2$, we have
$$
\ba
\int_{\mathbb{T}^2}\int_{-1}^1|w-\bar w|^2dx_1dx'\leq \f12\int_{\mathbb{T}^2}\int_{-1}^1(1-x_1)^2|\pa_{x_1}w|^2dx_1dx'
+\f{1}{4\pi^2}\int_{\mathbb{T}^2}\int_{-1}^1|\na_{x'}w|^2dx_1dx'.
\ea
$$
By a change of variable $y_1:=\frac{x_1+1}{2}, y':=x'$ and $W(y_1,y'):=w(2y_1-1,y')=w(x_1,x')$, we have
$$
\ba
\int_{\mathbb{T}^2}\int_{0}^1|W-\bar W|^2dy_1dy'\leq \f12\int_{\mathbb{T}^2}\int_{0}^1y_1(1-y_1)|\pa_{y_1}W|^2dy_1dy'
+\f{1}{4\pi^2}\int_{\mathbb{T}^2}\int_{0}^1|\na_{y'}W|^2dy_1dy',
\ea
$$
where $\bar W:=\int_{\mathbb{T}^2}\int_0^1 W dy_1dy'$. Notice that
$0\leq y_1(1-y_1)\leq \f14$ for $y_1\in[0,1]$, and so
$$
\frac{1}{y_1(1-y_1)}\geq 4,
$$
then we can prove \eqref{poin-2} .
\hfill $\Box$

\

Now we list a 3D Gagliardo-Nirenberg inequality in the domain $\Omega=\mathbb{R}\times \mathbb{T}^2$, whose proof can be
found in \cite{LWW-FP, WW}.
\begin{lemma}\label{sobolev}
It holds for $g(x)\in H^2(\Omega)$ with $x=(x_1,x_2,x_3)\in \Omega:=\mathbb{R}\times \mathbb{T}^2$ that
\begin{equation}\label{sobolev-inequality}
\|g\|_{L^{\infty}(\Omega)}\leq \sqrt{2}\|g\|^{\f12}_{L^2(\Omega)}\|\pa_{x_1}g\|^{\f12}_{L^2(\Omega)}
+ C\|\na_xg\|^{\f12}_{L^2(\Omega)}\|\na_x^2 g\|^{\f12}_{L^2(\Omega)},
\end{equation}
where $C>0$ is a positive constant.
\end{lemma}

Then we list several estimates on the relative quantities. For any function $F$ defined on $\mathbb{R}_+$, we define the associated relative quantity for $v$, $w\in\mathbb{R}_+$ as
$$
F(v|w)=F(v)-F(w)-F'(w)(v-w).
$$
We gather, in the following lemma, some useful inequalities on the relative quantities associated to the pressure $p(v)=v^{-\gamma}$ and the internal energy
$Q(v)=\f{v^{-\gamma+1}}{\gamma-1}$. The proofs are based on Taylor expansions and can be found in \cite{KV-1}.

\begin{lemma}\label{inequality-Q}
For given constants $\gamma>1$, and $v_->0$, their exist constants $C$, $\d_*>0$,  such that the following holds true.
\begin{itemize}

\item[(1)] For any $v$, $w$ such that $0<w<2v_-$, $0<v\leq 3v_-$,
\be\label{use-4}
|v-w|^2\leq C Q(v|w),\qquad |v-w|^2\leq C p(v|w).
\ee
\item[(2)] For any $v$, $w>\f{v_-}{2}$,
\be
|p(v)-p(w)|\leq C|v-w|.
\ee
\item[(3)] For any $0<\d<\d_*$, and for any $(v,w)\in\mathbb{R}^2_+$ satisfying $|p(v)-p(w)|<\d$, and $|p(w)-p(v_-)|<\d$, the following holds true:
\be\label{use-1}
p(v|w)\leq \Big(\f{\gamma+1}{2\gamma}\f{1}{p(w)}+C\d\Big)|p(v)-p(w)|^2,
\ee
\be\label{use-2}
Q(v|w)\geq \f{|p(v)-p(w)|^2}{2\gamma p^{1+\f{1}{\gamma}}(w)}
-\f{1+\gamma}{3\gamma^2}\f{(p(v)-p(w))^3}{p^{2+\f{1}{\gamma}}(w)},
\ee
\be\label{use-3}
Q(v|w)\leq\Big(\f{1}{2\gamma p^{1+\f{1}{\gamma}}(w)}+C\d\Big)|p(v)-p(w)|^2.
\ee

\end{itemize}
\end{lemma}

Finally we give an estimate involving the inverse of the pressure function $p(v)=v^{-\gamma}$, while its proof can be found in \cite{KV-1}.
\begin{lemma}\label{inver-pressure}
Fix $v_->0$. Then there exist $\d_0>0$ and $C>0$ such that for any $v_+>0$,
such that $0<\d:=p(v_-)-p(v_+)\leq \d_0$, $v_-\leq v\leq v_+$, we have
$$
\left|\f{v-v_-}{p(v)-p(v_-)}+\f{v-v_+}{p(v_+)-p(v)}+\f{1}{2}\f{p''(v_-)}{p'(v_-)^2}(v_--v_+)\right|\leq C\d^2.
$$
\end{lemma}

\

\subsection{Construction of shift function $\mathbf{X}(t)$}\label{sec-shift}

For notational simplification, denote
$$
\big((v^s)^{-\mathbf{X}},(u^s)^{-\mathbf{X}}\big)(\xi_1):=\big(v^s(\xi_1-\mathbf{X}(t)),u^s(\xi_1-\mathbf{X}(t))\big) 
$$
with the shift function $\mathbf{X}(t)$ to be defined in \eqref{X}.

The definition of the shift function $\mb{X}(t)$ depends on the weight function
$a:\mathbb{R}\rightarrow\mathbb{R}$ defined in \eqref{a}. For now we will only assume the fact that
$\|a\|_{C^1(\mathbb{R})}\leq 2$. Then we can define the shift $\mathbf{X}(t)$ as a solution to the ODE:
\be\label{X}
\left\{
\ba
\dot{\mathbf{X}}(t)&=-\f{M}{\d}\Big[
\int_{\mathbb{T}^2}\int_{\mathbb{R}}\f{a^{-\mathbf{X}}(\xi_1)}{\s_* }\r(h_1^s)^{-\mathbf{X}}_{_{\xi_1}}\big(p(v)-p((v^s)^{-\mathbf{X}})\big)d\xi_1d\xi'\\
&\qquad-\int_{\mathbb{T}^2}\int_{\mathbb{R}}a^{-\mathbf{X}}(\xi_1)\r p'((v^s)^{-\mathbf{X}})(v-(v^s)^{-\mathbf{X}})(v^s)^{-\mathbf{X}}_{_{\xi_1}}d\xi_1d\xi'\Big],\\
\mathbf{X}(0)&=0,
\ea
\right.
\ee
where the function $h^s_1:=u^s_1-(2\mu+\lambda)\pa_{_{\xi_1}}v^s$ as defined in \eqref{h-s} and the constant $M:=\f54\f{\gamma+1}{2\gamma}\f{\s_-^3v_-^2}{p(v_-)}$ with $\s_-=\sqrt{-p'(v_-)}$.





Let $F(t,\mathbf{X}(t))$ be the right-hand side of the ODE $\eqref{X}_1$.
Thanking to the facts that $\|a\|_{C^1(\mathbb{R})}\leq 2$, $\|v^s\|_{C^2(\mathbb{R})}\leq v_+$,
and $\|v^s_{\xi_1}\|_{L^1(\mathbb{R})}\leq C\d$, we can find some constant $C>0$, such that
\be\label{F(X)}
\sup_{\mathbf{X}\in\mathbb{R}}|F(t,\mathbf{X})|\leq\f{C}{\d}\|a\|_{C^1}
(\|v\|_{L^{\infty}}+\|(v^s)^{-\mb{X}}\|_{L^{\infty}})\int_{\mathbb{T}^2}\int_{\mathbb{R}}|(v^s)^{-\mathbf{X}}_{_{\xi_1}}|d\xi_1d\xi'\leq C,
\ee
and
\be
\sup_{\mathbf{X}\in\mathbb{R}}|\pa_{\mathbf{X}}F(t,\mathbf{X})|\leq \f{C}{\d}\|a\|_{C^1}
(\|v\|_{L^{\infty}}+\|(v^s)^{-\mb{X}}\|_{L^{\infty}})\\
\int_{\mathbb{T}^2}\int_{\mathbb{R}}|(v^s)^{-\mathbf{X}}_{_{\xi_1}}|d\xi_1d\xi'\leq C.
\ee
Then ODE \eqref{X} has a unique absolutely continuous solution $\mathbf{X}(t)$ defined on any interval in time $[0,T]$ by the well-known Cauchy-Lipschitz theorem.
In particular, since $|\dot{\mathbf{X}}(t)|\leq C$ by \eqref{F(X)}, we can obtain
$$
|\mathbf{X}(t)|\leq Ct, \qquad \forall t\in[0, T].
$$

\hfill $\Box$

\subsection{Proof of Theorem \ref{thm}}\label{pf-main}

In order to prove Theorem \ref{thm}, we shall combine a local existence result together with a-priori estimates by continuation arguments.

\begin{proposition}\label{local}
{\rm (Local existence)} Let $(v^s, u^s )(x_1-\s t)$ be the planar 2-viscous shock wave with $u^s(x_1-\s t):=(u_1^s(x_1-\s t),0,0)^t$.
For any $\Xi>0$, suppose the initial data $(v_0,u_0)$ satisfies
$$
\|(v_0(x)-v^s(x_1),u_0(x)-u^s(x_1))\|_{H^2(\mathbb{R}\times\mathbb{T}^2)}\leq \Xi.
$$
Then there exists a positive constant $T_0$ depending on $\Xi$ such that the 3D compressible Navier-Stokes system \eqref{ns-xi} has a unique solution $(v,u)$ on $(0,T_0)$ satisfying
$$
\ba
&v-v^s\in C([0,T_0]; H^2(\mathbb{R}\times\mathbb{T}^2)),\quad \na_{_\xi}(v-v^s)\in L^2(0,T_0;H^1(\mathbb{R}\times\mathbb{T}^2)),\\
&u-u^s\in C([0,T_0]; H^2(\mathbb{R}\times\mathbb{T}^2)),\quad \na_{_\xi}(u-u^s)\in L^2(0,T_0;H^2(\mathbb{R}\times\mathbb{T}^2)),
\ea
$$
and for $t\in[0,T_0]$, it holds
\begin{equation}\label{local-es}
\ba
&\sup_{\tau\in[0,t]}\|(v-v^s,u-u^s)(\tau)\|^2_{H^2}
+\int_0^t\big(\|\nabla_{_{\xi}}(v-v^s)\|^2_{H^1}+\|\nabla_{_{\xi}}(u-u^s)\|_{H^2}^2\big)d\tau\\
&\leq 4\|(v_0-v^s,u_0-u^s)\|^2_{H^2}.
\ea
\end{equation}
\end{proposition}

\begin{proposition}\label{priori}
{\rm (A priori estimates)} Suppose that $(v,u)$ is the solution to \eqref{ns-xi} on $[0,T]$ for some $T>0$,
and $((v^s)^{-\mathbf{X}},(u^s)^{-\mathbf{X}})$ is the solution of \eqref{stationary-1} with the shift function $\mathbf{X}$, which is
an absolutely continuous solution to \eqref{X}. Then there exist positive constants $\d_0\leqq1$, $\chi_0\leqq1$
and $C_0$ independent of $T$, such that if the shock wave strength $\d<\d_0$ and
$$
\ba
&v-(v^s)^{-\mathbf{X}}\in C([0,T]; H^2(\mathbb{R}\times\mathbb{T}^2)),\quad \na_{_\xi}(v-(v^s)^{-\mathbf{X}})\in L^2(0,T;H^1(\mathbb{R}\times\mathbb{T}^2)),\\
&u-(u^s)^{-\mathbf{X}}\in C([0,T]; H^2(\mathbb{R}\times\mathbb{T}^2)),\quad \na_{_\xi}(u-(u^s)^{-\mathbf{X}})\in L^2(0,T;H^2(\mathbb{R}\times\mathbb{T}^2)),
\ea
$$
with
\begin{equation}\label{assumption}
\ba
\chi :=\sup_{0\leq t \leq T}\|(v-(v^s)^{-\mathbf{X}},u-(u^s)^{-\mathbf{X}})(t,\cdot)\|_{H^2}\leq \chi_0,
\ea
\end{equation}
then the following estimate holds,
\begin{equation}\label{full-es}
\ba
&\sup_{0\leq t\leq T}\|(v-(v^s)^{-\mathbf{X}},u-(u^s)^{-\mathbf{X}})(t,\cdot)\|^2_{H^2}
+\d\int_0^T|\dot{\mathbf{X}}(t)|^2dt\\
&+\int_0^T\Big(\|\sqrt{|(v^s)_{\xi_1}^{-\mathbf{X}}|}(v-(v^s)^{-\mathbf{X}})\|^2+\|\na_{_\xi}(v-(v^s)^{-\mathbf{X}})\|^2_{H^1}
+\|\na_{_\xi}(u-(u^s)^{-\mathbf{X}})\|_{H^2}^2\Big)dt\\
&\leq C_0\|(v_0-v^s,u_0-u^s)\|^2_{H^2}.
\ea
\end{equation}
In addition, by \eqref{X}, we have
\be\label{point-X}
|\dot{ \mathbf{X}}(t)|\leq C_0\|v-(v^s)^{-\mathbf{X}}(t,\cdot)\|_{L^{\infty}}, \qquad \forall t\leq T.
\ee
\end{proposition}

Proposition \ref{local} can be proved by a standard way, see \cite{Solo}, we omit it. Proposition \ref{priori} will be proved in section \ref{estimate-(v,h)}.
Here we show Theorem \ref{thm} by the continuation arguments based on Propositions \ref{local} and \ref{priori}.

{\textbf{\emph{Proof of Theorem \ref{thm}}}:}
We first prove \eqref{space} in Theorem \ref{thm} by the continuation method based on Propositions \ref{local} and \ref{priori}.
Considering the maximal existence time of the solution
\be\label{T-max}
T_{\max}:=\Big\{t>0\Big|\sup_{\tau\in[0,t]}\|(v-(v^s)^{-\mathbf{X}},u-(u^s)^{-\mathbf{X}})(\tau)\|_{H^2}\leq\chi_0\Big\}.
\ee
We shall show the maximal existence time $T_{\max}=+\infty$ by the following steps. We define
$$
\varepsilon_0=\min\Big\{\f{\chi_0}{4},\f{\chi_0}{8\sqrt{C_0}}\Big\},\quad
\Xi=\f{\chi_0}{4},
$$
where $\chi_0$ and $C_0$ are given in Proposition \ref{priori}.

{\textbf{\emph{Step 1}}:} Suppose $\|(v_0-v^s,u_0-u^s)\|_{H^2}\leq \varepsilon_0\leq \f{\chi_0}{4}(=\Xi)$,  by local existence result in Proposition \ref{local}, there is a positive constant $T_0=T_0(\Xi)$ such that a unique solution exists on $[0,T_0]$ and satisfies $\|(v-v^s,u-u^s)(t)\|_{H^2}\leq 2\|(v_0-v^s,u_0-u^s)\|_{H^2}\leq 2\Xi=\f{\chi_0}{2}$ for $t\in[0,T_0]$. Without loss of generality, we can assume $T_0\leq1$. Then Sobolev inequality implies $\|(v-v^s)(t)\|_{L^{\infty}}\leq C\chi_0$ for $t\in[0,T_0]$. Using $v_-<v^s(\xi_1)<v_+$ and the smallness of $\chi_0$ in Proposition \ref{priori}, it holds $\f{v_-}{2}<v(t,\xi)<2v_+$ for $(t,\xi)\in[0,T_0]\times\Omega$. Therefore, we can see that \eqref{F(X)} holds for $t\in[0,T_0]$, and we can deduce from \eqref{X} that $|\mathbf{X}(t)|\leq Ct$ for $t\in[0,T_0]$. Then by the mean value theorem, it holds
$$
\|(v^s-(v^s)^{-\mathbf{X}},u^s-(u^s)^{-\mathbf{X}})(t)\|_{H^2}=|\mathbf{X}(t)|\,\|(v^s_{_{\xi_1}},u^s_{_{\xi_1}})\|_{H^2}
\leq C\d^{3/2}t\leq C\d_0\leq \f{\chi_0}{8}
$$
for suitably small $\d_0$. Therefore, it holds for $t\in[0,T_0]$ that
$$
\ba
\|(v-(v^s)^{-\mathbf{X}},u-(u^s)^{-\mathbf{X}})(t)\|_{H^2}
&\leq
\|(v-v^s,u-u^s)(t)\|_{H^2}+\|(v^s-(v^s)^{-\mathbf{X}},u^s-(u^s)^{-\mathbf{X}})(t)\|_{H^2}\\
&\leq\f{\chi_0}{2}+\f{\chi_0}{8}<\chi_0.
\ea
$$
Hence, we can apply the a-priori estimates in Proposition \ref{priori} with $T=T_0$ and get the estimate
$$
\|(v-(v^s)^{-\mathbf{X}},u-(u^s)^{-\mathbf{X}})(t)\|_{H^2}\leq
\sqrt{C_0}\|(v_0-v^s,u_0-u^s)\|_{H^2}\leq \sqrt{C_0}\varepsilon_0\leq\f{\chi_0}{8}
$$
for $t\in[0,T_0]$.

{\textbf{\emph{Step 2}}:} If the maximal existence time $T_{\max}<+\infty$, then there is a positive integer $N\geq1$, which may depend on $\chi_0$, such that $T_{\max}\in((N-1)T_0,NT_0]$. We can choose the small constant $\d_0$ satisfying $\sqrt{\d_0}\leq \f{1}{N+1}$.
We know from Step 1 that
$$
\ba
&\|(v-v^s,u-u^s)(T_0)\|_{H^2}\\
&\leq\|(v-(v^s)^{-\mathbf{X}},u-(u^s)^{-\mathbf{X}})(T_0)\|_{H^2}
+\|(v^s-(v^s)^{-\mathbf{X}},u^s-(u^s)^{-\mathbf{X}})(T_0)\|_{H^2}\\
&\leq \f{\chi_0}{8}+\f{\chi_0}{8}=\f{\chi_0}{4}(=\Xi).
\ea
$$
Hence, we can apply local existence result in Proposition \ref{local} by taking $t=T_0$ as the new initial time. Then we have a unique solution on $[T_0,2T_0]$ with the estimate $\|(v-v^s,u-u^s)(t)\|_{H^2}\leq 2\|(v-v^s,u-u^s)(T_0)\|_{H^2} \leq 2\Xi=\f{\chi_0}{2}$ for $t\in[T_0,2T_0]$. This together with Step 1 implies that $\|(v-v^s,u-u^s)(t)\|_{H^2}\leq \f{\chi_0}{2}$ holds on $t\in[0,2T_0]$. Same as Step 1, we can show that $|\mathbf{X}(t)|\leq Ct$ holds on $t\in[0,2T_0]$.
Using the smallness of $\d_0$, we have
$$
\ba
&\|(v^s-(v^s)^{-\mathbf{X}},u^s-(u^s)^{-\mathbf{X}})(t)\|_{H^2}=|\mathbf{X}(t)|\,\|(v^s_{_{\xi_1}},u^s_{_{\xi_1}})\|_{H^2}\\
&\leq C_*\d^{3/2}t\leq 2C_*\d_0 \sqrt{\d_0}\leq C_*\d_0\leq \f{\chi_0}{8}.
\ea
$$
Therefore, it holds for $t\in[0,2T_0]$ that
$$
\ba
\|(v-(v^s)^{-\mathbf{X}},u-(u^s)^{-\mathbf{X}})(t)\|_{H^2}
&\leq
\|(v-v^s,u-u^s)(t)\|_{H^2}+\|(v^s-(v^s)^{-\mathbf{X}},u^s-(u^s)^{-\mathbf{X}})(t)\|_{H^2}\\
&\leq\f{\chi_0}{2}+\f{\chi_0}{8}<\chi_0.
\ea
$$
Hence, we can apply the a-priori estimates in Proposition \ref{priori} again with $T=2T_0$ and get the estimate
$$
\|(v-(v^s)^{-\mathbf{X}},u-(u^s)^{-\mathbf{X}})(t)\|_{H^2}\leq
\sqrt{C_0}\|(v_0-v^s,u_0-u^s)\|_{H^2}\leq \sqrt{C_0}\varepsilon_0\leq\f{\chi_0}{8}
$$
for $t\in[0,2T_0]$.

{\textbf{\emph{Step 3}}:} Thus, repeating this continuation process, we can extend the solution to the interval $[0,NT_0]$ successively. At the time $t=NT_0$, it holds
$$
\ba
&\|(v-v^s,u-u^s)(NT_0)\|_{H^2}\\
&\leq\|(v-(v^s)^{-\mathbf{X}},u-(u^s)^{-\mathbf{X}})(NT_0)\|_{H^2}
+\|(v^s-(v^s)^{-\mathbf{X}},u^s-(u^s)^{-\mathbf{X}})(NT_0)\|_{H^2}\\
&\leq \f{\chi_0}{8}+\f{\chi_0}{8}=\f{\chi_0}{4}(=\Xi).
\ea
$$
Hence, we can apply Proposition \ref{local} by taking $t=NT_0$ as the new initial time. Then we have a unique solution on $[NT_0,(N+1)T_0]$ with the estimate $\|(v-v^s,u-u^s)(t)\|_{H^2}\leq 2\|(v-v^s,u-u^s)(NT_0)\|_{H^2}\leq 2\Xi=\f{\chi_0}{2}$ for $t\in[NT_0,(N+1)T_0]$, which implies that $\|(v-v^s,u-u^s)(t)\|_{H^2}\leq \f{\chi_0}{2}$ holds on $t\in[0,(N+1)T_0]$.
Meanwhile, we can also show that $|\mathbf{X}(t)|\leq Ct$ holds on  $t\in[0,(N+1)T_0]$. Using the smallness of $\d_0$, we have
$$
\ba
&\|(v^s-(v^s)^{-\mathbf{X}},u^s-(u^s)^{-\mathbf{X}})(t)\|_{H^2}=|\mathbf{X}(t)|\,\|(v^s_{_{\xi_1}},u^s_{_{\xi_1}})\|_{H^2}\\
&\leq C_*\d^{3/2}t\leq C_*\d_0 \sqrt{\d_0}(N+1)\leq C_*\d_0 \leq \f{\chi_0}{8}.
\ea
$$
Therefore, it holds for $t\in[0,(N+1)T_0]$ that
$$
\ba
\|(v-(v^s)^{-\mathbf{X}},u-(u^s)^{-\mathbf{X}})(t)\|_{H^2}
&\leq
\|(v-v^s,u-u^s)(t)\|_{H^2}+\|(v^s-(v^s)^{-\mathbf{X}},u^s-(u^s)^{-\mathbf{X}})(t)\|_{H^2}\\
&\leq\f{\chi_0}{2}+\f{\chi_0}{8}<\chi_0.
\ea
$$
Hence, we can apply Proposition \ref{priori} again with $T=(N+1)T_0$ and get the estimate
$$
\|(v-(v^s)^{-\mathbf{X}},u-(u^s)^{-\mathbf{X}})(t)\|_{H^2}\leq
\sqrt{C_0}\|(v_0-v^s,u_0-u^s)\|_{H^2}\leq \sqrt{C_0}\varepsilon_0\leq\f{\chi_0}{8}
$$
for $t\in[0,(N+1)T_0]$. This indicates that the solution has been extended to the interval $[0,(N+1)T_0]$, which contradicts that $T_{\max}(\leq NT_0)$ is the maximum existence time. Therefore, the maximum existence time defined in \eqref{T-max} is infinity, that is, $T_{\max}=+\infty$.

To complete the proof of Theorem \ref{thm}, we remain to justify
the time-asymptotic behaviors \eqref{large-time} and \eqref{large-X}. We use the same notations as in Lemma \ref{le-2nd-v} that $\phi:=v-(v^s)^{-\mathbf{X}}$ and
$\psi:=u-(u^s)^{-\mathbf{X}}$. Set
$$
g(t):=\|\na_{_\xi}\phi(t)\|^2+\|\na_{_\xi}\psi(t)\|^2.
$$
The aim is to show that
\be\label{4.1}
\begin{array}{l}
\di\int_0^{+\infty}\big(|g(t)|+|g'(t)|\big)dt<\infty,
\end{array}
\ee
which implies
\be\label{4.2}
\di \lim_{t\rightarrow+\infty}g(t)=
\lim_{t\rightarrow+\infty}(\|\na_{_\xi}\phi(t)\|^2+\|\na_{_\xi}\psi(t)\|^2)
=0.
\ee
First, we can deduce from \eqref{full-es} that $\di\int_0^{+\infty}|g(t)|dt<\infty$. Then we apply $\na_{_\xi}$ to the equation \eqref{re-perturb-vu}$_1$ to get
$$
\ba
\int_0^{t}\|\na_{_\xi}\pa_{_t}\phi\|^2d\tau&\leq C\int_0^t\|\na_{_\xi}^2(\phi,\psi)\|^2d\tau
+C(\d+\chi)\int_0^t\|\na_{_\xi}(\phi,\psi)\|^2d\tau\\
&+C\d^2\int_0^t\Big(|\dot{\mathbf{X}}(\tau)|^2+G_3(\tau)+G^s(\tau)+D(\tau)\Big)d\tau\leq C.
\ea
$$
Meanwhile, it follows from \eqref{dao-perturb-vu}$_2$ that
$$
\ba
\int_0^t\|\na_{_\xi}\pa_{_t}\psi\|^2d\tau&\leq C\int_0^t(\|\na_{_\xi}^2(\phi,\psi)\|^2+\|\na_{_\xi}^3\psi\|^2)d\tau
+C(\d+\chi)\int_0^t\|\na_{_\xi}(\phi,\psi)\|^2d\tau\\
&+C\d^2\int_0^t\Big(|\dot{\mathbf{X}}(\tau)|^2G_3(\tau)+G^s(\tau)+D(\tau)\Big)d\tau\leq C.
\ea
$$
Using the above two facts and Cauchy inequality, we have
$$
\ba
\int_0^{+\infty}|g'(t)|dt&=\int_0^{+\infty}\int_{\mathbb{T}^2}\int_{\mathbb{R}}(2|\na_{_\xi}\phi||\na_{_\xi}\pa_{_t}\phi|
+2|\na_{_\xi}\psi||\na_{_\xi}\pa_{_t}\psi|)d\xi_1d\xi'dt\\
&\leq 2\int_0^{+\infty}(\|\na_{_\xi}\phi\|\|\na_{_\xi}\pa_{_t}\phi\|+\|\na_{_\xi}\psi\|\|\na_{_\xi}\pa_{_t}\psi\|)dt
< \infty.
\ea
$$
By Gagliardo-Nirenberg inequality in Lemma \ref{sobolev} and \eqref{4.2}, we have
$$
\lim_{t\rightarrow+\infty}\|(\phi,\psi)\|_{L^{\infty}}\leq \lim_{t\rightarrow+\infty}(\sqrt{2}\|(\phi,\psi)\|^{\f12}\|\pa_{_{\xi_1}}(\phi,\psi)\|^{\f12}
+C\|\na_{_{\xi}}(\phi,\psi)\|^{\f12}\|\na^2_{_{\xi}}(\phi,\psi)\|^{\f12})=0,
$$
which proves \eqref{large-time}.
In addition, by \eqref{point-X} and the above large-time behavior, it holds
$$
|\dot{\mathbf{X}}(t)|\leq C\|v-(v^s)^{-\mathbf{X}}(t,\cdot)\|_{L^{\infty}}
\rightarrow 0, \quad {\rm as}\quad t\rightarrow+\infty,
$$
which proves \eqref{large-X}. Thus the proof of Theorem \ref{thm} is completed.

\hfill $\Box$

%
%
%
%

\section{Uniform-in-time $H^2$-estimates}\label{estimate-(v,h)}
\setcounter{equation}{0}

Throughout this section, $C$ denotes a positive constant which may change from
line to line, but which stays independent on $\d$ (the shock strength) and $\nu$ (the
total variation of the function $a(\xi_1)$). We will consider two smallness conditions,
one on $\d$, and the other on $\d/\nu$. In the argument, $\d$ will be far smaller than $\d/\nu$.

\subsection{Reformulation of the problem}\label{set up}

We introduce a new multi-dimensional effective velocity
\be\label{h}
h:=u-(2\mu+\lambda)\na_{_\xi}v.
\ee
Then the system \eqref{ns-xi} is transformed into
\be\label{new-ns}
\left\{
\begin{array}{l}
\r(\pa_t v-\s\pa_{_{\xi_1}} v+u\cdot\na_{_\xi}v)-\div_{_\xi}h=(2\mu+\lambda)\Delta_{_\xi}v,\\[1mm]
\r(\pa_t h-\s\pa_{_{\xi_1}}h+u\cdot\na_{_\xi}h)+\na_{_\xi}p(v)=R,
\end{array}
\right.
\ee
where
\be\label{R}
R=\f{2\mu+\lambda}{v}(\na_{_\xi} u\cdot\na_{_\xi} v-\div_{_\xi} u\na_{_\xi}v)
-\mu\na_{_\xi}\times\na_{_\xi}\times u.
\ee
We also set
\be\label{h-s}
h^s_1:=u^s_1-(2\mu+\lambda)\pa_{_{\xi_1}}v^s,\qquad h^s:=(h^s_1,0,0)^t.
\ee
We use here a change of variable $\xi_1\rightarrow\xi_1-\mathbf{X}(t)$, the system \eqref{stationary-1} can be rewritten as
\be\label{new-shock}
\left\{
\begin{array}{l}
(\r^s)^{-\mathbf{X}}\big(-\s\pa_{_{\xi_1}}(v^s)^{-\mathbf{X}}+(u^s_1)^{-\mathbf{X}}\pa_{_{\xi_1}}(v^s)^{-\mathbf{X}}\big)
-\pa_{_{\xi_1}}(h^s_1)^{-\mathbf{X}}=(2\mu+\lambda)\pa_{_{\xi_1}}^2(v^s)^{-\mathbf{X}},\\[2mm]
(\r^s)^{-\mathbf{X}}\big(-\s\pa_{_{\xi_1}}(h^s_1)^{-\mathbf{X}}+(u^s_1)^{-\mathbf{X}}\pa_{_{\xi_1}}(h^s_1)^{-\mathbf{X}}\big)
+\pa_{_{\xi_1}}p((v^s)^{-\mathbf{X}})=0.
\end{array}
\right.
\ee
It follows from \eqref{new-ns} and \eqref{new-shock} that,
\be\label{perturb}
\left\{
\begin{array}{l}
\r\pa_t (v-(v^s)^{-\mathbf{X}})-\s\r\pa_{_{\xi_1}}(v-(v^s)^{-\mathbf{X}})+\r u\cdot\na_{_\xi}(v-(v^s)^{-\mathbf{X}})-\div_{_\xi}(h-(h^s)^{-\mathbf{X}})\\[2mm]
\qquad-\dot{\mathbf{X}}(t)\r\pa_{_{\xi_1}}(v^s)^{-\mathbf{X}}+F\pa_{_{\xi_1}}(v^s)^{-\mathbf{X}}=(2\mu+\lambda)\Delta_{_\xi}(v-(v^s)^{-\mathbf{X}}),\\[2mm]
\r\pa_t (h-(h^s)^{-\mathbf{X}})-\s\r\pa_{_{\xi_1}}(h-(h^s)^{-\mathbf{X}})+\r u\cdot\na_{_\xi}(h-(h^s)^{-\mathbf{X}})+\na_{_\xi}\big(p(v)-p((v^s)^{-\mathbf{X}})\big)\\[2mm]
\qquad-\dot{\mathbf{X}}(t)\r\pa_{_{\xi_1}}(h^s)^{-\mathbf{X}}+F\pa_{_{\xi_1}}(h^s)^{-\mathbf{X}}=R,
\end{array}
\right.
\ee
where
\be\label{F}
\ba
F&=-\s(\r-(\r^s)^{-\mathbf{X}})+\r u_1-(\r^s u^s_1)^{-\mathbf{X}}\\
&=-\f{\s_* }{(\r^s)^{-\mathbf{X}}}(\r-(\r^s)^{-\mathbf{X}})+\r(u_1-(u_1^s)^{-\mathbf{X}})\\
&=\s_* \f{v-(v^s)^{-\mathbf{X}}}{v}+\f{h_1-(h_1^s)^{-\mathbf{X}}}{v}
+(2\mu+\lambda)\f{\pa_{_{\xi_1}}(v-(v^s)^{-\mathbf{X}})}{v}.
\ea
\ee
We define the weight function $a(\xi_1)$ by
\be\label{a}
a(\xi_1)=1+\f{\nu}{\d}\big(p(v_-)-p(v^s(\xi_1))\big),
\ee
where the constant $\nu$ is chosen to be small but far bigger than $\d$ such that
\be\label{small}
\d\ll\nu\leq C\sqrt{\d}.
\ee
For definiteness and simplicity, we can choose
$$
\nu=\sqrt{\d}.
$$
Notice that
\be\label{a>1}
1<a(\xi_1)<1+\nu,
\ee
and
\be\label{a-v}
a'(\xi_1)=-\f{\nu}{\d}p'(v^s)v^s_{_{\xi_1}}>0,\qquad
|a'|\sim\f{\nu}{\d}|v^s_{_{\xi_1}}|.
\ee

\begin{lemma}\label{le-1}
Let $a(\xi_1)$ be the weighted function defined by \eqref{a}, it holds
\be\label{ba-expression}
\f{d}{dt}\int_{\mathbb{T}^2}\int_{\mathbb{R}}
a^{-\mathbf{X}}\r \Big(Q(v|(v^s)^{-\mathbf{X}})+\f1{2}{|h-(h^s)^{-\mathbf{X}}|^2}\Big)d\xi_1d\xi'\\
=\dot{\mathbf{X}}(t)\mathbf{Y}(t)+\mathbf{B}(t)-\mathbf{G}(t)-\mathbf{D}(t),
\ee
where
$$
\ba
\mathbf{Y}(t):=&-\int_{\mathbb{T}^2}\int_{\mathbb{R}}a^{-\mathbf{X}}_{_{\xi_1}}\r\Big(
Q(v|(v^s)^{-\mathbf{X}})+\f1{2}{|h-(h^s)^{-\mathbf{X}}|^2}\Big) d\xi_1d\xi'\\
&-\int_{\mathbb{T}^2}\int_{\mathbb{R}}a^{-\mathbf{X}}\r p'((v^s)^{-\mathbf{X}})(v-(v^s)^{-\mathbf{X}})(v^s)^{-\mathbf{X}}_{_{\xi_1}}d\xi_1d\xi'\\
&+\int_{\mathbb{T}^2}\int_{\mathbb{R}}a^{-\mathbf{X}}\r (h^s_1)^{-\mathbf{X}}_{_{\xi_1}}(h_1-(h^s_1)^{-\mathbf{X}})d\xi_1d\xi',\\[4mm]
\mathbf{B}(t):=&\sum_{i=1}^9\mathbf{B}_i(t),
\ea
$$
with
$$
\ba
\mathbf{B}_1(t)&:=
\f{1}{2\s_* }\int_{\mathbb{T}^2}\int_{\mathbb{R}}a^{-\mathbf{X}}_{_{\xi_1}}|p(v)-p((v^s)^{-\mathbf{X}})|^2d\xi_1d\xi',\\
\mathbf{B}_2(t)&:=\s_* \int_{\mathbb{T}^2}\int_{\mathbb{R}}a^{-\mathbf{X}}p(v|(v^s)^{-\mathbf{X}})(v^s)_{_{\xi_1}}^{-\mathbf{X}}d\xi_1d\xi',\\
\mathbf{B}_3(t)&:=\f{\d}{\nu}\int_{\mathbb{T}^2}\int_{\mathbb{R}}\f{a^{-\mathbf{X}}}{\s_*  v}a^{-\mathbf{X}}_{_{\xi_1}}(h_1-(h^s_1)^{-\mathbf{X}})^2d\xi_1d\xi',\\
\mathbf{B}_4(t)&:=\int_{\mathbb{T}^2}\int_{\mathbb{R}}Fa^{-\mathbf{X}}_{_{\xi_1}}\big(Q(v|(v^s)^{-\mathbf{X}})+\f{|h-(h^s)^{-\mathbf{X}}|^2}{2}\big)d\xi_1d\xi',\\
\mathbf{B}_5(t)&:=\int_{\mathbb{T}^2}\int_{\mathbb{R}}a^{-\mathbf{X}}\f{2\mu+\lambda}{v}p'((v^s)^{-\mathbf{X}})(v^s)^{-\mathbf{X}}_{_{\xi_1}}\\
&\qquad\qquad\qquad\quad
\cdot\pa_{_{\xi_1}}(v-(v^s)^{-\mathbf{X}})\big(v-(v^s)^{-\mathbf{X}}-\f{h_1-(h^s_1)^{-\mathbf{X}}}{\s_* }\big)d\xi_1d\xi',\\
\mathbf{B}_6(t)&:=-(2\mu+\lambda)\int_{\mathbb{T}^2}\int_{\mathbb{R}} a^{-\mathbf{X}}\pa_{_{\xi_1}}\big(p(v)-p((v^s)^{-\mathbf{X}})\big)\\
&\qquad\qquad\qquad\quad \cdot\pa_{_{\xi_1}}p((v^s)^{-\mathbf{X}})
\gamma^{-1} \Big(p^{-1-\f{1}{\gamma}}(v)-p^{-1-\f{1}{\gamma}}((v^s)^{-\mathbf{X}})\Big)
d\xi_1d\xi',\\
\mathbf{B}_7(t)&:=-(2\mu+\lambda)\int_{\mathbb{T}^2}\int_{\mathbb{R}}
 a^{-\mathbf{X}}_{_{\xi_1}}\gamma^{-1} p^{-1-\f{1}{\gamma}}(v)\\
 &\qquad\qquad\qquad\qquad\quad \cdot\big(p(v)-p((v^s)^{-\mathbf{X}})\big)
\pa_{_{\xi_1}}\big(p(v)-p((v^s)^{-\mathbf{X}})\big)d\xi_1d\xi',\\
\mathbf{B}_8(t)&:=-(2\mu+\lambda)\int_{\mathbb{T}^2}\int_{\mathbb{R}} a_{_{\xi_1}}^{-\mathbf{X}}\big(p(v)-p((v^s)^{-\mathbf{X}})\big)\pa_{_{\xi_1}}p((v^s)^{-\mathbf{X}})\\
&\qquad\qquad\qquad\qquad\quad
\cdot\gamma^{-1} \Big(p^{-1-\f{1}{\gamma}}(v)-p^{-1-\f{1}{\gamma}}((v^s)^{-\mathbf{X}})\Big)
d\xi_1d\xi',\\
\mathbf{B}_9(t)&:=\int_{\mathbb{T}^2}\int_{\mathbb{R}}a^{-\mathbf{X}}(h-(h^s)^{-\mathbf{X}})\cdot Rd\xi_1d\xi',
\ea
$$
and
$$
\ba
\mathbf{G}(t)&=\s_* \int_{\mathbb{T}^2}\int_{\mathbb{R}}a^{-\mathbf{X}}_{_{\xi_1}}
Q(v|(v^s)^{-\mathbf{X}})d\xi_1d\xi'+\s_* \int_{\mathbb{T}^2}\int_{\mathbb{R}}a^{-\mathbf{X}}_{_{\xi_1}}\f{h_2^3+h_3^2}{2}d\xi_1d\xi'\\
&\quad+\f{\s_* }{2}\int_{\mathbb{T}^2}\int_{\mathbb{R}}a^{-\mathbf{X}}_{_{\xi_1}}
\Big|h_1-(h^s_1)^{-\mathbf{X}}\f{p(v)-p((v^s)^{-\mathbf{X}})}{\s_* }\Big|^2d\xi_1d\xi'\\
&\quad+\int_{\mathbb{T}^2}\int_{\mathbb{R}}a^{-\mathbf{X}}\f{\s_* }{v}|p'((v^s)^{-\mathbf{X}})|(v^s)^{-\mathbf{X}}_{_{\xi_1}}(v-(v^s)^{-\mathbf{X}})^2d\xi_1d\xi'\\
&\quad+(2\mu+\lambda)\int_{\mathbb{T}^2}\int_{\mathbb{R}}a^{-\mathbf{X}}\gamma^{-1} p^{-1-\f{1}{\gamma}}(v)|\na_{_\xi}\big(p(v)-p((v^s)^{-\mathbf{X}})\big)|^2d\xi_1d\xi'\\
&:=\sum_{i=1}^4\mathbf{G}_i(t)+\mathbf{D}(t).
\ea
$$

\end{lemma}

\begin{remark}
Since $\s_*  a^{-\mathbf{X}}_{_{\xi_1}}>0$ and $a^{-\mathbf{X}}>1$, $-\mathbf{G}(t)$ consists of five terms with good sign, while $\mathbf{B}(t)$ consists of bad terms.
\end{remark}

{\textbf{\emph{Proof}}:}
We set $a^{-\mathbf{X}}:=a(\xi_1-\mathbf{X}(t))$. Multiplying \eqref{perturb}$_1$ by $-a^{-\mathbf{X}}\big(p(v)-p((v^s)^{-\mathbf{X}})\big)$, it holds
\be\label{ba-1}
\ba
&\pa_t\big(a^{-\mathbf{X}}\r Q(v|(v^s)^{-\mathbf{X}})\big)-\s\pa_{_{\xi_1}}\big(a^{-\mathbf{X}}\r Q(v|(v^s)^{-\mathbf{X}})\big)
+\div_{_\xi}\big(a^{-\mathbf{X}}\r u Q(v|(v^s)^{-\mathbf{X}})\big)\\
&+a^{-\mathbf{X}}\big(p(v)-p((v^s)^{-\mathbf{X}})\big)\div_{_\xi}(h-(h^s)^{-\mathbf{X}})\\
=&-\dot{\mathbf{X}}(t)a^{-\mathbf{X}}_{_{\xi_1}}\r Q(v|(v^s)^{-\mathbf{X}})-\dot{\mathbf{X}}(t)a^{-\mathbf{X}}\r p'((v^s)^{-\mathbf{X}})(v-(v^s)^{-\mathbf{X}})(v^s)^{-\mathbf{X}}_{_{\xi_1}}\\
&+(-\s\r+\r u_1)a_{_{\xi_1}}^{-\mathbf{X}}Q(v|(v^s)^{-\mathbf{X}})-(-\s\r+\r u_1)a^{-\mathbf{X}} p(v|(v^s)^{-\mathbf{X}})(v^s)^{-\mathbf{X}}_{_{\xi_1}}\\
&+a^{-\mathbf{X}}F\big(p(v)-p((v^s)^{-\mathbf{X}})\big)(v^s)^{-\mathbf{X}}_{_{\xi_1}}\\
&-(2\mu+\lambda)\div_{_\xi}\big( a^{-\mathbf{X}}\big(p(v)-p((v^s)^{-\mathbf{X}})\big)\na_{_\xi}(v-(v^s)^{-\mathbf{X}})\big)\\
&+(2\mu+\lambda)\na_{_\xi}\big( a^{-\mathbf{X}}\big(p(v)-p((v^s)^{-\mathbf{X}})\big)\cdot\na_{_\xi}(v-(v^s)^{-\mathbf{X}}).
\ea
\ee
Using \eqref{F} the definition of $F$, and \eqref{bar-sigma} the definition of $\s_* $, it holds
\be\label{ba-5}
-\s\r+\r u_1=-\s(\r^s)^{-\mathbf{X}}+(\r^su^s_1)^{-\mathbf{X}}
-\s(\r-(\r^s)^{-\mathbf{X}})+(\r u_1-(\r^su^s_1)^{-\mathbf{X}})=-\s_* +F.
\ee
Thus, we have
$$
(-\s\r+\r u_1)a_{_{\xi_1}}^{-\mathbf{X}}Q(v|(v^s)^{-\mathbf{X}})=(-\s_* +F)a_{_{\xi_1}}^{-\mathbf{X}}Q(v|(v^s)^{-\mathbf{X}})
$$
and
$$
\ba
&-(-\s\r+\r u_1)a^{-\mathbf{X}} p(v|(v^s)^{-\mathbf{X}})(v^s)^{-\mathbf{X}}_{_{\xi_1}}
+a^{-\mathbf{X}}F\big(p(v)-p((v^s)^{-\mathbf{X}})\big)(v^s)^{-\mathbf{X}}_{_{\xi_1}}\\
&=\s_*  a^{-\mathbf{X}}p(v|(v^s)^{-\mathbf{X}})(v^s)^{-\mathbf{X}}_{_{\xi_1}}+Fa^{-\mathbf{X}}p'((v^s)^{-\mathbf{X}})(v-(v^s)^{-\mathbf{X}})(v^s)^{-\mathbf{X}}_{_{\xi_1}}.
\ea
$$
Notice that
$$
\na_{_\xi}v=\f{\na_{_\xi} p(v)}{p'(v)}=\f{\na_{_\xi} p(v)}{-\gamma p^{1+\f{1}{\gamma}}(v)}.
$$
Hence, the last term on the right hand side of \eqref{ba-1} holds
\begin{align}\label{cal-viscosity}  \nm
&(2\mu+\lambda)\na_{_\xi}\big( a^{-\mathbf{X}}\big(p(v)-p((v^s)^{-\mathbf{X}})\big)\cdot\na_{_\xi}(v-(v^s)^{-\mathbf{X}})\\  \nm
=&(2\mu+\lambda) a^{-\mathbf{X}}\na_{_\xi}\big(p(v)-p((v^s)^{-\mathbf{X}})\big)\cdot\Big(\f{\na_{_\xi} p(v)}{-\gamma p^{1+\f{1}{\gamma}}(v)}-\f{\na_{_\xi} p((v^s)^{-\mathbf{X}})}{-\gamma p^{1+\f{1}{\gamma}}((v^s)^{-\mathbf{X}})}\Big)\\  \nm
&+(2\mu+\lambda) a^{-\mathbf{X}}_{_{\xi_1}}\big(p(v)-p((v^s)^{-\mathbf{X}})\big)\Big(\f{\pa_{_{\xi_1}} p(v)}{-\gamma p^{1+\f{1}{\gamma}}(v)}-\f{\pa_{_{\xi_1}} p((v^s)^{-\mathbf{X}})}{-\gamma p^{1+\f{1}{\gamma}}((v^s)^{-\mathbf{X}})}\Big)\\
=&-(2\mu+\lambda) a^{-\mathbf{X}}\f{|\na_{_\xi}\big(p(v)-p((v^s)^{-\mathbf{X}})\big)|^2}{\gamma p^{1+\f{1}{\gamma}}(v)}\\  \nm
&-(2\mu+\lambda) a^{-\mathbf{X}}\pa_{_{\xi_1}}\big(p(v)-p((v^s)^{-\mathbf{X}})\big)\pa_{_{\xi_1}}p((v^s)^{-\mathbf{X}})
\Big(\f{1}{\gamma p^{1+\f{1}{\gamma}}(v)}-\f{1}{\gamma p^{1+\f{1}{\gamma}}((v^s)^{-\mathbf{X}})}\Big)\\  \nm
&-(2\mu+\lambda) a^{-\mathbf{X}}_{_{\xi_1}}\big(p(v)-p((v^s)^{-\mathbf{X}})\big)
\f{\pa_{_{\xi_1}}\big(p(v)-p((v^s)^{-\mathbf{X}})\big)}{\gamma p^{1+\f{1}{\gamma}}(v)}\\  \nm
&-(2\mu+\lambda) a_{_{\xi_1}}^{-\mathbf{X}}\big(p(v)-p((v^s)^{-\mathbf{X}})\big)\pa_{_{\xi_1}}p((v^s)^{-\mathbf{X}})
\Big(\f{1}{\gamma p^{1+\f{1}{\gamma}}(v)}-\f{1}{\gamma p^{1+\f{1}{\gamma}}((v^s)^{-\mathbf{X}})}\Big).
\end{align}

Multiplying \eqref{perturb}$_2$ by $a^{-\mathbf{X}}(h-(h^s)^{-\mathbf{X}})$, we have
\be\label{ba-2}
\ba
&\pa_t\big(a^{-\mathbf{X}}\r\f{|h-(h^s)^{-\mathbf{X}}|^2}{2}\big)-\s\pa_{_{\xi_1}}\big(a^{-\mathbf{X}}\r\f{|h-(h^s)^{-\mathbf{X}}|^2}{2}\big)
+\div_{_\xi}\big(a^{-\mathbf{X}}\r u\f{|h-(h^s)^{-\mathbf{X}}|^2}{2}\big)\\
&+\div_{_\xi}\big(a^{-\mathbf{X}}\big(p(v)-p((v^s)^{-\mathbf{X}})\big)(h-(h^s)^{-\mathbf{X}})\big)
-a^{-\mathbf{X}}\big(p(v)-p((v^s)^{-\mathbf{X}})\big)\div_{_\xi}(h-(h^s)^{-\mathbf{X}})\\
=&-\dot{\mathbf{X}}(t)a_{_{\xi_1}}^{-\mathbf{X}}\r\f{|h-(h^s)^{-\mathbf{X}}|^2}{2}+\dot{\mathbf{X}}(t) a^{-\mathbf{X}}\r (h^s_1)^{-\mathbf{X}}_{_{\xi_1}}(h_1-(h^s_1)^{-\mathbf{X}})\\
&+a_{_{\xi_1}}^{-\mathbf{X}}\big(p(v)-p((v^s)^{-\mathbf{X}})\big)(h_1-(h^s_1)^{-\mathbf{X}})
+(-\s\r+\r u_1)a_{_{\xi_1}}^{-\mathbf{X}}\f{|h-(h^s)^{-\mathbf{X}}|^2}{2}\\
&-Fa^{-\mathbf{X}}(h^s_1)^{-\mathbf{X}}_{_{\xi_1}}(h_1-(h_1^s)^{-\mathbf{X}})+a^{-\mathbf{X}}(h-(h^s)^{-\mathbf{X}})\cdot R.
\ea
\ee

Before we add \eqref{ba-1} and \eqref{ba-2} together, direct calculations yield
\be\label{ba-3}
\ba
&-\s_*  a_{_{\xi_1}}^{-\mathbf{X}}\f{|h-(h^s)^{-\mathbf{X}}|^2}{2}
+a_{_{\xi_1}}^{-\mathbf{X}}\big(p(v)-p((v^s)^{-\mathbf{X}})\big)(h_1-(h^s_1)^{-\mathbf{X}})\\
=&-\f{\s_* }{2}a^{-\mathbf{X}}_{_{\xi_1}}\Big|h_1-(h^s_1)^{-\mathbf{X}}-\f{p(v)-p((v^s)^{-\mathbf{X}})}{\s_* }\Big|^2
+a^{-\mathbf{X}}_{_{\xi_1}}\f{|p(v)-p((v^s)^{-\mathbf{X}})|^2}{2\s_* }
-\s_*  a^{-\mathbf{X}}_{_{\xi_1}}\f{h_2^2+h_3^2}{2}.
\ea
\ee

We treat the perturbed flux term in the Eulerian coordinates along the shock wave propagation direction,
which is different from that of in the Lagrangian coordinates. It follows from \eqref{new-shock}$_2$ that
$\s_*  \pa_{_{\xi_1}}(h^s_1)^{-\mathbf{X}}=\pa_{_{\xi_1}}p((v^s)^{-\mathbf{X}})$. Hence, using \eqref{F}, we have
\be\label{ba-4}
\ba
&\quad Fa^{-\mathbf{X}}p'((v^s)^{-\mathbf{X}})(v-(v^s)^{-\mathbf{X}})(v^s)^{-\mathbf{X}}_{_{\xi_1}}-Fa^{-\mathbf{X}}(h^s_1)^{-\mathbf{X}}_{_{\xi_1}}(h_1-(h_1^s)^{-\mathbf{X}})\\
&=a^{-\mathbf{X}}\f{\s_* }{v}p'((v^s)^{-\mathbf{X}})(v^s)^{-\mathbf{X}}_{_{\xi_1}}(v-(v^s)^{-\mathbf{X}})^2
+\f{\d}{\nu}\f{a^{-\mathbf{X}}}{\s_*  v}a^{-\mathbf{X}}_{_{\xi_1}}(h_1-(h^s_1)^{-\mathbf{X}})^2\\
&+a^{-\mathbf{X}}\f{2\mu+\lambda}{v}p'((v^s)^{-\mathbf{X}})(v^s)^{-\mathbf{X}}_{_{\xi_1}}\pa_{_{\xi_1}}(v-(v^s)^{-\mathbf{X}})
\Big(v-(v^s)^{-\mathbf{X}}-\f{h_1-(h^s_1)^{-\mathbf{X}}}{\s_* }\Big).
\ea
\ee
Adding \eqref{ba-1} and \eqref{ba-2} together, integrating the resultant equation by parts over $\Omega:=\mathbb{R}\times\mathbb{T}^2$,
and using \eqref{cal-viscosity}, \eqref{ba-3} and \eqref{ba-4}, we can obtain \eqref{ba-expression}. The proof of Lemma \ref{le-1} is completed.

\hfill $\Box$

In order to derive the a-contraction property of the viscous shock wave, we decompose the function $\mathbf{Y}(t)$ in Lemma \ref{le-1} as
$$
\mathbf{Y}(t):=\sum_{i=1}^5\mathbf{Y}_i(t),
$$
where
\begin{align} \nm
\mathbf{Y}_1(t)&:=\int_{\mathbb{T}^2}\int_{\mathbb{R}}\f{a^{-\mathbf{X}}}{\s_* }\r(h_1^s)^{-\mathbf{X}}_{_{\xi_1}}\big(p(v)-p((v^s)^{-\mathbf{X}})\big)d\xi_1d\xi',\\ \nm
\mathbf{Y}_2(t)&:=-\int_{\mathbb{T}^2}\int_{\mathbb{R}}a^{-\mathbf{X}}\r p'((v^s)^{-\mathbf{X}})(v-(v^s)^{-\mathbf{X}})(v^s)^{-\mathbf{X}}_{_{\xi_1}}d\xi_1d\xi',\\  \nm
\mathbf{Y}_3(t)&:=\int_{\mathbb{T}^2}\int_{\mathbb{R}}a^{-\mathbf{X}}\r(h_1^s)^{-\mathbf{X}}_{_{\xi_1}}\Big(
h_1-(h_1^s)^{-\mathbf{X}}-\f{p(v)-p((v^s)^{-\mathbf{X}})}{\s_* }\Big)d\xi_1d\xi',\\  \nm
\mathbf{Y}_4(t)&:=-\f12\int_{\mathbb{T}^2}\int_{\mathbb{R}}a^{-\mathbf{X}}_{_{\xi_1}}\r
\Big(h_1-(h_1^s)^{-\mathbf{X}}-\f{p(v)-p((v^s)^{-\mathbf{X}})}{\s_* }\Big)\\  \nm
&\qquad\qquad\qquad\cdot\Big(h_1-(h_1^s)^{-\mathbf{X}}+\f{p(v)-p((v^s)^{-\mathbf{X}})}{\s_* }\Big)d\xi_1d\xi',\\  \nm
\mathbf{Y}_5(t)&:=-\int_{\mathbb{T}^2}\int_{\mathbb{R}}a^{-\mathbf{X}}_{_{\xi_1}}\r\big(
Q(v|(v^s)^{-\mathbf{X}})+\f{h_2^2+h_3^2}{2}\big)d\xi_1d\xi'\\  \nm
&\quad\,\,-\int_{\mathbb{T}^2}\int_{\mathbb{R}}a^{-\mathbf{X}}_{_{\xi_1}}\r\f{|p(v)-p((v^s)^{-\mathbf{X}})|^2}{2\s_* ^2}d\xi_1d\xi'.
\end{align}
Notice that
\be\label{X=Y1+Y2}
\dot{\mathbf{X}}(t)=-\f{M}{\d}(\mathbf{Y}_1(t)+\mathbf{Y}_2(t)),
\ee
and so
\be
\dot{\mathbf{X}}(t)\mathbf{Y}(t)=-\f{\d}{M}|\dot{\mathbf{X}}(t)|^2+\dot{\mathbf{X}}(t)\sum_{i=3}^5\mathbf{Y}_i(t).
\ee

Then we have the following Lemma.

\begin{lemma}\label{es-good}
There exits uniform-in-time $C>0$ such that for $\forall t\in[0,T]$
\be\label{es-key}
\ba
&-\f{\d}{2M}|\dot{\mathbf{X}}(t)|^2+\mathbf{B}_1(t)+\mathbf{B}_2(t)+\mathbf{B}_3(t)-\mathbf{G}_{1}(t)-\mathbf{G}_4(t)-\f34\mathbf{D}(t)\\
&\leq -C\int_{\mathbb{T}^2}\int_{\mathbb{R}}|(v^s)^{-\mathbf{X}}_{_{\xi_1}}||p(v)-p((v^s)^{-\mathbf{X}})|^2d\xi_1d\xi'\\
&+C\int_{\mathbb{T}^2}\int_{\mathbb{R}}a^{-\mathbf{X}}_{_{\xi_1}}|p(v)-p((v^s)^{-\mathbf{X}})|^3d\xi_1d\xi'
+\f{1}{40}\mathbf{G}_3(t).
\ea
\ee
\end{lemma}

{\textbf{\emph{Proof}}:} We now rewrite the above functions with respect to the following variables
\be\label{w}
w:=p(v)-p((v^s)^{-\mathbf{X}}),\quad y_1:=\f{p(v_-)-p((v^s(\xi_1))^{-\mathbf{X}})}{\d},\quad y'=\xi',\quad {\rm i.e.},\quad(y_2,y_3)=(\xi_2,\xi_3).
\ee
We use a change of variable $\xi_1\in\mathbb{R}\mapsto y_1\in[0,1]$. Then it follows from \eqref{a} that $a^{-\mathbf{X}}(\xi_1)=1+\nu y_1$ and
\be\label{change-variable}
\f{dy_1}{d\xi_1}=-\f{1}{\d}p((v^s)^{-\mathbf{X}})_{_{\xi_1}},\quad a_{_{\xi_1}}^{-\mathbf{X}}=\nu\f{dy_1}{d\xi_1},\quad |a^{-\mathbf{X}}-1|\leq \nu=\sqrt{\d}.
\ee
To perform the sharp estimates, we will consider the $O(1)-$constants:
$$
\s_-:=\sqrt{-p'(v_-)},\qquad \a_-:=\f{\gamma+1}{2\gamma\s_-p(v_-)},
$$
which are indeed independent of the small constant $\d$, since $\f{v_+}{2}\leq v_-\leq v_+$. Note that
\be\label{dengjia-1}
|\s_* -\s_-|\leq C\d,
\ee
which together with $\s_-^2=-p'(v_-)=\gamma p^{1+\f{1}{\gamma}}(v_-)$ implies
\be\label{dengjia-2}
|\s_-^2+p'((v^s)^{-\mathbf{X}})|\leq C\d,\qquad \left|\f{1}{\s_-^2}-\f{1}{\gamma p^{1+\f{1}{\gamma}}((v^s)^{-\mathbf{X}})}\right|\leq C\d.
\ee

$\bullet$ {\textbf{Estimate on $-\f{\d}{2M}|\dot{\mathbf{X}}(t)|^2$}:} To do this, we will control $\mathbf{Y}_1(t)$ and $\mathbf{Y}_2(t)$ due to \eqref{X=Y1+Y2}. Using \eqref{bar-sigma}, system \eqref{new-shock} is transformed into
\be\label{new-shock-2}
\left\{
\begin{array}{l}
-\s_* (v^s)^{-\mathbf{X}}_{_{\xi_1}}-(h^s_1)^{-\mathbf{X}}_{_{\xi_1}}=(2\mu+\lambda)(v^s)^{-\mathbf{X}}_{_{\xi_1\xi_1}},\\[2mm]
-\s_* (h^s_1)^{-\mathbf{X}}_{_{\xi_1}}+p((v^s)^{-\mathbf{X}})_{_{\xi_1}}=0.
\end{array}
\right.
\ee
Using \eqref{new-shock-2}$_2$ and the new variable \eqref{change-variable}, we have
$$
\mathbf{Y}_1(t)=\int_{\mathbb{T}^2}\int_{\mathbb{R}}\f{a^{-\mathbf{X}}}{\s_* ^2 v}p((v^s)^{-\mathbf{X}})_{_{\xi_1}}\big(p(v)-p((v^s)^{-\mathbf{X}})\big)d\xi_1d\xi'
=-\f{\d}{\s_* ^2}\int_{\mathbb{T}^2}\int_0^1a^{-\mathbf{X}}\f{w}{v}dy_1dy'.
$$
Using \eqref{dengjia-1} and $|a^{-\mathbf{X}}-1|\leq \nu$, it holds
\be\label{es-Y1}
\left|\mathbf{Y}_1(t)+\f{\d}{\s_- ^2 v_-}\int_{\mathbb{T}^2}\int_0^1 w dy_1dy'\right|
\leq C\d(\nu+\d)\int_{\mathbb{T}^2}\int_0^1|w|dy_1dy'.
\ee
For
$$
\mathbf{Y}_2(t)=-\int_{\mathbb{T}^2}\int_{\mathbb{R}}a^{-\mathbf{X}}\r p((v^s)^{-\mathbf{X}})_{_{\xi_1}}(v-(v^s)^{-\mathbf{X}})d\xi_1d\xi'
=\d\int_{\mathbb{T}^2}\int_0^1 a^{-\mathbf{X}}\f{v-(v^s)^{-\mathbf{X}}}{v}dy_1dy',
$$
we observe that since (by considering $v=p(v)^{-\f{1}{\gamma}}$)
$$
\left|v-(v^s)^{-\mathbf{X}}+\f{p(v)-p((v^s)^{-\mathbf{X}})}{\gamma p^{1+\f{1}{\gamma}}((v^s)^{-\mathbf{X}})}\right|\leq C|p(v)-p((v^s)^{-\mathbf{X}})|^2,
$$
then it holds
$$
\left|v-(v^s)^{-\mathbf{X}}+\f{1}{\s_-^2}\big(p(v)-p((v^s)^{-\mathbf{X}})\big)\right|\leq C(\d+\chi)|p(v)-p((v^s)^{-\mathbf{X}})|.
$$
This implies
\be\label{es-Y2}
\left|\mathbf{Y}_2(t)+\f{\d}{\s_- ^2 v_-}\int_{\mathbb{T}^2}\int_0^1 w dy_1dy'\right|
\leq C\d(\nu+\d+\chi)\int_{\mathbb{T}^2}\int_0^1|w|dy_1dy'.
\ee
By \eqref{X=Y1+Y2}, \eqref{es-Y1} and \eqref{es-Y2}, we have
$$
\ba
\left|\dot{\mathbf{X}}(t)-\f{2M}{\s_-^2v_-}\int_{\mathbb{T}^2}\int_0^1 w dy_1 dy'\right|
&=\left|\sum_{i=1}^2\f{M}{\d}\Big(\mathbf{Y}_i(t)+\f{\d}{\s_-^2v_-}\int_{\mathbb{T}^2}\int_0^1 w dy_1 dy'\Big)\right|\\
&\leq C(\nu+\d+\chi)\int_{\mathbb{T}^2}\int_0^1 |w| dy_1 dy',
\ea
$$
which yields
$$
\ba
\left(\left|\f{2M}{\s_-^2 v_-}\int_{\mathbb{T}^2}\int_0^1 w dy_1 dy'\right|
-|\dot {\mathbf{X}}(t)|\right)^2&\leq C(\nu+\d+\chi)^2\left(\int_{\mathbb{T}^2}\int_0^1 |w| dy_1 dy'\right)^2\\
&\leq C(\nu+\d+\chi)^2\int_{\mathbb{T}^2}\int_0^1 |w|^2 dy_1 dy',
\ea
$$
which together with the algebraic inequality $\f{p^2}{2}-q^2\leq (p-q)^2$ for all $p,q\geq0$ indicate
$$
\f{2M^2}{\s_-^4v_-^2}\left(\int_{\mathbb{T}^2}\int_0^1 w dy_1 dy'\right)^2
-|\dot{\mathbf{X}}(t)|^2\leq C(\nu+\d+\chi)^2\int_{\mathbb{T}^2}\int_0^1 |w|^2 dy_1 dy'.
$$
Thus, we can get
\be\label{es-X}
-\f{\d}{2M}|\dot{\mathbf{X}}(t)|^2\leq -\f{M\d}{\s_-^4 v_-^2}\left(\int_{\mathbb{T}^2}\int_0^1 w dy_1 dy'\right)^2
+C\d(\nu+\d+\chi)^2\int_{\mathbb{T}^2}\int_0^1 |w|^2 dy_1 dy'.
\ee



$\bullet$ \textbf{Change of variable for $\mathbf{B}_i(t)~(i=1,2,3)$}: By the change of variable,
using \eqref{dengjia-1}, we have
\be\label{B1}
\ba
\mathbf{B}_1(t)&=\f{\nu}{2\s_* }\int_{\mathbb{T}^2}\int_0^1w^2dy_1dy'
=\f{\nu}{2\s_-}\int_{\mathbb{T}^2}\int_0^1w^2dy_1dy'
+\f{\nu}{2}\Big(\f{1}{\s_* }-\f{1}{\s_-}\Big)\int_{\mathbb{T}^2}\int_0^1w^2dy_1dy'\\
&\leq\f{\nu}{2\s_-}\int_{\mathbb{T}^2}\int_0^1w^2dy_1dy'+C\nu\d\int_{\mathbb{T}^2}\int_0^1w^2dy_1dy'.
\ea
\ee
For $\mathbf{B}_2(t)$, by the change of variable, and using \eqref{use-1}, it holds
\be\label{B2}
\ba
\mathbf{B}_2(t)&=\s_* \d\int_{\mathbb{T}^2}\int_0^1(1+\nu y_1)p(v|(v^s)^{-\mathbf{X}})\f{1}{|p'((v^s)^{-\mathbf{X}})|}dy_1dy'\\
&\leq \s_* \d(1+\nu)\int_{\mathbb{T}^2}\int_0^1\Big(\f{\gamma+1}{2\gamma}\f{1}{p((v^s)^{-\mathbf{X}})}+C\chi\Big)\f{|p(v)-p((v^s)^{-\mathbf{X}})|^2}{|p'((v^s)^{-\mathbf{X}})|}dy_1dy'\\
&\leq \s_* \d(1+\nu)\int_{\mathbb{T}^2}\int_0^1\Big(\a_-\f{\s_- p(v_-)}{p((v^s)^{-\mathbf{X}})}+C\chi\Big)\f{w^2}{|p'((v^s)^{-\mathbf{X}})|}dy_1dy'\\
&\leq \d\a_-(1+C(\nu+\d+\chi))\int_{\mathbb{T}^2}\int_0^1 w^2 dy_1dy',
\ea
\ee
where in the last inequality we have used \eqref{dengjia-1} and \eqref{dengjia-2}.

For $\mathbf{B}_3(t)$, using the algebraic inequality $p^2=(q+p-q)^2\leq (1+\kappa)q^2+(1+\f{1}{\kappa})(p-q)^2$ for
$\k>0$, it holds
$$
\ba
\mathbf{B}_3(t)&\leq (1+\k)\f{\d}{\nu}\int_{\mathbb{T}^2}\int_{\mathbb{R}}\f{a^{-\mathbf{X}}}{\s_* ^3 v}a^{-\mathbf{X}}_{_{\xi_1}}
|p(v)-p((v^s)^{-\mathbf{X}})|^2d\xi_1d\xi'+C(1+\f{1}{\k})\f{\d}{\nu}\mathbf{G}_3(t).
\ea
$$
Since
$$
\ba
(1+\k)\f{1}{\s_-^3v_-\a_-}=(1+\k)\f{1}{\s_-^3v_-}\f{2\gamma\s_-p(v_-)}{\gamma+1}
=\f{2}{\gamma+1}(1+\k).
\ea
$$
By the change of variable,
and using \eqref{dengjia-1}, we have
$$
\ba
&\f{\d}{\nu}(1+\k)\int_{\mathbb{T}^2}\int_{\mathbb{R}}\f{a^{-\mathbf{X}}}{\s_* ^3 v}a^{-\mathbf{X}}_{_{\xi_1}}
|p(v)-p((v^s)^{-\mathbf{X}})|^2d\xi_1d\xi'\\
=&\d(1+\k)\int_{\mathbb{T}^2}\int_0^1\f{(1+\nu y_1)}{\s_-^3v_-}\f{\s_-^3v_-}{\s_* ^3 v}w^2dy_1dy'\\
\leq &\d(1+\k)(1+C(\nu+\d+\chi))\int_{\mathbb{T}^2}\int_0^1\f{1}{\s_-^3v_-}w^2dy_1dy'\\
\leq &\f{2}{\gamma+1}(1+\k)\d\a_-(1+C(\nu+\d+\chi))\int_{\mathbb{T}^2}\int_0^1w^2dy_1dy'.
\ea
$$
Thus, we have
\be\label{B3}
\mathbf{B}_3(t)\leq \f{2}{\gamma+1}(1+\k)\d\a_-(1+C(\nu+\d+\chi))\int_{\mathbb{T}^2}\int_0^1w^2dy_1dy'
+C(1+\f{1}{\k})\f{\d}{\nu}\mathbf{G}_3(t).
\ee


$\bullet$ \textbf{Change of variable for $\mathbf{G}_1(t)$, $\mathbf{G}_4(t)$}: For $\mathbf{G}_1(t)$,
we first use \eqref{use-2} to split it into two parts:
\be\label{G1}
\ba
\mathbf{G}_1(t)&\geq\underbrace{\s_* \int_{\mathbb{T}^2}\int_{\mathbb{R}}a^{-\mathbf{X}}_{_{\xi_1}}
\f{|p(v)-p((v^s)^{-\mathbf{X}})|^2}{2\gamma p^{1+\f{1}{\gamma}}((v^s)^{-\mathbf{X}})}d\xi_1d\xi'}_{\mathbf{G}_{1,1}(t)}\\
&-\s_* \int_{\mathbb{T}^2}\int_{\mathbb{R}}a^{-\mathbf{X}}_{_{\xi_1}}
\f{1+\gamma}{3\gamma^2}\f{(p(v)-p((v^s)^{-\mathbf{X}}))^3}{p^{2+\f{1}{\gamma}}((v^s)^{-\mathbf{X}})}d\xi_1d\xi'.
\ea
\ee
We only need to do the change of variable for the good term $\mathbf{G}_{1,1}(t)$ as follow: By \eqref{dengjia-1},
\eqref{dengjia-2} and the change of variable,
$$
\ba
\mathbf{G}_{1,1}(t)&\geq \f{\s_* }{2\s_-^2}(1-C\d)
\int_{\mathbb{T}^2}\int_{\mathbb{R}}a^{-\mathbf{X}}_{_{\xi_1}}|p(v)-p((v^s)^{-\mathbf{X}})|^2d\xi_1d\xi'\\
&\geq \f{\nu}{2\s_-}(1-C\d)
\int_{\mathbb{T}^2}\int_0^1 w^2 dy_1dy',
\ea
$$
which together with \eqref{B1} yields
\be\label{es-B1+G21}
\mathbf{B}_1(t)-\mathbf{G}_{1,1}(t)\leq C\nu\d\int_{\mathbb{T}^2}\int_0^1 w^2 dy_1dy'.
\ee

For $\mathbf{G}_4(t)$, using the mean value theorem,
$$
v-(v^s)^{-\mathbf{X}}=\f{p(v)-p((v^s)^{-\mathbf{X}})}{p'(\z)},\quad {\rm for}\,\  \z \,\  {\rm between}\,\  v\, \ {\rm and}\, \ (v^s)^{-\mathbf{X}}.
$$
Using \eqref{change-variable} and the change of variable, it holds
\be\label{G4}
\ba
\mathbf{G}_4(t)&=\f{\d}{\nu}\int_{\mathbb{T}^2}\int_{\mathbb{R}}a^{-\mathbf{X}}a^{-\mathbf{X}}_{_{\xi_1}}\f{\s_* }{v}|v-(v^s)^{-\mathbf{X}}|^2d\xi_1d\xi'
\geq\d\int_{\mathbb{T}^2}\int_0^1\f{\s_* }{v}\f{w^2}{|p'(v_-)|^2}\Big|\f{p'(v_-)}{p'(\z)}\Big|^2dy_1dy'\\
&\geq\d(1-C(\d+\chi))\f{1}{\s_-^3v_-}\int_{\mathbb{T}^2}\int_0^1w^2dy_1dy'
=\f{2}{\gamma+1}\d\a_-(1-C(\d+\chi))\int_{\mathbb{T}^2}\int_0^1w^2dy_1dy'.
\ea
\ee


$\bullet$ {\textbf{Change of variable for $\mathbf{D}(t)$}}: First, using \eqref{a>1} $a^{-\mathbf{X}}>1$ and the change of
variable, it holds
$$
\ba
\mathbf{D}(t)&\geq (2\mu+\lambda)\int_{\mathbb{T}^2}\int_{\mathbb{R}}
\f{|\pa_{_{\xi_1}}\big(p(v)-p((v^s)^{-\mathbf{X}})\big)|^2}{\gamma p^{1+\f{1}{\gamma}}(v)}d\xi_1d\xi'\\
&+(2\mu+\lambda)\int_{\mathbb{T}^2}\int_{\mathbb{R}}\f{|\na_{_{\xi'}}\big(p(v)-p((v^s)^{-\mathbf{X}})\big)|^2}{\gamma p^{1+\f{1}{\gamma}}(v)}d\xi_1d\xi'
\\
&=\underbrace{(2\mu+\lambda)\int_{\mathbb{T}^2}\int_{0}^1\f{|\pa_{y_1}w|^2}{\gamma p^{1+\f{1}{\gamma}}(v)}
\Big(\f{dy_1}{d\xi_1}\Big)dy_1dy'}_{\mathbf{D}_{\rm I}(t)}
+\underbrace{(2\mu+\lambda)\int_{\mathbb{T}^2}\int_{0}^1\f{|\na_{y'}w|^2}{\gamma p^{1+\f{1}{\gamma}}(v)}
\Big(\f{d\xi_1}{dy_1}\Big)dy_1dy'}_{\mathbf{D}_{\rm II}(t)}.
\ea
$$
On one hand, integrating \eqref{new-shock-2} over $(-\infty,\xi]$ yields
$$
(2\mu+\lambda)(v^s)^{-\mathbf{X}}_{_{\xi_1}}=-\s_* ((v^s)^{-\mathbf{X}}-v_-)-\f{1}{\s_* }\big(p((v^s)^{-\mathbf{X}})-p(v_-)\big).
$$
On the other hand,
$$
(v^s)^{-\mathbf{X}}_{_{\xi_1}}=\f{p((v^s)^{-\mathbf{X}})_{_{\xi_1}}}{p'((v^s)^{-\mathbf{X}})}=\f{\d}{\gamma p^{1+\f{1}{\gamma}}((v^s)^{-\mathbf{X}})}
\f{dy_1}{d\xi_1}.
$$
Hence, we have
$$
\ba
(2\mu+\lambda)\f{\d}{\gamma p^{1+\f{1}{\gamma}}((v^s)^{-\mathbf{X}})}
\f{dy_1}{d\xi_1}&=-\s_* ((v^s)^{-\mathbf{X}}-v_-)-\f{1}{\s_* }\big(p((v^s)^{-\mathbf{X}})-p(v_-)\big)\\
&=\f{-1}{\s_* }\Big(\s_* ^2((v^s)^{-\mathbf{X}}-v_-)+\big(p((v^s)^{-\mathbf{X}})-p(v_-)\big)\Big),
\ea
$$
which together with $\s_* ^2=-\f{p(v_-)-p(v_+)}{v_--v_+}$ leads to
$$
\ba
&(2\mu+\lambda)\f{\d}{\gamma p^{1+\f{1}{\gamma}}((v^s)^{-\mathbf{X}})}
\f{dy_1}{d\xi_1}\\
=&\f{-1}{\s_* (v_+-v_-)}\Big(\big(p(v_-)-p(v_+)\big)((v^s)^{-\mathbf{X}}-v_-)
+\big(p((v^s)^{-\mathbf{X}})-p(v_-)\big)(v_+-v_-)\Big)\\
=&\f{-1}{\s_* (v_+-v_-)}\Big(\big(p((v^s)^{-\mathbf{X}})-p(v_+)\big)((v^s)^{-\mathbf{X}}-v_-)
-((v^s)^{-\mathbf{X}}-v_+)\big(p((v^s)^{-\mathbf{X}})-p(v_-)\big)\Big).
\ea
$$
Recall that $y_1=\f{p(v_-)-p((v^s)^{-\mathbf{X}})}{\d}$ and $1-y_1=\f{p((v^s)^{-\mathbf{X}})-p(v_+)}{\d}$, it holds
$$
\f{1}{y_1(1-y_1)}\f{2\mu+\lambda}{\gamma p^{1+\f{1}{\gamma}}((v^s)^{-\mathbf{X}})}\f{dy_1}{d\xi_1}
=\f{\d}{\s_* (v_+-v_-)}\left(\f{(v^s)^{-\mathbf{X}}-v_-}{p((v^s)^{-\mathbf{X}})-p(v_-)}
-\f{(v^s)^{-\mathbf{X}}-v_+}{p((v^s)^{-\mathbf{X}})-p(v_+)}\right).
$$
Then
$$
\ba
&\left|\f{1}{y_1(1-y_1)}\f{2\mu+\lambda}{\gamma p^{1+\f{1}{\gamma}}((v^s)^{-\mathbf{X}})}\f{dy_1}{d\xi_1}
-\f{\d p''(v_-)}{2\s_-(p'(v_-))^2}\right|\\
\leq& \left|\f{1}{y_1(1-y_1)}\f{2\mu+\lambda}{\gamma p^{1+\f{1}{\gamma}}((v^s)^{-\mathbf{X}})}\f{dy_1}{d\xi_1}
-\f{\d p''(v_-)}{2 \s_*(p'(v_-))^2}\right|
+\f{\d p''(v_-)}{2(p'(v_-))^2}\left|\f{1}{\s_* }-\f{1}{\s_-}\right|=:I_1+I_2.
\ea
$$
Using Lemma \ref{inver-pressure}, we have
$$
\ba
I_1=\left|\f{\d}{\s_* (v_+-v_-)}\left(\f{(v^s)^{-\mathbf{X}}-v_-}{p((v^s)^{-\mathbf{X}})-p(v_-)}
+\f{(v^s)^{-\mathbf{X}}-v_+}{p(v_+)-p((v^s)^{-\mathbf{X}})}+\f12\f{p''(v_-)}{(p'(v_-))^2}(v_--v_+)
\right)\right|\leq C\d^2.
\ea
$$
Since it follows from \eqref{dengjia-1} that $I_2\leq C\d^2$, it holds
\be\label{coeffi-1}
\left|\f{1}{y_1(1-y_1)}\f{2\mu+\lambda}{\gamma p^{1+\f{1}{\gamma}}((v^s)^{-\mathbf{X}})}\f{dy_1}{d\xi_1}
-\f{\d p''(v_-)}{2\s_-(p'(v_-))^2}\right|\leq C\d^2.
\ee
In addition,
$$
\left|\left(\f{p((v^s)^{-\mathbf{X}})}{p(v)}\right)^{1+\f{1}{\gamma}}-1\right|\leq C|v-(v^s)^{-\mathbf{X}}|\leq C\chi.
$$
Thus, it holds
$$
\ba
\mathbf{D}_{\rm I}(t)&=\int_{\mathbb{T}^2}\int_{0}^1y_1(1-y_1)|\pa_{y_1}w|^2
\left(\f{p((v^s)^{-\mathbf{X}})}{p(v)}\right)^{1+\f{1}{\gamma}}
\f{1}{y_1(1-y_1)}\f{2\mu+\lambda}{\gamma p^{1+\f{1}{\gamma}}((v^s)^{-\mathbf{X}})}
\Big(\f{dy_1}{d\xi_1}\Big)dy_1dy'\\
&\geq (1-C\chi)\left(\f{\d p''(v_-)}{2\s_-(p'(v_-))^2}-C\d^2\right)
\int_{\mathbb{T}^2}\int_0^1y_1(1-y_1)|\pa_{y_1}w|^2dy_1dy'.
\ea
$$
Since
$$
\f{p''(v_-)}{2\s_-(p'(v_-))^2}=\f{\gamma+1}{2\gamma\s_-p(v_-)}=\a_-,
$$
we have
$$
\mathbf{D}_{\rm I}(t)\geq\a_-\d(1-C(\d+\chi))\int_{\mathbb{T}^2}\int_0^1y_1(1-y_1)|\pa_{y_1}w|^2dy_1dy'.
$$
We can deduce from \eqref{coeffi-1} that
$$
y_1(1-y_1)\f{d\xi_1}{dy_1}\geq\f{2\mu+\lambda}{\gamma p^{1+\f{1}{\gamma}}((v^s)^{-\mathbf{X}})}\f{1}{\a_-\d+C\d^2}
\geq \f{2\mu+\lambda}{2\a_-\d|p'(v_-)|}.
$$
Hence, we have
$$
\ba
\mathbf{D}_{\rm II}(t)&=(2\mu+\lambda)\int_{\mathbb{T}^2}\int_{0}^1\f{|\na_{y'}w|^2}{y_1(1-y_1)}
\left(\f{p(v_-)}{p(v)}\right)^{1+\f{1}{\gamma}}
\f{y_1(1-y_1)}{\gamma p^{1+\f{1}{\gamma}}(v_-)}
\Big(\f{d\xi_1}{dy_1}\Big)dy_1dy'\\
&\geq (1-C(\d+\chi))(2\mu+\lambda)\int_{\mathbb{T}^2}\int_{0}^1\f{|\na_{y'}w|^2}{y_1(1-y_1)}
\f{2\mu+\lambda}{|p'(v_-)|^2}\f{1}{2\a_-\d}dy_1dy'\\
&\geq (1-C(\d+\chi))\f{\s_-}{\d}\f{(2\mu+\lambda)^2}{p''(v_-)}
\int_{\mathbb{T}^2}\int_0^1\f{|\na_{y'}w|^2}{y_1(1-y_1)}dy_1dy'.
\ea
$$
Combining the estimates on $\mathbf{D}_{\rm I}(t)$ and $\mathbf{D}_{\rm II}(t)$, we have
\be\label{es-D}
\ba
\mathbf{D}(t)&\geq \a_-\d(1-C(\d+\chi))\int_{\mathbb{T}^2}\int_0^1y_1(1-y_1)|\pa_{y_1}w|^2dy_1dy'\\
&+(1-C(\d+\chi))\f{\s_-}{\d}\f{(2\mu+\lambda)^2}{p''(v_-)}
\int_{\mathbb{T}^2}\int_0^1\f{|\na_{y'}w|^2}{y_1(1-y_1)}dy_1dy'.
\ea
\ee


$\bullet$ \textbf{Proof of Lemma \ref{es-good}}: First, by \eqref{B2}, \eqref{B3}, \eqref{es-B1+G21}, \eqref{G4} and \eqref{es-D}, it holds
$$
\ba
&\quad\mathbf{B}_1(t)+\mathbf{B}_2(t)+\mathbf{B}_3(t)-\mathbf{G}_{1,1}(t)-\mathbf{G}_4(t)-\f34\mathbf{D}(t)\\
&\leq\d\a_-\Big(1+C(\nu+\d+\chi)+\f{2(1+\k)}{\gamma+1}(1+C(\nu+\d+\chi))-\f{2}{\gamma+1}(1-C(\d+\chi))\Big)
\int_{\mathbb{T}^2}\int_0^1 w^2 dy_1dy'\\
&-\f34\d\a_-(1-C(\d+\chi))\int_{\mathbb{T}^2}\int_0^1y_1(1-y_1)|\pa_{y_1}w|^2dy_1dy'
+C\big(1+\f{1}{\k}\big)\f{\d}{\nu}\mathbf{G}_3(t)\\
&-\f34(1-C(\d+\chi))\f{\s_-}{\d}\f{(2\mu+\lambda)^2}{p''(v_-)}
\int_{\mathbb{T}^2}\int_0^1\f{|\na_{y'}w|^2}{y_1(1-y_1)}dy_1dy'.
\ea
$$
Choosing $\nu$, $\d$, $\chi$ and $\k$ suitably small such that, we have
$$
\ba
&\quad\mathbf{B}_1(t)+\mathbf{B}_2(t)+\mathbf{B}_3(t)-\mathbf{G}_{1,1}(t)-\mathbf{G}_4(t)-\f34\mathbf{D}(t)\\
&\leq \f65\d\a_-\int_{\mathbb{T}^2}\int_0^1 w^2 dy_1dy'
-\f58\d\a_-\int_{\mathbb{T}^2}\int_0^1y_1(1-y_1)|\pa_{y_1}w|^2dy_1dy'\\
&-\f58\f{\s_-}{\d}\f{(2\mu+\lambda)^2}{p''(v_-)}
\int_{\mathbb{T}^2}\int_0^1\f{|\na_{y'}w|^2}{y_1(1-y_1)}dy_1dy'+C\f{\d}{\nu}\mathbf{G}_3(t).
\ea
$$
Using \eqref{poin-2} and the fact that for $\bar w:=\int_{\mathbb{T}^2}\int_0^1 w dy_1dy'$, it holds
$$
\int_{\mathbb{T}^2}\int_0^1|w-\bar w|^2dy_1dy'=\int_{\mathbb{T}^2}\int_0^1w^2dy_1dy'
-\bar w^2.
$$
By \eqref{poin-2}, we have
\be\label{B1+...}
\ba
&\quad\mathbf{B}_1(t)+\mathbf{B}_2(t)+\mathbf{B}_3(t)-\mathbf{G}_{1,1}(t)-\mathbf{G}_4(t)-\f34\mathbf{D}(t)\\
&\leq \f65\d\a_-\int_{\mathbb{T}^2}\int_0^1 w^2 dy_1dy'
-\f54\d\a_-\int_{\mathbb{T}^2}\int_0^1|w-\bar w|^2dy_1dy'\\
&\quad-\Big(\f{5\s_-}{8\d}\f{(2\mu+\lambda)^2}{p''(v_-)}-\f{5\d\a_-}{64\pi^2}\Big)
\int_{\mathbb{T}^2}\int_0^1\f{|\na_{y'}w|^2}{y_1(1-y_1)}dy_1dy'+C\f{\d}{\nu}\mathbf{G}_3(t)\\
&=-\f{\d\a_-}{20}\int_{\mathbb{T}^2}\int_0^1 w^2 dy_1dy'
+\f54\d\a_-\left(\int_{\mathbb{T}^2}\int_0^1 \,w\, dy_1dy'\right)^2\\
&\quad-\f58\Big(\f{\s_-}{\d}\f{(2\mu+\lambda)^2}{p''(v_-)}-\f{\d\a_-}{8\pi^2}\Big)
\int_{\mathbb{T}^2}\int_0^1\f{|\na_{y'}w|^2}{y_1(1-y_1)}dy_1dy'+C\f{\d}{\nu}\mathbf{G}_3(t).
\ea
\ee
Choosing $M=\f54\a_-\s_-^4v_-^2$, and combining \eqref{es-X}, \eqref{G1} and \eqref{B1+...}, we obtain
\begin{align}  \nm
&-\f{\d}{2M}|\dot{\mathbf{X}}(t)|^2+\mathbf{B}_1(t)+\mathbf{B}_2(t)+\mathbf{B}_3(t)-\mathbf{G}_{1}(t)-\mathbf{G}_4(t)-\f34\mathbf{D}(t)\\  \nm
&\leq
\Big(-\f{\d\a_-}{20}+C\d(\nu+\d+\chi)^2\Big)\int_{\mathbb{T}^2}\int_0^1 w^2 dy_1dy'\\  \nm
&-\f58\Big(\f{\s_-}{\d}\f{(2\mu+\lambda)^2}{p''(v_-)}-\f{\d\a_-}{8\pi^2}\Big)
\int_{\mathbb{T}^2}\int_0^1\f{|\na_{y'}w|^2}{y_1(1-y_1)}dy_1dy'\\  \nm
&+\s_* \int_{\mathbb{T}^2}\int_{\mathbb{R}}a^{-\mathbf{X}}_{_{\xi_1}}
\f{1+\gamma}{3\gamma^2}\f{|p(v)-p((v^s)^{-\mathbf{X}})|^3}{p^{2+\f{1}{\gamma}}((v^s)^{-\mathbf{X}})}d\xi_1d\xi'
+C\f{\d}{\nu}\mathbf{G}_3(t),
\end{align}
which indicates the desired inequality \eqref{es-key} by using \eqref{small}. The proof of Lemma \ref{es-good} is completed.

\hfill $\Box$

%
%
%
%

\begin{lemma} \label{le-basic-new}
Under the hypotheses of Proposition \ref{priori}, there exists constant $C>0$ independent of $\nu$, $\d$, $\chi$ and $T$,
such that for all $t\in[0,T]$, it holds
\begin{equation}\label{basic}
\ba
&\int_{\mathbb{T}^2}\int_{\mathbb{R}}
\r \Big(Q(v|(v^s)^{-\mathbf{X}})+\f{|h-(h^s)^{-\mathbf{X}}|^2}{2}\Big)d\xi_1d\xi'\\
&\quad
+\d\int_0^t|\dot{\mathbf{X}}(\tau)|^2d\tau+\int_0^t(G_2(\tau)+G_3(\tau)+G^s(\tau)+D(\tau))d\tau\\
&\leq C\int_{\mathbb{T}^2}\int_{\mathbb{R}}
 \Big(Q(v_0|v^s)+\f{|h_0-h^s|^2}{2}\Big)d\xi_1d\xi'
+C(\d+\chi)\int_0^t\|\na_{_\xi}(u-(u^s)^{-\mathbf{X}})\|_{H^1}^2d\tau,
\ea
\end{equation}
where $h_0(\xi)=u_0(\xi)-(2\mu+\lambda)\na_{_\xi}v_0(\xi)$ and
\be\label{new-G...}
\ba
G_2(t)&:=\f{\nu}{\d}\int_{\mathbb{T}^2}\int_{\mathbb{R}}|(v^s)^{-\mathbf{X}}_{_{\xi_1}}|(h_2^2+h_3^2)d\xi_1d\xi',\\
G_3(t)&:=\f{\nu}{\d}\int_{\mathbb{T}^2}\int_{\mathbb{R}}|(v^s)^{-\mathbf{X}}_{_{\xi_1}}|\Big|h_1-(h^s_1)^{-\mathbf{X}}-\f{p(v)-p((v^s)^{-\mathbf{X}})}{\s_* }\Big|^2d\xi_1d\xi',\\
G^s(t)&:=\int_{\mathbb{T}^2}\int_{\mathbb{R}}|(v^s)^{-\mathbf{X}}_{_{\xi_1}}||p(v)-p((v^s)^{-\mathbf{X}})|^2d\xi_1d\xi',\\
D(t)&:=\int_{\mathbb{T}^2}\int_{\mathbb{R}}|\na_{_\xi}\big(p(v)-p((v^s)^{-\mathbf{X}})\big)|^2d\xi_1d\xi'.
\ea
\ee
\end{lemma}
Note that by \eqref{a>1}, \eqref{a-v} and the uniform lower and upper boundedness of the volume function $v$, we have
$$
\mathbf{G}_2(t)\sim G_2(t), \quad \mathbf{G}_3(t)\sim G_3(t), \quad \mathbf{D}(t)\sim D(t),
$$
uniform in time $t\in[0,T]$.

\begin{proof}

First of all, we use \eqref{ba-expression} to have
$$
\ba
&\f{d}{dt}\int_{\mathbb{T}^2}\int_{\mathbb{R}}
a^{-\mathbf{X}}\r \Big(Q(v|(v^s)^{-\mathbf{X}})+\f{|h-(h^s)^{-\mathbf{X}}|^2}{2}\Big)d\xi_1d\xi'\\
=&-\f{\d}{2M}|\dot{\mathbf{X}}(t)|^2+\mathbf{B}_1(t)+\mathbf{B}_2(t)+\mathbf{B}_3(t)-\mathbf{G}_{1}(t)-\mathbf{G}_4(t)-\f34\mathbf{D}(t)\\
&-\f{\d}{2M}|\dot{\mathbf{X}}(t)|^2+\dot{\mathbf{X}}(t)\sum_{i=3}^5\mathbf{Y}_i(t)+\sum_{i=4}^9\mathbf{B}_i(t)-\sum_{i=2}^3\mathbf{G}_i(t)-\f14\mathbf{D}(t).
\ea
$$
Using Lemma \ref{es-good} and Cauchy's inequality, we find that there exist positive constants $C_1$ and $C$ such that
\begin{align}  \nm
&\f{d}{dt}\int_{\mathbb{T}^2}\int_{\mathbb{R}}
a^{-\mathbf{X}}\r \Big(Q(v|(v^s)^{-\mathbf{X}})+\f{|h-(h^s)^{-\mathbf{X}}|^2}{2}\Big)d\xi_1d\xi'\\  \nm
\leq &-C_1\underbrace{\int_{\mathbb{T}^2}\int_{\mathbb{R}}(v^s)^{-\mathbf{X}}_{_{\xi_1}}|p(v)-p((v^s)^{-\mathbf{X}})|^2d\xi_1d\xi'}_{\mathbf{G}^s(t)}
+C\int_{\mathbb{T}^2}\int_{\mathbb{R}}a^{-\mathbf{X}}_{_{\xi_1}}|p(v)-p((v^s)^{-\mathbf{X}})|^3d\xi_1d\xi'\\  \nm
&-\f{\d}{4M}|\dot{\mathbf{X}}(t)|^2+\f{C}{\d}\sum_{i=3}^5|\mathbf{Y}_i(t)|^2+\sum_{i=4}^9\mathbf{B}_i(t)-\sum_{i=2}^3\mathbf{G}_i(t)-\f14\mathbf{D}(t)
+\f{1}{40}\mathbf{G}_3(t).
\end{align}
In what follows, to control the above bad terms, we will use the above good terms $\mathbf{G}_i(t)~(i=2,3)$,
$\mathbf{D}(t)$ and $\mathbf{G}^s(t)$. In the following, we control the terms on the right hand side of the above
inequality one by one. First, for simplicity, we use the notation $w=p(v)-p((v^s)^{-\mathbf{X}})$ as \eqref{w}. Using
\eqref{change-variable}, the 3D Gagliardo-Nirenberg inequality \eqref{sobolev-inequality} in strip domain and the assumption \eqref{assumption}, it holds
\be\label{cubic}
\ba
&C\int_{\mathbb{T}^2}\int_{\mathbb{R}}a^{-\mathbf{X}}_{_{\xi_1}}|p(v)-p((v^s)^{-\mathbf{X}})|^3d\xi_1d\xi'
\leq  C\f{\nu}{\d}\|w\|^2_{L^{\infty}}\int_{\mathbb{T}^2}\int_{\mathbb{R}}(v^s)_{_{\xi_1}}^{-\mathbf{X}}|w|d\xi_1d\xi'\\
\leq & C\f{\nu}{\d}(\|w\|\|\pa_{_{\xi_1}}w\|+\|\na_{_\xi}w\|\|\na^2_{_\xi}w\|)
\Big(\int_{\mathbb{T}^2}\int_{\mathbb{R}}(v^s)_{_{\xi_1}}^{-\mathbf{X}}w^2d\xi_1d\xi'\Big)^{\f12}
\Big(\int_{\mathbb{T}^2}\int_{\mathbb{R}}(v^s)_{_{\xi_1}}^{-\mathbf{X}}d\xi_1d\xi'\Big)^{\f12}\\
\leq & C\f{\nu}{\sqrt{\d}}(\|w\|+\|\na_{_\xi}^2w\|)\|\na_{_\xi}w\|
\Big(\int_{\mathbb{T}^2}\int_{\mathbb{R}}(v^s)_{_{\xi_1}}^{-\mathbf{X}}w^2d\xi_1d\xi'\Big)^{\f12}
\leq C\f{\nu}{\sqrt{\d}}\chi\|\na_{_\xi}w\|\sqrt{\mathbf{G}^s(t)}\\
\leq &\f{1}{80}(\mathbf{D}(t)+C_1\mathbf{G}^s(t)).
\ea
\ee


$\bullet$ \textbf{Estimates on the terms $\mathbf{Y}_i(t)~(i=3,4,5)$}: Since
$$
|\mathbf{Y}_3(t)|\leq C\f{\d}{\nu}\int_{\mathbb{T}^2}\int_{\mathbb{R}}a^{-\mathbf{X}}_{_{\xi_1}}
\left|h_1-(h^s_1)^{-\mathbf{X}}-\f{p(v)-p((v^s)^{-\mathbf{X}})}{\s_* }\right|d\xi_1d\xi'
\leq C\f{\d}{\sqrt{\nu}}\sqrt{\mathbf{G}_3(t)},
$$
which leads to
$$
\f{C}{\d}|\mathbf{Y}_3(t)|^2\leq C\f{\d}{\nu}\mathbf{G}_3(t)\leq \f{1}{40}\mathbf{G}_3(t).
$$
For $\mathbf{Y}_4(t)$, we first use the definition of $h$  and $h^s=(h^s_1,0,0)^t$ to estimate $h-(h^s)^{-\mathbf{X}}$
in terms of $u-(u^s)^{-\mathbf{X}}$ and $v-(v^s)^{-\mathbf{X}}$ as follow.
\be\label{h-h^s}
\ba
h-(h^s)^{-\mathbf{X}}=u-(u^s)^{-\mathbf{X}}-(2\mu+\lambda)\na_{_\xi}(v-(v^s)^{-\mathbf{X}}),
\ea
\ee
which together with the assumption \eqref{assumption} implies
$$
\|h-(h^s)^{-\mathbf{X}}\|\leq C(\|u-(u^s)^{-\mathbf{X}}\|+\|\na_{_\xi}(v-(v^s)^{-\mathbf{X}})\|)\leq C\chi.
$$
This together with $\|a^{-\mathbf{X}}_{_{\xi_1}}\|_{L^{\infty}}\leq C\nu\d$ and the assumption \eqref{assumption} lead to
$$
|\mathbf{Y}_4(t)|\leq C\sqrt{\mathbf{G}_3(t)}\|a^{-\mathbf{X}}_{_{\xi_1}}\|^{\f12}_{L^{\infty}}
(\|u_1-(u^s_1)^{-\mathbf{X}}\|+\|v-(v^s)^{-\mathbf{X}}\|_{H^1})\leq C\chi(\nu\d)^{\f12}\sqrt{\mathbf{G}_3(t)},
$$
which implies
$$
\f{C}{\d}|\mathbf{Y}_4(t)|^2\leq C\chi^2\nu\mathbf{G}_3(t)\leq \f{1}{40}\mathbf{G}_3(t).
$$
Using \eqref{use-3} and assumption \eqref{assumption}, we have
$$
\ba
\f{C}{\d}|\mathbf{Y}_5(t)|^2&\leq \f{C}{\d}\left(\int_{\mathbb{T}^2}\int_{\mathbb{R}}a^{-\mathbf{X}}_{_{\xi_1}}w^2d\xi_1d\xi'\right)^2
+\f{C}{\d}\left(\int_{\mathbb{T}^2}\int_{\mathbb{R}}a^{-\mathbf{X}}_{_{\xi_1}}\f{h_2^2+h_3^2}{2}d\xi_1d\xi'\right)^2\\
&\leq C\f{\nu^2}{\d^3}\left(\int_{\mathbb{T}^2}\int_{\mathbb{R}}(v^s)^{-\mathbf{X}}_{_{\xi_1}}w^2d\xi_1d\xi'\right)^2
+\f{C}{\d}\left(\int_{\mathbb{T}^2}\int_{\mathbb{R}}a^{-\mathbf{X}}_{_{\xi_1}}\f{h_2^2+h_3^2}{2}d\xi_1d\xi'\right)^2\\
&\leq C\f{\nu^2}{\d}\|w\|^2\int_{\mathbb{T}^2}\int_{\mathbb{R}}(v^s)^{-\mathbf{X}}_{_{\xi_1}}w^2d\xi_1d\xi'
+C\nu\|h-(h^s)^{-\mathbf{X}}\|^2\int_{\mathbb{T}^2}\int_{\mathbb{R}}a^{-\mathbf{X}}_{_{\xi_1}}\f{h_2^2+h_3^2}{2}d\xi_1d\xi'\\
&\leq C\chi^2\mathbf{G}^s(t)+C\nu\chi^2\mathbf{G}_2(t)\leq \f{1}{40}(C_1\mathbf{G}^s(t)+\mathbf{G}_2(t)).
\ea
$$


$\bullet$ \textbf{Estimates on the terms $\mathbf{B}_i(t)$ $(i=4,\cdots,9)$ }:
Recalling the definition of $F$ and the assumption \eqref{assumption}, it holds
\be\label{F-small}
\|F\|_{L^{\infty}}\leq C(\|v-(v^s)^{-\mathbf{X}}\|_{L^{\infty}}+\|u_1-(u^s_1)^{-\mathbf{X}}\|_{L^{\infty}})\leq C\chi.
\ee
For $\mathbf{B}_4(t)$,  using \eqref{use-3} and \eqref{F-small}, it holds
$$
\ba
\mathbf{B}_4(t)&=\int_{\mathbb{T}^2}\int_{\mathbb{R}}F a^{-\mathbf{X}}_{_{\xi_1}}\Big(
Q(v|(v^s)^{-\mathbf{X}})+\f{(h_1-(h^s_1)^{-\mathbf{X}})^2}{2}+\f{h_2^2+h_3^2}{2}\Big)d\xi_1d\xi'\\
&\leq C\int_{\mathbb{T}^2}\int_{\mathbb{R}}|F| a^{-\mathbf{X}}_{_{\xi_1}}\big(
Q(v|(v^s)^{-\mathbf{X}})+|p(v)-p((v^s)^{-\mathbf{X}})|^2\big)d\xi_1d\xi'\\
&+C\int_{\mathbb{T}^2}\int_{\mathbb{R}}|F| a^{-\mathbf{X}}_{_{\xi_1}}\Big(\f{h_2^2+h_3^2}{2}
+\Big|h_1-(h^s_1)^{-\mathbf{X}}-\f{p(v)-p((v^s)^{-\mathbf{X}})}{\s_* }\Big|^2\Big)d\xi_1d\xi'\\
&\leq C\int_{\mathbb{T}^2}\int_{\mathbb{R}}|F| a^{-\mathbf{X}}_{_{\xi_1}}w^2d\xi_1d\xi'
+C\chi(\mathbf{G}_2(t)+\mathbf{G}_3(t)).
\ea
$$
By the definition of $F$ \eqref{F}, we have
$$
\ba
C\int_{\mathbb{T}^2}\int_{\mathbb{R}}|F| a^{-\mathbf{X}}_{_{\xi_1}}w^2d\xi_1d\xi'
&\leq \underbrace{C\int_{\mathbb{T}^2}\int_{\mathbb{R}}a^{-\mathbf{X}}_{_{\xi_1}}|w|^3d\xi_1d\xi'}_{\mathbf{B}_{4,1}(t)}
+\underbrace{C\int_{\mathbb{T}^2}\int_{\mathbb{R}}a^{-\mathbf{X}}_{_{\xi_1}}|h_1-(h^s_1)^{-\mathbf{X}}|w^2d\xi_1d\xi'}_{\mathbf{B}_{4,2}(t)}\\
&+\underbrace{C\int_{\mathbb{T}^2}\int_{\mathbb{R}}a^{-\mathbf{X}}_{_{\xi_1}}|\pa_{_{\xi_1}}(v-(v^s)^{-\mathbf{X}})|w^2d\xi_1d\xi'}_{\mathbf{B}_{4,3}(t)}.
\ea
$$
The same as \eqref{cubic}, it holds
$$
\mathbf{B}_{4,1}(t)\leq \f{1}{80}(\mathbf{D}(t)+C_1\mathbf{G}^s(t)).
$$
Using the Sobolev inequality in strip domain and Lemma \ref{le-shock} , we have
$$
\ba
\mathbf{B}_{4,2}(t)&\leq C\int_{\mathbb{T}^2}\int_{\mathbb{R}}a^{-\mathbf{X}}_{_{\xi_1}}\Big|h_1-(h^s_1)^{-\mathbf{X}}-\f{p(v)-p((v^s)^{-\mathbf{X}})}{\s_* }\Big|w^2d\xi_1d\xi'
+C\int_{\mathbb{T}^2}\int_{\mathbb{R}}a^{-\mathbf{X}}_{_{\xi_1}}|w|^3d\xi_1d\xi'\\
&\leq C\|w\|_{L^{\infty}}^2\sqrt{\mathbf{G}_3(t)}\sqrt{\f{\nu}{\d}}\Big(\int_{\mathbb{T}^2}\int_{\mathbb{R}}(v^s)^{-\mathbf{X}}_{_{\xi_1}} d\xi_1d\xi'\Big)^{\f12}
+\f{1}{80}(\mathbf{D}(t)+C_1\mathbf{G}^s(t))\\
&\leq C(\|w\|+\|\na_{_\xi}^2w\|)\|\na_{_\xi}w\|\sqrt{\mathbf{G}_3(t)}\sqrt{\nu}+\f{1}{80}(\mathbf{D}(t)+C_1\mathbf{G}^s(t))\\
&\leq C\chi\sqrt{\nu}\sqrt{\mathbf{D}(t)}\sqrt{\mathbf{G}_3(t)}+\f{1}{80}(\mathbf{D}(t)+C_1\mathbf{G}^s(t))
\leq \f{1}{40}\mathbf{D}(t)+\f{1}{80}(\mathbf{G}_3(t)+C_1\mathbf{G}^s(t)).
\ea
$$
Using the fact that
$$
\na_{_\xi}\big(p(v)-p((v^s)^{-\mathbf{X}})\big)=p'(v)\na_{_\xi}(v-(v^s)^{-\mathbf{X}})+\na_{_\xi}(v^s)^{-\mathbf{X}}\big(p'(v)-p'((v^s)^{-\mathbf{X}})\big),
$$
hence, using \eqref{use-4} and \eqref{use-1}, we have
\be\label{daoshu}
\ba
|\na_{_\xi}(v-(v^s)^{-\mathbf{X}})|&\leq C|\na_{_\xi}\big(p(v)-p((v^s)^{-\mathbf{X}})\big)|+C|(v^s)^{-\mathbf{X}}_{_{\xi_1}}||v-(v^s)^{-\mathbf{X}}|\\
&\leq C|\na_{_\xi}\big(p(v)-p((v^s)^{-\mathbf{X}})\big)|+C|(v^s)^{-\mathbf{X}}_{_{\xi_1}}||p(v)-p((v^s)^{-\mathbf{X}})|.
\ea
\ee
Using \eqref{daoshu}, we have
$$
\ba
\mathbf{B}_{4,3}(t)&\leq C\int_{\mathbb{T}^2}\int_{\mathbb{R}}a^{-\mathbf{X}}_{_{\xi_1}}\big(|\pa_{_{\xi_1}}w|+(v^s)^{-\mathbf{X}}_{_{\xi_1}}|w|\big)w^2d\xi_1d\xi'\\
&\leq \|w\|_{L^{\infty}}^2\|\pa_{_{\xi_1}}w\|\|a^{-\mathbf{X}}_{_{\xi_1}}\|+C\|w\|_{L^{\infty}}\mathbf{G}^s(t)\\
&\leq C(\|w\|+\|\na_{_\xi}^2w\|)\|\na_{_\xi}w\|^2\sqrt{\d}\nu+C\chi\mathbf{G}^s(t)\\
&\leq C\chi\sqrt{\d}\nu\mathbf{D}(t)+C\chi\mathbf{G}^s(t)
\leq \f{1}{80}(\mathbf{D}(t)+C_1\mathbf{G}^s(t)).
\ea
$$
Combining the above estimates, we have
$$
\mathbf{B}_4(t)\leq \f{1}{16}\mathbf{D}(t)+\f{1}{20}(\mathbf{G}_2(t)+\mathbf{G}_3(t)+C_1\mathbf{G}^s(t)).
$$
Notice that
$$
\pa_{_{\xi_1}}(v-(v^s)^{-\mathbf{X}})=\f{\pa_{_{\xi_1}}\big(p(v)-p((v^s)^{-\mathbf{X}})\big)}{p'(v)}
+\pa_{_{\xi_1}}p((v^s)^{-\mathbf{X}})\left(\f{1}{p'(v)}-\f{1}{p'((v^s)^{-\mathbf{X}})}\right).
$$
Using \eqref{a-v}, it holds
$$
\ba
\mathbf{B}_5(t)&\leq C\int_{\mathbb{T}^2}\int_{\mathbb{R}}(v^s)^{-\mathbf{X}}_{_{\xi_1}}|\pa_{_{\xi_1}}w|
\left(|w|+\left|h_1-(h^s_1)^{-\mathbf{X}}-\f{p(v)-p((v^s)^{-\mathbf{X}})}{\s_* }\right|\right)d\xi_1d\xi'\\
&\quad+C\int_{\mathbb{T}^2}\int_{\mathbb{R}}|(v^s)^{-\mathbf{X}}_{_{\xi_1}}|^2|w|
\left(|w|+\left|h_1-(h^s_1)^{-\mathbf{X}}-\f{p(v)-p((v^s)^{-\mathbf{X}})}{\s_* }\right|\right)d\xi_1d\xi'\\
&\leq C\d\|\pa_{_{\xi_1}}w\|\sqrt{\mathbf{G}^s(t)}+C\d\|\pa_{_{\xi_1}}w\|\sqrt{\f{\d}{\nu}}\sqrt{\mathbf{G}_3(t)}
+C\d^2\mathbf{G}^s(t)+C\d^2\sqrt{\mathbf{G}^s(t)}\sqrt{\f{\d}{\nu}}\sqrt{\mathbf{G}_3(t)}\\
&\leq \f{1}{80}(\mathbf{D}(t)+C_1\mathbf{G}^s(t)+\mathbf{G}_3(t)).
\ea
$$
Similarly, we have
$$
\ba
\mathbf{B}_6(t)+\mathbf{B}_7(t)+\mathbf{B}_8(t)
&\leq  C\int_{\mathbb{T}^2}\int_{\mathbb{R}}|(v^s)^{-\mathbf{X}}_{_{\xi_1}}||\pa_{_{\xi_1}}w||w|d\xi_1d\xi'
+C\f{\nu}{\d}\int_{\mathbb{T}^2}\int_{\mathbb{R}}|(v^s)^{-\mathbf{X}}_{_{\xi_1}}||\pa_{_{\xi_1}}w||w|d\xi_1d\xi'\\
&\quad+C\f{\nu}{\d}\int_{\mathbb{T}^2}\int_{\mathbb{R}}|(v^s)^{-\mathbf{X}}_{_{\xi_1}}|^2w^2d\xi_1d\xi'\\
&\leq C\d\|\pa_{_{\xi_1}}w\|\sqrt{\mathbf{G}^s(t)}+C\nu\|\pa_{_{\xi_1}}w\|\sqrt{\mathbf{G}^s(t)}
+C\nu\d\mathbf{G}^s(t)\\
&\leq \f{1}{80}(\mathbf{D}(t)+C_1\mathbf{G}^s(t)).
\ea
$$
Finally, we control the last term $\mathbf{B}_9(t)$. Using the definition of $R$ in \eqref{R}, it holds
$$
\ba
\mathbf{B}_9(t)&=(2\mu+\lambda)\int_{\mathbb{T}^2}\int_{\mathbb{R}}a^{-\mathbf{X}}(h-(h^s)^{-\mathbf{X}})\cdot\f{\na_{_\xi} u\cdot\na_{_\xi} v-\div_{_\xi} u\na_{_\xi}v}{v}d\xi_1d\xi'\\
&-\mu\int_{\mathbb{T}^2}\int_{\mathbb{R}}a^{-\mathbf{X}}(h-(h^s)^{-\mathbf{X}})\cdot\na_{_\xi}\times\na_{_\xi}\times u\,d\xi_1d\xi'
=:\mathbf{B}_{9,1}(t)+\mathbf{B}_{9,2}(t).
\ea
$$
To control $\mathbf{B}_{9,1}(t)$, we set $u'=(u_2,u_3)$, $\na_{_{\xi'}}=(\pa_{_{_{\xi_2}}},\pa_{_{_{\xi_3}}})$ and
$\na_{_{\xi'}}\cdot u'=\pa_{_{_{\xi_2}}}u_2+\pa_{_{_{\xi_3}}}u_3$. Notice that the first component of
$\na_{_\xi}u\cdot\na_{_\xi}v-\div_{_\xi}u\na_{_\xi}v$ is
$$
\ba
&\pa_{_{\xi_1}}u\cdot\na_{_\xi}v-\div_{_\xi}u\pa_{_{\xi_1}}v=\pa_{_{\xi_1}}u'\cdot\na_{_{\xi'}}v-\na_{_{\xi'}}\cdot u'\pa_{_{\xi_1}}v\\
=&\pa_{_{\xi_1}}u'\cdot\f{\na_{_{\xi'}}p(v)}{-\gamma p^{1+\f{1}{\gamma}}(v)}
-\na_{_{\xi'}}\cdot u'\f{\pa_{_{\xi_1}}\big(p(v)-p((v^s)^{-\mathbf{X}})\big)}{-\gamma p^{1+\f{1}{\gamma}}(v)}
-\na_{_{\xi'}}\cdot u'\f{\pa_{_{\xi_1}}p((v^s)^{-\mathbf{X}})}{-\gamma p^{1+\f{1}{\gamma}}(v)}.
\ea
$$
By assumption \eqref{assumption} and Sobolev inequality, we have
\be\label{h-L3}
\|h_1-(h^s_1)^{-\mathbf{X}}\|_{L^3}\leq C\|h_1-(h^s_1)^{-\mathbf{X}}\|_{H^1}\leq C(\|u-(u^s)^{-\mathbf{X}}\|_{H^1}
+\|v-(v^s)^{-\mathbf{X}}\|_{H^2})\leq C\chi.
\ee
Thus, using \eqref{h-L3}, the first part of the integrand in $\mathbf{B}_{9,1}(t)$ can be controlled as
$$
\ba
&(2\mu+\lambda)\int_{\mathbb{T}^2}\int_{\mathbb{R}}\f{a^{-\mathbf{X}}}{v}(h_1-(h^s_1)^{-\mathbf{X}})
(\pa_{_{\xi_1}}u\cdot\na_{_\xi}v-\div_{_\xi}u\pa_{_{\xi_1}}v)d\xi_1d\xi'\\
\leq & C\|h_1-(h^s_1)^{-\mathbf{X}}\|_{L^3}\|\na_{_\xi}(u-(u^s)^{-\mathbf{X}})\|_{L^6}\sqrt{\mathbf{D}(t)}
+C\d\Big(\sqrt{\f{\d}{\nu}}\sqrt{\mathbf{G}_3(t)}+\sqrt{\mathbf{G}^s(t)}\Big)\|\na_{_{\xi'}}(u-(u^{s})^{-\mathbf{X}})\|\\
\leq &C(\chi+\d)(\mathbf{G}_3(t)+C_1\mathbf{G}^s(t)+\mathbf{D}(t)+\|\na_{_\xi}(u-(u^s)^{-\mathbf{X}})\|^2_{H^1})\\
\leq &\f{1}{80}(\mathbf{G}_3(t)+C_1\mathbf{G}^s(t)+\mathbf{D}(t))+C(\chi+\d)\|\na_{_\xi}(u-(u^s)^{-\mathbf{X}})\|^2_{H^1}.
\ea
$$
Similarly, the second and third parts of integrand in $\mathbf{B}_{9,1}(t)$ can be treated as
$$
\ba
&(2\mu+\lambda)\int_{\mathbb{T}^2}\int_{\mathbb{R}}\f{a^{-\mathbf{X}}}{v}h'\cdot
(\na_{_{\xi'}}u\cdot\na_{_\xi}v-\div_{_\xi}u\na_{_{\xi'}}v)d\xi_1d\xi'\\
\leq & C\|h'\|_{L^3}\|\na_{_\xi}(u-(u^s)^{-\mathbf{X}})\|_{L^6}\sqrt{\mathbf{D}(t)}
+C\d\sqrt{\f{\d}{\nu}}\sqrt{\mathbf{G}_2(t)}(\|\na_{_{\xi'}}(u-(u^{s})^{-\mathbf{X}})\|+\sqrt{\mathbf{D}(t)})\\
\leq &C(\chi+\d)(\mathbf{G}_2(t)+\mathbf{D}(t)+\|\na_{_\xi}(u-(u^s)^{-\mathbf{X}})\|^2_{H^1})\\
\leq &\f{1}{80}(\mathbf{G}_2(t)+\mathbf{D}(t))+C(\chi+\d)\|\na_{_\xi}(u-(u^s)^{-\mathbf{X}})\|^2_{H^1}.
\ea
$$
Thus, we have
\be\label{es-R1}
\mathbf{B}_{9,1}(t)\leq \f{1}{40}(\mathbf{G}_2(t)+\mathbf{G}_3(t)+C_1\mathbf{G}^s(t)+\mathbf{D}(t))
+C(\chi+\d)\|\na_{_\xi}(u-(u^s)^{-\mathbf{X}})\|^2_{H^1}.
\ee

Using \eqref{h-h^s}, and direct calculations yield
$$
\ba
\mathbf{B}_{9,2}(t)&=-\mu\int_{\mathbb{T}^2}\int_{\mathbb{R}}a^{-\mathbf{X}}(u-(u^s)^{-\mathbf{X}})\cdot\na_{_\xi}\times\na_{_\xi}\times u\,d\xi_1d\xi'\\
&\quad+\mu(2\mu+\lambda)\int_{\mathbb{T}^2}\int_{\mathbb{R}}a^{-\mathbf{X}}\na_{_\xi}(v-(v^s)^{-\mathbf{X}})\cdot\na_{_\xi}\times\na_{_\xi}\times u\,d\xi_1d\xi'\\
&
=:\mathbf{B}^1_{9,2}(t)+\mathbf{B}^2_{9,2}(t).
\ea
$$
By integration by parts over $\mathbb{R}\times\mathbb{T}^2$, it holds
$$
\ba
\mathbf{B}^1_{9,2}(t)&=-\mu\int_{\mathbb{T}^2}\int_{\mathbb{R}}a^{-\mathbf{X}}|\na_{_\xi}\times(u-(u^s)^{-\mathbf{X}})|^2d\xi_1d\xi'\\
&\quad-\mu\int_{\mathbb{T}^2}\int_{\mathbb{R}}a^{-\mathbf{X}}_{_{\xi_1}}\big(u_2(\pa_{_{\xi_1}}u_2-\pa_{_{\xi_2}}u_1)
-u_3(\pa_{_{\xi_3}}u_1-\pa_{_{\xi_1}}u_3)\big)d\xi_1d\xi'\\
&\leq -\f{3\mu}{4}\int_{\mathbb{T}^2}\int_{\mathbb{R}}a^{-\mathbf{X}}|\na_{_\xi}\times(u-(u^s)^{-\mathbf{X}})|^2d\xi_1d\xi'\\
&\quad+C\int_{\mathbb{T}^2}\int_{\mathbb{R}}|a^{-\mathbf{X}}_{_{\xi_1}}|^2(h_2^2+h_3^2)d\xi_1d\xi'
+C\int_{\mathbb{T}^2}\int_{\mathbb{R}}|a^{-\mathbf{X}}_{_{\xi_1}}|^2|\na_{_{\xi'}}v|^2d\xi_1d\xi'\\
&\leq -\f{3\mu}{4}\int_{\mathbb{T}^2}\int_{\mathbb{R}}a^{-\mathbf{X}}|\na_{_\xi}\times(u-(u^s)^{-\mathbf{X}})|^2d\xi_1d\xi'
+C\nu\d\mathbf{G}_2(t)+C(\nu\d)^2\mathbf{D}(t)\\
&\leq -\f{3\mu}{4}\int_{\mathbb{T}^2}\int_{\mathbb{R}}a^{-\mathbf{X}}|\na_{_\xi}\times(u-(u^s)^{-\mathbf{X}})|^2d\xi_1d\xi'
+\f{1}{80}(\mathbf{G}_2(t)+\mathbf{D}(t)).
\ea
$$
By integration by parts over $\mathbb{R}\times\mathbb{T}^2$, we have
$$
\ba
\mathbf{B}^2_{9,2}(t)&=-\mu(2\mu+\lambda)\int_{\mathbb{T}^2}\int_{\mathbb{R}}(v-(v^s)^{-\mathbf{X}})\na_{_\xi}\times\na_{_\xi}\times u\cdot\na_{_\xi}a^{-\mathbf{X}}d\xi_1d\xi'\\
&=-\mu(2\mu+\lambda)\int_{\mathbb{T}^2}\int_{\mathbb{R}}(v-(v^s)^{-\mathbf{X}})
\big(\pa_{_{\xi_2}}(\pa_{_{\xi_1}}u_2-\pa_{_{\xi_2}}u_1)-\pa_{_{\xi_3}}(\pa_{_{\xi_3}}u_1-\pa_{_{\xi_1}}u_3)\big)a^{-\mathbf{X}}_{_{\xi_1}}d\xi_1d\xi'\\
&=\mu(2\mu+\lambda)\int_{\mathbb{T}^2}\int_{\mathbb{R}}\big(
\pa_{_{\xi_2}}(v-(v^s)^{-\mathbf{X}})(\pa_{_{\xi_1}}u_2-\pa_{_{\xi_2}}u_1)
-\pa_{_{\xi_3}}(v-(v^s)^{-\mathbf{X}})(\pa_{_{\xi_3}}u_1-\pa_{_{\xi_1}}u_3)\big)a^{-\mathbf{X}}_{_{\xi_1}}d\xi_1d\xi'\\
&\leq C\nu\d\mu\int_{\mathbb{T}^2}\int_{\mathbb{R}}a^{-\mathbf{X}}|\na_{_\xi}\times(u-(u^s)^{-\mathbf{X}})|^2d\xi_1d\xi'
+C\nu\d\int_{\mathbb{T}^2}\int_{\mathbb{R}}|\na_{_{\xi'}}(v-(v^s)^{-\mathbf{X}})|^2d\xi_1d\xi'\\
&\leq \f{\mu}{4}\int_{\mathbb{T}^2}\int_{\mathbb{R}}a^{-\mathbf{X}}|\na_{_\xi}\times(u-(u^s)^{-\mathbf{X}})|^2d\xi_1d\xi'
+\f{1}{80}\mathbf{D}(t).
\ea
$$
Combining the above two estimates, we have
\be\label{es-R2}
\mathbf{B}_{9,2}(t)\leq -\f{\mu}{2}\int_{\mathbb{T}^2}\int_{\mathbb{R}}a^{-\mathbf{X}}|\na_{_\xi}\times(u-(u^s)^{-\mathbf{X}})|^2d\xi_1d\xi'
+\f{1}{80}\mathbf{G}_2(t)+\f{1}{40}\mathbf{D}(t).
\ee
Combination \eqref{es-R1} and \eqref{es-R2} yields
$$
\mathbf{B}_9(t)\leq \f{3}{80}(\mathbf{G}_2(t)+\mathbf{G}_3(t)+C_1\mathbf{G}^s(t))+\f{1}{20}\mathbf{D}(t)
+C(\chi+\d)\|\na_{_\xi}(u-(u^s)^{-\mathbf{X}})\|^2_{H^1}.
$$


$\bullet$ \textbf{Conclusion}: Combining the above estimates, we have
$$
\ba
&\f{d}{dt}\int_{\mathbb{T}^2}\int_{\mathbb{R}}
a^{-\mathbf{X}}\r \Big(Q(v|(v^s)^{-\mathbf{X}})+\f{|h-(h^s)^{-\mathbf{X}}|^2}{2}\Big)d\xi_1d\xi'\\
\leq &-\f{\d}{4M}|\dot{\mathbf{X}}(t)|^2-\f{1}{2}\mathbf{G}_2(t)-\f{1}{2}\mathbf{G}_3(t)
-\f{1}{16}\mathbf{D}(t)-\f12 C_1\mathbf{G}^s(t)+C(\chi+\d)\|\na_{_\xi}(u-(u^s)^{-\mathbf{X}})\|^2_{H^1}.
\ea
$$
Integrating the above inequality over $[0,t]$ for any $t\in[0,T]$, using \eqref{daoshu},
we can obtain the desired inequality \eqref{basic} with the new notations \eqref{new-G...}.
\end{proof}

%
%
%
%

\subsection{Estimates for $\|u-(u^s)^{-\mathbf{X}}\|$ and $\|v-(v^s)^{-\mathbf{X}}\|_{H^1}$}

In this subsection, we shall obtain the zero-th order energy estimates for function $(v,u)$.

\begin{lemma}\label{le-1st-v}
Under the hypotheses of Proposition \ref{priori}, there exists constant $C>0$ independent of $\nu$, $\d$, $\chi$ and $T$,
such that for all $t\in[0,T]$, it holds
\be\label{es-1st}
\ba
&\|(v-(v^s)^{-\mathbf{X}})(t)\|_{H^1}^2+\|(u-(u^s)^{-\mathbf{X}})(t)\|^2+\d\int_0^t|\dot{\mathbf{X}}(\tau)|^2d\tau\\
&+\int_0^t(G^s(\tau)+D(\tau)+\|\na_{_\xi}(u-(u^s)^{-\mathbf{X}})\|^2)d\tau\\
&\leq C\big(\|v_0-v^s\|_{H^1}^2+\|u_0-u^s\|^2\big)
+C(\d+\chi)\int_0^t\|\na^2_{_\xi}(u-(u^s)^{-\mathbf{X}})\|^2d\tau,
\ea
\ee
where $G^s$ and $D$ are as in \eqref{new-G...}.
\end{lemma}
\textbf{\emph{Proof}}: From systems \eqref{ns-xi} and \eqref{stationary-1}, we can get the perturbed system for
$(v-(v^s)^{-\mathbf{X}}, u-(u^s)^{-\mathbf{X}})$ as
\be\label{perturb-vu}
\left\{
\begin{array}{l}
\r\pa_t (v-(v^s)^{-\mathbf{X}})-\s\r\pa_{_{\xi_1}}(v-(v^s)^{-\mathbf{X}})+\r u\cdot\na_{_\xi}(v-(v^s)^{-\mathbf{X}})\\[2mm]
\quad-\dot{\mathbf{X}}(t)\r\pa_{_{\xi_1}}(v^s)^{-\mathbf{X}}+F\pa_{_{\xi_1}}(v^s)^{-\mathbf{X}}=\div_{_\xi}(u-(u^s)^{-\mathbf{X}}),\\[2mm]
\r\pa_t (u-(u^s)^{-\mathbf{X}})-\s\r\pa_{_{\xi_1}}(u-(u^s)^{-\mathbf{X}})+\r u\cdot\na_{_\xi}(u-(u^s)^{-\mathbf{X}})+\na_{_\xi}\big(p(v)-p((v^s)^{-\mathbf{X}})\big)\\[2mm]
\quad-\dot{\mathbf{X}}(t)\r\pa_{_{\xi_1}}(u^s)^{-\mathbf{X}}+F\pa_{_{\xi_1}}(u^s)^{-\mathbf{X}}=\mu\Delta_{_\xi}(u-(u^s)^{-\mathbf{X}})+(\mu+\lambda)\na_{_\xi}\div_{_\xi}(u-(u^s)^{-\mathbf{X}}),
\end{array}
\right.
\ee
where $F$ is defined in \eqref{F}. Multiplying \eqref{perturb-vu}$_1$ by $-(p(v)-p((v^s)^{-\mathbf{X}}))$, using \eqref{ba-5}, it holds
\be\label{vu-1}
\ba
&\pa_t\big(\r Q(v|(v^s)^{-\mathbf{X}})\big)-\s\pa_{_{\xi_1}}\big(\r Q(v|(v^s)^{-\mathbf{X}})\big)
+\div_{_\xi}\big(\r u Q(v|(v^s)^{-\mathbf{X}})\big)\\
=&-\big(p(v)-p((v^s)^{-\mathbf{X}})\big)\div_{_\xi}(u-(u^s)^{-\mathbf{X}})
-\dot{\mathbf{X}}(t)\r p'((v^s)^{-\mathbf{X}})(v-(v^s)^{-\mathbf{X}})(v^s)^{-\mathbf{X}}_{_{\xi_1}}\\
&+\s_*  p(v|(v^s)^{-\mathbf{X}})(v^s)_{_{\xi_1}}^{-\mathbf{X}}
+F p'((v^s)^{-\mathbf{X}})(v-(v^s)^{-\mathbf{X}})(v^s)^{-\mathbf{X}}_{_{\xi_1}}.
\ea
\ee
Multiplying \eqref{vu-1}$_2$ by $u-(u^s)^{-\mathbf{X}}$, using \eqref{ba-5}, it holds
\be\label{vu-2}
\ba
&\pa_t\big(\r\f{|u-(u^s)^{-\mathbf{X}}|^2}{2}\big)-\s\pa_{_{\xi_1}}\big(\r\f{|u-(u^s)^{-\mathbf{X}}|^2}{2}\big)
+\div_{_\xi}\big(\r u\f{|u-(u^s)^{-\mathbf{X}}|^2}{2}\big)\\
&+\div_{_\xi}\big((p(v)-p((v^s)^{-\mathbf{X}}))(u-(u^s)^{-\mathbf{X}})\big)
-\div_{_\xi}(u-(u^s)^{-\mathbf{X}})\big(p(v)-p((v^s)^{-\mathbf{X}})\big)\\
=&\dot{\mathbf{X}}(t)\r (u^s_1)^{-\mathbf{X}}_{_{\xi_1}}(u_1-(u^s_1)^{-\mathbf{X}})
-F(u^s_1)^{-\mathbf{X}}_{_{\xi_1}}(u_1-(u^s_1)^{-\mathbf{X}})
-\mu|\na_{_\xi}(u-(u^s)^{-\mathbf{X}})|^2\\
&-(\mu+\lambda)(\div_{_\xi}(u-(u^s)^{-\mathbf{X}}))^2
+\mu\div_{_\xi}\big(\na_{_\xi}(u-(u^s)^{-\mathbf{X}})\cdot(u-(u^s)^{-\mathbf{X}})\big)\\
&+(\mu+\lambda)\div_{_\xi}\big(\div_{_\xi}(u-(u^s)^{-\mathbf{X}})(u-(u^s)^{-\mathbf{X}})\big).
\ea
\ee
Adding \eqref{vu-1} and \eqref{vu-2} together, and integrating the resultant equation over
$\mathbb{R}\times\mathbb{T}^2$, we have
\be\label{ba-expre-2}
\ba
&\quad\f{d}{dt}\int_{\mathbb{T}^2}\int_{\mathbb{R}}
\r \Big(Q(v|(v^s)^{-\mathbf{X}})+\f{|u-(u^s)^{-\mathbf{X}}|^2}{2}\Big)d\xi_1d\xi'\\
&+\underbrace{\int_{\mathbb{T}^2}\int_{\mathbb{R}}\big(\mu|\na_{_\xi}(u-(u^s)^{-\mathbf{X}})|^2
+(\mu+\lambda)(\div_{_\xi}(u-(u^s)^{-\mathbf{X}}))^2\big)d\xi_1d\xi'}_{\mathbf{D}_1(t)}\\
&=\dot{\mathbf{X}}(t)\mathcal{Y}(t)+\sum_{i=1}^3I_i(t),
\ea
\ee
where
$$
\ba
\mathcal{Y}(t)=
&-\int_{\mathbb{T}^2}\int_{\mathbb{R}}\r p'((v^s)^{-\mathbf{X}})(v-(v^s)^{-\mathbf{X}})(v^s)^{-\mathbf{X}}_{_{\xi_1}}d\xi_1d\xi'\\
&+\int_{\mathbb{T}^2}\int_{\mathbb{R}}\r (u^s_1)^{-\mathbf{X}}_{_{\xi_1}}(u_1-(u^s_1)^{-\mathbf{X}})d\xi_1d\xi'=:\mathcal{Y}_1(t)+\mathcal{Y}_2(t),
\ea
$$
and
$$
\ba
I_1(t)&=\s_* \int_{\mathbb{T}^2}\int_{\mathbb{R}}p(v|(v^s)^{-\mathbf{X}})(v^s)^{-\mathbf{X}}_{_{\xi_1}}d\xi_1d\xi',\\
I_2(t)&=\int_{\mathbb{T}^2}\int_{\mathbb{R}}Fp'((v^s)^{-\mathbf{X}})(v-(v^s)^{-\mathbf{X}})(v^s)^{-\mathbf{X}}_{_{\xi_1}}d\xi_1d\xi',\\
I_3(t)&=-\int_{\mathbb{T}^2}\int_{\mathbb{R}}F(u^s_1)^{-\mathbf{X}}_{_{\xi_1}}(u_1-(u^s_1)^{-\mathbf{X}})d\xi_1d\xi'.
\ea
$$
Since \eqref{use-4} and \eqref{use-1}, it holds
$$
|\mathcal{Y}_1(t)|\leq\Big(\int_{\mathbb{T}^2}\int_{\mathbb{R}}(v^s)^{-\mathbf{X}}_{_{\xi_1}}d\xi_1d\xi'\Big)^{\f12}
\Big(\int_{\mathbb{T}^2}\int_{\mathbb{R}}(v^s)^{-\mathbf{X}}_{_{\xi_1}}|v-(v^s)^{-\mathbf{X}}|^2d\xi_1d\xi'\Big)^{\f12}
\leq C\sqrt{\d}\sqrt{G^s(t)},
$$
hence, we have
$$
|\dot{\mathbf{X}} (t)||\mathcal{Y}_1(t)|\leq \f{\d}{8}|\dot{\mathbf{X}}(t)|^2+CG^s(t).
$$

To control $\mathcal{Y}_2(t)$, by \eqref{h-h^s}, it holds
\be\label{u1-u1^s}
u_1-(u^s_1)^{-\mathbf{X}}=h_1-(h^s_1)^{-\mathbf{X}}+(2\mu+\lambda)\pa_{_{\xi_1}}(v-(v^s)^{-\mathbf{X}}),
\ee
and using \eqref{stationary-2}$_1$ and \eqref{daoshu}, we have
$$
\ba
|\mathcal{Y}_2(t)|&\leq C\int_{\mathbb{T}^2}\int_{\mathbb{R}}|(v^s)^{-\mathbf{X}}_{_{\xi_1}}|\Big(\Big|h_1-(h^s_1)^{-\mathbf{X}}-\f{p(v)-p((v^s)^{-\mathbf{X}})}{\s_* }\Big|
+|p(v)-p((v^s)^{-\mathbf{X}})|\\
&\quad+\big|\pa_{_{\xi_1}}\big(p(v)-p((v^s)^{-\mathbf{X}})\big)\big|+|(v^s)^{-\mathbf{X}}_{_{\xi_1}}||v-(v^s)^{-\mathbf{X}}| \Big)d\xi_1d\xi'\\
&\leq C\Big(\f{\d}{\sqrt{\nu}}\sqrt{G_3(t)}+\sqrt{\d}\sqrt{G^s(t)}+\d^{\f32}\sqrt{D(t)}\Big).
\ea
$$
Thus, it holds
$$
\ba
|\dot{\mathbf{X}} (t)||\mathcal{Y}_2(t)|\leq \f{\d}{8}|\dot{\mathbf{X}}(t)|^2+C\f{\d}{\nu}G_3(t)+CG^s(t)+C\d^2D(t).
\ea
$$
For $I_1(t)$, using \eqref{use-1}, it holds
$$
|I_1(t)|\leq CG^s(t).
$$
For $I_2(t)$, using Lemma \ref{le-shock}, and the definition of $F$ \eqref{F}, it holds
$$
\ba
I_2(t) &\leq C\int_{\mathbb{T}^2}\int_{\mathbb{R}}(v^s)^{-\mathbf{X}}_{_{\xi_1}}\Big(\Big|h_1-(h^s_1)^{-\mathbf{X}}-\f{p(v)-p((v^s)^{-\mathbf{X}})}{\s_* }\Big|
+|p(v)-p((v^s)^{-\mathbf{X}})|\\
&\quad+\big|\pa_{_{\xi_1}}\big(p(v)-p((v^s)^{-\mathbf{X}})\big)\big|+(v^s)^{-\mathbf{X}}_{_{\xi_1}}|v-(v^s)^{-\mathbf{X}}| \Big)|p(v)-p((v^s)^{-\mathbf{X}})|d\xi_1d\xi'\\
&\leq C\Big(\sqrt{\f{\d}{\nu}}\sqrt{G_3(t)}\sqrt{G^s(t)}+G^s(t)+\d\sqrt{D(t)}\sqrt{G^s(t)}\Big)\\
&\leq C\f{\d}{\nu}G_3(t)+CG^s(t)+C\d^2D(t).
\ea
$$
Using \eqref{u1-u1^s} and \eqref{stationary-2}$_1$, similarly as $I_2(t)$, it holds
$$
I_3(t)\leq  C\f{\d}{\nu}G_3(t)+CG^s(t)+C\d^2D(t).
$$
Integrating \eqref{ba-expre-2} over $[0,t]$ for any $t\leq T$, and combining the above estimates, we can find that for some constant $C_2>0$,
$$
\ba
&\int_{\mathbb{T}^2}\int_{\mathbb{R}}
\r \Big(Q(v|(v^s)^{-\mathbf{X}})+\f{|u-(u^s)^{-\mathbf{X}}|^2}{2}\Big)d\xi_1d\xi'+\int_0^t\mathbf{D}_1(\tau)d\tau\\
&\leq C\int_{\mathbb{T}^2}\int_{\mathbb{R}}
 \rho_0\Big(Q(v_0|v^s)+\f{|u_0-u^s|^2}{2}\Big)d\xi_1d\xi'
+\int_0^t\Big(\f{\d}{4}|\dot{\mathbf{X}}(\tau)|^2+C\f{\d}{\nu}G_3(\tau)+C_2G^s(\tau)+C\d^2D(\tau)\Big)d\tau.
\ea
$$
Therefore, multiplying the above inequality by $\f{1}{2\max\{1,C_2\}}$, then adding the resultant equations to \eqref{basic} together,
choosing $\f{\d}{\nu}$, $\d$ and $\chi$ suitably small, we have
\be\label{es-u-L2}
\ba
&\|(v-(v^s)^{-\mathbf{X}})(\tau)\|^2+\|(h-(h^s)^{-\mathbf{X}})(\tau)\|^2+\|(u-(u^s)^{-\mathbf{X}})(\tau)\|^2\\
&+\d\int_0^t|\dot{\mathbf{X}}(\tau)|^2d\tau+\int_0^t(G_2(\tau)+G_3(\tau)+G^s(
\tau)+D(\tau)+\mathbf{D}_1(\tau))d\tau\\
&\leq C(\|v_0-v^s\|^2+\|h_0-h^s\|^2+\|u_0-u^s\|^2)
+C(\d+\chi)\int_0^t\|\na^2_{_\xi}(u-(u^s)^{-\mathbf{X}})\|^2d\tau,
\ea
\ee
where by Lemma \ref{inequality-Q}, we have used the fact
$$
C^{-1}|v-(v^s)^{-\mathbf{X}}|^2\leq Q(v|(v^s)^{-\mathbf{X}})\leq C|v-(v^s)^{-\mathbf{X}}|^2.
$$
Finally, from \eqref{h-h^s}, it holds
\be\label{v-1st-dao}
\|\na_{_\xi}(v-(v^s)^{-\mathbf{X}})\|\leq C(\|h-(h^s)^{-\mathbf{X}}\|+\|u-(u^s)^{-\mathbf{X}}\|)
\ee
and
\be\label{h0}
\|h_0-h^s\|\leq C(\|u_0-u^s\|+\|\na_{_\xi}(v_0-v^s\|).
\ee
Thus, combining \eqref{es-u-L2}, \eqref{v-1st-dao} and \eqref{h0}, and using the fact
$$
\|\na_{\xi}(u-(u^s)^{-\mathbf{X}})\|^2\sim\mathbf{D}_1,
$$
we can obtain  the desired inequality \eqref{es-1st}.
The proof of Lemma \ref{le-1st-v} is completed.

\hfill $\Box$

%
%
%
%

\subsection{Estimates for $\|\na_{_\xi}(u-(u^s)^{-\mathbf{X}})\|$}

\begin{lemma}\label{le-1st-u}
Under the hypotheses of Proposition \ref{priori}, there exists constant $C>0$ independent of $\nu$, $\d$, $\chi$ and $T$,
such that for all $t\in[0,T]$, it holds
\be\label{es-u-1st}
\|\na_{_\xi}(u-(u^s)^{-\mathbf{X}})(t)\|^2+\int_0^t\|\na_{_\xi}^2(u-(u^s)^{-\mathbf{X}})\|^2d\tau
\leq C\big(\|v_0-v^s\|_{H^1}^2+\|u_0-u^s\|_{H^1}^2\big).
\ee
\end{lemma}
\textbf{\emph{Proof}}: Multiplying \eqref{perturb-vu}$_2$ by $-v\Delta_{_\xi}(u-(u^s)^{-\mathbf{X}})$,
and integrating the resultant equations over $\mathbb{R}\times\mathbb{T}^2$, it holds
$$
\ba
&\f{d}{dt}\int_{\mathbb{T}^2}\int_{\mathbb{R}}\f{|\na_{_\xi}(u-(u^s)^{-\mathbf{X}})|^2}{2}d\xi_1d\xi'\\
&+\underbrace{\mu\int_{\mathbb{T}^2}\int_{\mathbb{R}}v|\Delta_{_\xi}(u-(u^s)^{-\mathbf{X}})|^2d\xi_1d\xi'
+(\mu+\lambda)\int_{\mathbb{T}^2}\int_{\mathbb{R}}v|\na_{_\xi}\div_{_\xi}(u-(u^s)^{-\mathbf{X}})|^2d\xi_1d\xi'}_{\mathbf{D}_2(t)}\\
=&-\int_{\mathbb{T}^2}\int_{\mathbb{R}}\Big(\na_{_\xi} u\cdot\na_{_\xi}(u-(u^s)^{-\mathbf{X}})\cdot\na_{_\xi}(u-(u^s)^{-\mathbf{X}})
-\div_{_\xi}u\f{|\na_{_\xi}(u-(u^s)^{-\mathbf{X}})|^2}{2}\Big)d\xi_1d\xi'\\
&-\dot{\mathbf{X}}(t)\int_{\mathbb{T}^2}\int_{\mathbb{R}}(u^s_1)^{-\mathbf{X}}_{_{\xi_1}}\Delta_{_\xi}(u_1-(u^s_1)^{-\mathbf{X}})d\xi_1d\xi'
+\int_{\mathbb{T}^2}\int_{\mathbb{R}}vF(u^s_1)^{-\mathbf{X}}_{_{\xi_1}}\Delta_{_\xi}(u_1-(u^s_1)^{-\mathbf{X}})d\xi_1d\xi'\\
&+(\mu+\lambda)\int_{\mathbb{T}^2}\int_{\mathbb{R}}\big(\Delta_{_\xi}(u-(u^s)^{-\mathbf{X}})-\na_{_\xi}\div_{_\xi}(u-(u^s)^{-\mathbf{X}})\big)
\cdot\na_{_\xi}v\div_{_\xi}(u-(u^s)^{-\mathbf{X}})d\xi_1d\xi'\\
&+\int_{\mathbb{T}^2}\int_{\mathbb{R}}v\Delta_{_\xi}(u-(u^s)^{-\mathbf{X}})\cdot\na_{_\xi}\big(p(v)-p((v^s)^{-\mathbf{X}})\big)d\xi_1d\xi'
=:\sum_{i=1}^5J_i(t).
\ea
$$
Using assumption \eqref{assumption} and Sobolev inequality, we have $\|\na_{_\xi}(u-(u^s)^{-\mathbf{X}})\|_{L^3}\leq C\|\na_{_\xi}(u-(u^s)^{-\mathbf{X}})\|_{H^1}\leq C\chi$.
Thus, it holds
$$
\ba
J_1(t)&\leq C\|\na_{_\xi}(u-(u^s)^{-\mathbf{X}})\|_{L^3}\|\na_{_\xi}(u-(u^s)^{-\mathbf{X}})\|_{L^6}\|\na_{_\xi}(u-(u^s)^{-\mathbf{X}})\|
+C\d^2\|\na_{_\xi}(u-(u^s)^{-\mathbf{X}})\|^2\\
&\leq C\chi\|\na_{_\xi}(u-(u^s)^{-\mathbf{X}})\|_{H^1}\|\na_{_\xi}(u-(u^s)^{-\mathbf{X}})\|+C\d^2\|\na_{_\xi}(u-(u^s)^{-\mathbf{X}})\|^2\\
&\leq C(\d+\chi)(\mathbf{D}_2(t)+\|\na_{_\xi}(u-(u^s)^{-\mathbf{X}})\|^2).
\ea
$$
We notice $-\s_* (v^s)^{-\mathbf{X}}_{_{\xi_1}}=(u^s_1)^{-\mathbf{X}}_{_{\xi_1}}$ from \eqref{stationary-2}$_2$, and use Lemma \ref{le-shock},
then we have
$$
\ba
J_2(t)&\leq |\dot{\mathbf{X}}(t)| \|(u^s_1)^{-\mathbf{X}}_{_{\xi_1}}\|_{L^2(\mathbb{R})}\sqrt{\mathbf{D}_2(t)}
\leq|\dot{\mathbf{X}}(t)| \d^{\f32}\sqrt{\mathbf{D}_2(t)}\leq C\d^2|\dot{\mathbf{X}}(t)|^2+C\d\mathbf{D}_2(t).
\ea
$$
For $J_3(t)$, using $-\s_* (v^s)^{-\mathbf{X}}_{_{\xi_1}}=(u^s_1)^{-\mathbf{X}}_{_{\xi_1}}$ and Lemma \ref{le-shock}, and the definition of $F$ \eqref{F}, it holds
$$
\ba
J_3(t) &\leq C\int_{\mathbb{T}^2}\int_{\mathbb{R}}(v^s)^{-\mathbf{X}}_{_{\xi_1}}\Big(\Big|h_1-(h^s_1)^{-\mathbf{X}}-\f{p(v)-p((v^s)^{-\mathbf{X}})}{\s_* }\Big|
+|p(v)-p((v^s)^{-\mathbf{X}})|\\
&\quad+\big|\pa_{_{\xi_1}}\big(p(v)-p((v^s)^{-\mathbf{X}})\big)\big|+(v^s)^{-\mathbf{X}}_{_{\xi_1}}|v-(v^s)^{-\mathbf{X}}| \Big)|\Delta_{_\xi}(u_1-(u^s_1)^{-\mathbf{X}})|d\xi_1d\xi'\\
&\leq C\Big(\d\sqrt{\f{\d}{\nu}}\sqrt{G_3(t)}+\d\sqrt{G^s(t)}+\d^2\sqrt{D(t)}\Big)\sqrt{\mathbf{D}_2(t)}\\
&
\leq C\f{\d}{\nu}G_3(t)+CG^s(t)+C\d^2D(t)+C\d^2\mathbf{D}_2(t).
\ea
$$
For $J_4(t)$, using assumption we have $\|\na_{_\xi}(v-(v^s)^{-\mathbf{X}})\|_{L^3}\leq C\|\na_{_\xi}(v-(v^s)^{-\mathbf{X}})\|_{H^1}\leq C\chi$, it holds
$$
\ba
J_4(t)&\leq C(\|\Delta_{_\xi}(u-(u^s)^{-\mathbf{X}})\|+\|\na_{_\xi}\div_{_\xi}(u-(u^s)^{-\mathbf{X}})\|)\|(v^s)^{-\mathbf{X}}_{_{\xi_1}}\|_{L^{\infty}}
\|\div_{_\xi}(u-(u^s)^{-\mathbf{X}})\|\\
&+C(\|\Delta_{_\xi}(u-(u^s)^{-\mathbf{X}})\|+\|\na_{_\xi}\div_{_\xi}(u-(u^s)^{-\mathbf{X}})\|)\|\na_{_\xi}(v-(v^s)^{-\mathbf{X}})\|_{L^3}
\|\div_{_\xi}(u-(u^s)^{-\mathbf{X}})\|_{L^6}\\
&\leq C\d^2\sqrt{\mathbf{D}_2(t)}\|\div_{_\xi}(u-(u^s)^{-\mathbf{X}})\|+C\chi\sqrt{\mathbf{D}_2(t)}\|\div_{_\xi}(u-(u^s)^{-\mathbf{X}})\|_{H^1}\\
&\leq C(\d+\chi)(\mathbf{D}_2(t)+\|\na_{_\xi}(u-(u^s)^{-\mathbf{X}})\|^2).
\ea
$$
By Cauchy inequality, we have
$$
J_5(t)\leq \f14\mathbf{D}_2(t)+C\|\na_{_\xi}\big(p(v)-p((v^s)^{-\mathbf{X}})\big)\|^2=\f14\mathbf{D}_2(t)+CD(t).
$$
Therefore, the combination of the above estimates yields
$$
\ba&
\f{d}{dt}\|\na_{_\xi}(u-(u^s)^{-\mathbf{X}})(t)\|^2+\mathbf{D}_2(t)\\
&\leq C\d^2|\dot{\mathbf{X}}(t)|^2+C\f{\d}{\nu}G_3(t)+CG^s(t)+CD(t)
+C(\d+\chi)\|\na_{_\xi}(u-(u^s)^{-\mathbf{X}})\|^2.
\ea
$$
Integrating the above inequality over $[0,t]$ for any $t\leq T$, and using Lemma \ref{le-basic-new},
Lemma \ref{le-1st-v} and the fact that
$$
\|\na_{_\xi}^2(u-(u^s)^{-\mathbf{X}})(t)\|^2 \sim \mathbf{D}_2(t),
$$
we can obtain the desired inequality \eqref{es-u-1st}. The proof of Lemma \ref{le-1st-u} is completed.

\hfill $\Box$

%
%
%
%

\subsection{Estimates for $\|\na^2_{_\xi}(v-(v^s)^{-\mathbf{X}})\|$}

\begin{lemma}\label{le-2nd-v}
Under the hypotheses of Proposition \ref{priori}, there exists constant $C>0$ independent of $\nu$, $\d$, $\chi$ and $T$,
such that for all $t\in[0,T]$, it holds
\be\label{es-v-2nd}
\ba
&\|\na_{_\xi}^2(v-(v^s)^{-\mathbf{X}})(t)\|^2+\int_0^t\|\na_{_\xi}^2(v-(v^s)^{-\mathbf{X}})\|^2d\tau\\
&\leq C\big(\|v_0-v^s\|_{H^2}^2+\|u_0-u^s\|_{H^1}^2\big)+C(\d+\chi)\int_0^t\|\na_{_\xi}^3(u-(u^s)^{-\mathbf{X}})\|^2d\tau.
\ea
\ee
\end{lemma}
\textbf{\emph{Proof}}: We set $\phi:=v-(v^s)^{-\mathbf{X}}$, $\psi:=u-(u^s)^{-\mathbf{X}}$ for notational simplicity, and rewrite \eqref{perturb-vu} as
\be\label{re-perturb-vu}
\left\{
\ba
\pa_t \phi&-\s\pa_{_{\xi_1}}\phi+ u\cdot\na_{_\xi}\phi
-\dot{\mathbf{X}}(t)(v^s)^{-\mathbf{X}}_{_{\xi_1}}+vF(v^s)^{-\mathbf{X}}_{_{\xi_1}}=v\div_{_\xi}\psi,\\
\pa_t \psi&-\s\pa_{_{\xi_1}}\psi+ u\cdot\na_{_\xi}\psi+vp'(v)\na_{_\xi}\phi
+v\big(p'(v)-p'((v^s)^{-\mathbf{X}})\big)\na_{_\xi}(v^s)^{-\mathbf{X}}\\
&-\dot{\mathbf{X}}(t)(u^s)^{-\mathbf{X}}_{_{\xi_1}}+vF(u^s)^{-\mathbf{X}}_{_{\xi_1}}
=\mu v\Delta_{_\xi}\psi+(\mu+\lambda)v\na_{_\xi}\div_{_\xi}\psi.
\ea
\right.
\ee
Applying $\na_{_\xi}\pa_{_{\xi_i}}~(i=1,2,3)$ to \eqref{re-perturb-vu}$_1$, and $\pa_{_{\xi_i}}~(i=1,2,3)$ to \eqref{re-perturb-vu}$_2$,
we have
\be\label{dao-perturb-vu}
\left\{
\ba
&\pa_t\na_{_\xi}\pa_{_{\xi_i}} \phi-\s\pa_{_{\xi_1}}\na_{_\xi}\pa_{_{\xi_i}}\phi+ u\cdot\na_{_\xi}(\na_{_\xi}\pa_{_{\xi_i}}\phi)
-\dot{\mathbf{X}}(t)\na_{_\xi}\pa_{_{\xi_i}}(v^s)^{-\mathbf{X}}_{_{\xi_1}}+vF\na_{_\xi}\pa_{_{\xi_i}}(v^s)^{-\mathbf{X}}_{_{\xi_1}}\\
&\quad+\na_{_\xi}\pa_{_{\xi_i}}u\cdot\na_{_\xi}\phi+\nabla_{_\xi} u\cdot\nabla_{_\xi}\partial_{_{\xi_i}}\phi+\partial_{_{\xi_i}} u\cdot\na_{_\xi}(\nabla_{_\xi}\phi)
+\na_{_\xi}\pa_{_{\xi_i}}(v F)(v^s)^{-\mathbf{X}}_{_{\xi_1}}
+\na_{_{\xi}}(vF)\pa_{_{\xi_i}}(v^s)^{-\mathbf{X}}_{_{\xi_1}}\\
&\quad+\pa_{_{\xi_i}}(vF)\na_{_{\xi}}(v^s)^{-\mathbf{X}}_{_{\xi_1}}
=v\na_{_\xi}\pa_{_{\xi_i}}\div_{_\xi}\psi+\na_{_\xi}\pa_{_{\xi_i}}v\div_{_\xi}\psi
+\partial_{_{\xi_i}}v\nabla_{_\xi}\div_{_\xi} \psi+\nabla_{_\xi} v \partial_{_{\xi_i}}\div_{_\xi} \psi,\\
&\pa_t\pa_{_{\xi_i}} \psi-\s\pa_{_{\xi_1}}\pa_{_{\xi_i}}\psi+ u\cdot\na_{_\xi}\pa_{_{\xi_i}}\psi+vp'(v)\na_{_\xi}\pa_{_{\xi_i}}\phi
+\pa_{_{\xi_i}}u\cdot\na_{_\xi}\psi+\pa_{_{\xi_i}}(vp'(v))\na_{_\xi}\phi\\
&\quad+\pa_{_{\xi_i}}\big(v(p'(v)-p'((v^s)^{-\mathbf{X}}))\big)\na_{_\xi}(v^s)^{-\mathbf{X}}
+v\big(p'(v)-p'((v^s)^{-\mathbf{X}})\big)\na_{_\xi}\pa_{_{\xi_i}}(v^s)^{-\mathbf{X}}\\
&\quad-\dot{\mathbf{X}}(t)\pa_{_{\xi_i}}(u^s)^{-\mathbf{X}}_{_{\xi_1}}
+vF\pa_{_{\xi_i}}(u^s)^{-\mathbf{X}}_{_{\xi_1}}+\pa_{_{\xi_i}}(v F)(u^s)^{-\mathbf{X}}_{_{\xi_1}}\\
&=\mu v\Delta_{_\xi}\pa_{_{\xi_i}}\psi+(\mu+\lambda)v\na_{_\xi}\pa_{_{\xi_i}}\div_{_\xi}\psi
+\pa_{_{\xi_i}}v\big(\mu \Delta_{_\xi}\psi+(\mu+\lambda)\na_{_\xi}\div_{_\xi}\psi\big).
\ea
\right.
\ee
Multiplying \eqref{dao-perturb-vu}$_1$ by $\r(2\mu+\lambda)\na_{_\xi}{\pa_{_{\xi_i}}}\phi$, and summating $i$ from 1 to 3,
then integrating the resultant equations over $\mathbb{R}\times\mathbb{T}^2$, it holds
\begin{align}\label{v-2nd-1} \nm
&(2\mu+\lambda)\f{d}{dt}\int_{\mathbb{T}^2}\int_{\mathbb{R}}\r\f{|\na_{_\xi}^2\phi|^2}{2}d\xi_1d\xi'
-(2\mu+\lambda)\sum_{i=1}^3\int_{\mathbb{T}^2}\int_{\mathbb{R}}\na_{_\xi}\pa_{_{\xi_i}}\phi\cdot\na_{_\xi}\pa_{_{\xi_i}}\div_{_\xi}\psi d\xi_1d\xi'\\ \nm
&=(2\mu+\lambda)\dot{\mathbf{X}}(t)\int_{\mathbb{T}^2}\int_{\mathbb{R}}\r\pa_{_{\xi_1}}^2\phi(v^s)^{-\mathbf{X}}_{_{\xi_1\xi_1\xi_1}} d\xi_1d\xi'
-(2\mu+\lambda)\int_{\mathbb{T}^2}\int_{\mathbb{R}}F\pa_{_{\xi_1}}^2\phi(v^s)^{-\mathbf{X}}_{_{\xi_1\xi_1\xi_1}} d\xi_1d\xi'\\
&-(2\mu+\lambda)\sum_{i=1}^3\int_{\mathbb{T}^2}\int_{\mathbb{R}}\r\na_{_\xi}\pa_{_{\xi_i}}\phi\cdot\big[\na_{_\xi}\pa_{_{\xi_i}}u\cdot\na_{_\xi}\phi+\nabla_{_\xi} u\cdot\nabla_{_\xi}\partial_{_{\xi_i}}\phi+\partial_{_{\xi_i}} u\cdot\na_{_\xi}(\nabla_{_\xi}\phi)\big] d\xi_1d\xi'\\ \nm
&-(2\mu+\lambda)\sum_{i=1}^3\int_{\mathbb{T}^2}\int_{\mathbb{R}}\r\big(\na_{_\xi}\pa_{_{\xi_i}}\phi\cdot\na_{_\xi}\pa_{_{\xi_i}}(vF)(v^s)^{-\mathbf{X}}_{_{\xi_1}} +\pa_{_{\xi_1}}\pa_{_{\xi_i}}\phi\pa_{_{\xi_i}}(vF)(v^s)^{-\mathbf{X}}_{_{\xi_1\xi_1}}\big) d\xi_1d\xi'\\ \nm
&-(2\mu+\lambda)\int_{\mathbb{T}^2}\int_{\mathbb{R}}\r\na_{_\xi}\pa_{_{\xi_1}}\phi\cdot\na_{_\xi}(vF)
(v^s)^{-\mathbf{X}}_{_{\xi_1\xi_1}}d\xi_1d\xi'\\\nm
&+(2\mu+\lambda)\sum_{i=1}^3\int_{\mathbb{T}^2}\int_{\mathbb{R}}\r\na_{_\xi}\pa_{_{\xi_i}}\phi\cdot\big[\na_{_\xi}\pa_{_{\xi_i}}v\div_{_\xi}\psi+\partial_{_{\xi_i}}v\nabla_{_\xi}\div_{_\xi} \psi+\nabla_{_\xi} v \partial_{_{\xi_i}}\div_{_\xi} \psi\big] d\xi_1d\xi'.\nm
\end{align}
Multiplying \eqref{dao-perturb-vu}$_2$ by $-\r\na_{_\xi}{\pa_{_{\xi_i}}}\phi$, and summating $i$ from 1 to 3,
then integrating the resultant equations over $\mathbb{R}\times\mathbb{T}^2$, it holds
\begin{align}\nm  \label{v-2nd-2}
&\int_{\mathbb{T}^2}\int_{\mathbb{R}}-p'(v)|\na_{_\xi}^2\phi|^2d\xi_1d\xi'
+(2\mu+\lambda)\sum_{i=1}^3\int_{\mathbb{T}^2}\int_{\mathbb{R}}\na_{_\xi}\pa_{_{\xi_i}}\phi\cdot\na_{_\xi}\pa_{_{\xi_i}}\div_{_\xi}\psi d\xi_1d\xi'\\ \nm
&=\f{d}{dt}\sum_{i=1}^3\int_{\mathbb{T}^2}\int_{\mathbb{R}}\r\pa_{_{\xi_i}}\psi\cdot\na_{_\xi}\pa_{_{\xi_i}}\phi d\xi_1d\xi'+\sum_{i=1}^3\int_{\mathbb{T}^2}\int_{\mathbb{R}}\r \na_{_\xi}\pa_{_{\xi_i}}\phi\cdot\pa_{_{\xi_i}}u\cdot\na_{_\xi}\psi d\xi_1d\xi'\\\nm
&-\sum_{i=1}^3\int_{\mathbb{T}^2}\int_{\mathbb{R}}\r\pa_{_{\xi_i}}\psi\cdot\big[\na_{_\xi}\pa_{_{\xi_i}}\pa_{_t}\phi-\s
\pa_{_{\xi_1}} \na_{_\xi}\pa_{_{\xi_i}}\phi+u\cdot\na_{_\xi}(\na_{_\xi}\pa_{_{\xi_i}}\phi)\big]d\xi_1d\xi'\\ \nm
&
+\sum_{i=1}^3\int_{\mathbb{T}^2}\int_{\mathbb{R}}\r \pa_{_{\xi_i}}(vp'(v))\na_{_\xi}\pa_{_{\xi_i}}\phi\cdot\na_{_\xi}\phi d\xi_1d\xi'-\dot{\mathbf{X}}(t)\int_{\mathbb{T}^2}\int_{\mathbb{R}}\r\pa_{_{\xi_1}}^2\phi(u^s_1)^{-\mathbf{X}}_{_{\xi_1\xi_1}} d\xi_1d\xi'\\
&+\sum_{i=1}^3\int_{\mathbb{T}^2}\int_{\mathbb{R}}\r \pa_{_{\xi_i}}\big(v(p'(v)-p'((v^s)^{-\mathbf{X}}))\big)\pa_{_{\xi_1}}\pa_{_{\xi_i}}\phi(v^s)^{-\mathbf{X}}_{_{\xi_1}}d\xi_1d\xi'\\  \nm
&+\int_{\mathbb{T}^2}\int_{\mathbb{R}}(p'(v)-p'((v^s)^{-\mathbf{X}}))\pa_{_{\xi_1}}^2\phi(v^s)^{-\mathbf{X}}_{_{\xi_1\xi_1}}d\xi_1d\xi'\\  \nm
&+\int_{\mathbb{T}^2}\int_{\mathbb{R}}F\pa_{_{\xi_1}}^2\phi(u^s_1)^{-\mathbf{X}}_{_{\xi_1\xi_1}} d\xi_1d\xi'
+\sum_{i=1}^3\int_{\mathbb{T}^2}\int_{\mathbb{R}}\r\pa_{_{\xi_i}}(vF)\pa_{_{\xi_1}}\pa_{_{\xi_i}}\phi(u^s_1)^{-\mathbf{X}}_{_{\xi_1}} d\xi_1d\xi'\\  \nm
&-\sum_{i=1}^3\int_{\mathbb{T}^2}\int_{\mathbb{R}}\r\pa_{_{\xi_i}}v(\mu\Delta_{_\xi}\psi+(\mu+\lambda)\na_{_\xi}\div_{_\xi}\psi)\cdot\na_{_\xi}\pa_{_{\xi_i}}\phi d\xi_1d\xi'.
\end{align}
Adding \eqref{v-2nd-1} and \eqref{v-2nd-2} together, and integrating the resultant equations over $[0,t]$ for any $t\in[0, T]$, we have
\be\label{v-2nd-3}
(2\mu+\lambda)\int_{\mathbb{T}^2}\int_{\mathbb{R}}\r\f{|\na_{_\xi}^2\phi|^2}{2}d\xi_1d\xi'\Big|_{\tau=0}^{\tau=t}
+\int_0^t\int_{\mathbb{T}^2}\int_{\mathbb{R}}|p'(v)||\na_{_\xi}^2\phi|^2d\xi_1d\xi'd\tau\\
=\sum_{j=1}^{8}K_j(t),
\ee
where
\begin{align}  \nm
K_1(t)&=\sum_{i=1}^3\int_{\mathbb{T}^2}\int_{\mathbb{R}}\r\pa_{_{\xi_i}}\psi\cdot\na_{_\xi}\pa_{_{\xi_i}}\phi d\xi_1d\xi'\Big|_{\tau=0}^{\tau=t},\quad\\ \nm
K_2(t)&=\int_0^t\dot{\mathbf{X}}(\tau)\int_{\mathbb{T}^2}\int_{\mathbb{R}}\r\pa_{_{\xi_1}}^2\phi\big[(2\mu+\lambda)(v^s)^{-\mathbf{X}}_{_{\xi_1\xi_1\xi_1}}-(u^s_1)^{-\mathbf{X}}_{_{\xi_1\xi_1}}\big] d\xi_1d\xi'd\tau,\\ \nm
K_3(t)&=-\int_0^t\int_{\mathbb{T}^2}\int_{\mathbb{R}}F\pa_{_{\xi_1}}^2\phi\big[(2\mu+\lambda)(v^s)^{-\mathbf{X}}_{_{\xi_1\xi_1\xi_1}} -(u^s_1)^{-\mathbf{X}}_{_{\xi_1\xi_1}}\big] d\xi_1d\xi'd\tau,\\  \nm
K_4(t)&=
-(2\mu+\lambda)\sum_{i=1}^3\int_0^t\int_{\mathbb{T}^2}\int_{\mathbb{R}}\r\big(\na_{_\xi}\pa_{_{\xi_i}}\phi\cdot\na_{_\xi}\pa_{_{\xi_i}}(vF)(v^s)^{-\mathbf{X}}_{_{\xi_1}} +\pa_{_{\xi_1}}\pa_{_{\xi_i}}\phi\pa_{_{\xi_i}}(vF)(v^s)^{-\mathbf{X}}_{_{\xi_1\xi_1}}\big) d\xi_1d\xi'd\tau\\ \nm
&-(2\mu+\lambda)\int_0^t\int_{\mathbb{T}^2}\int_{\mathbb{R}}\r\na_{_\xi}\pa_{_{\xi_1}}\phi\cdot\na_{_\xi}(vF)
(v^s)^{-\mathbf{X}}_{_{\xi_1\xi_1}}d\xi_1d\xi'd\tau\\ \nm
&+\sum_{i=1}^3\int_0^t\int_{\mathbb{T}^2}\int_{\mathbb{R}}\r\pa_{_{\xi_1}}\pa_{_{\xi_i}}\phi\pa_{_{\xi_i}}(vF)(u^s_1)^{-\mathbf{X}}_{_{\xi_1}} d\xi_1d\xi'd\tau,\\ \nm
K_5(t)&=-\sum_{i=1}^3\int_0^t\int_{\mathbb{T}^2}\int_{\mathbb{R}}\r\pa_{_{\xi_i}}\psi\cdot\big[\na_{_\xi}\pa_{_{\xi_i}}\pa_{_t}\phi-\s \pa_{_{\xi_1}}\na_{_\xi}\pa_{_{\xi_i}}\phi+u\cdot\na_{_\xi}(\na_{_\xi}\pa_{_{\xi_i}}\phi)\big]d\xi_1d\xi'd\tau,\\ \nm
K_6(t)&=-(2\mu+\lambda)\sum_{i=1}^3\int_0^t\int_{\mathbb{T}^2}\int_{\mathbb{R}}\r\na_{_\xi}\pa_{_{\xi_i}}\phi\cdot\big[\na_{_\xi}\pa_{_{\xi_i}}u\cdot\na_{_\xi}\phi+\nabla_{_\xi} u\cdot\nabla_{_\xi}\partial_{_{\xi_i}}\phi+\partial_{_{\xi_i}} u\cdot\na_{_\xi}(\nabla_{_\xi}\phi)\big]d\xi_1d\xi'd\tau\\ \nm
&+\sum_{i=1}^3\int_0^t\int_{\mathbb{T}^2}\int_{\mathbb{R}}\r \na_{_\xi}\pa_{_{\xi_i}}\phi\cdot\big[\pa_{_{\xi_i}}u\cdot\na_{_\xi}\psi+ \pa_{_{\xi_i}}(vp'(v))\na_{_\xi}\phi \big]d\xi_1d\xi'd\tau,\\  \nm
K_7(t)&=\sum_{i=1}^3\int_0^t\int_{\mathbb{T}^2}\int_{\mathbb{R}}\r \pa_{_{\xi_i}}\big(v(p'(v)-p'((v^s)^{-\mathbf{X}}))\big)\pa_{_{\xi_1}}\pa_{_{\xi_i}}\phi(v^s)^{-\mathbf{X}}_{_{\xi_1}}d\xi_1d\xi'd\tau\\ \nm
&+\int_0^t\int_{\mathbb{T}^2}\int_{\mathbb{R}}(p'(v)-p'((v^s)^{-\mathbf{X}}))\pa_{_{\xi_1}}^2\phi(v^s)^{-\mathbf{X}}_{_{\xi_1\xi_1}}d\xi_1d\xi'd\tau,\\ \nm
K_8(t)&=(2\mu+\lambda)\sum_{i=1}^3\int_0^t\int_{\mathbb{T}^2}\int_{\mathbb{R}}\r\na_{_\xi}\pa_{_{\xi_i}}\phi\cdot\big[\na_{_\xi}\pa_{_{\xi_i}}v\div_{_\xi}\psi+\partial_{_{\xi_i}}v\nabla_{_\xi}(\div_{_\xi} \psi)+\nabla_{_\xi} v \partial_{_{\xi_i}}(\div_{_\xi} \psi)\big] d\xi_1d\xi'd\tau\\\nm
&-\sum_{i=1}^3\int_0^t\int_{\mathbb{T}^2}\int_{\mathbb{R}}\r\pa_{_{\xi_i}}v \na_{_\xi}\pa_{_{\xi_i}}\phi\cdot\big[\mu\Delta_{_\xi}\psi+(\mu+\lambda)\na_{_\xi}\div_{_\xi}\psi\big] d\xi_1d\xi'd\tau.
\end{align}
Using Cauchy inequality, it holds
$$
K_1(t)\leq \f{2\mu+\lambda}{8}\|\sqrt{\r}\na_{_\xi}^2\phi(t)\|^2+C\|\na_{_\xi}\psi(t)\|^2+C(\|\na^2_{_\xi}\phi_0\|^2+\|\na_{_\xi}\psi_0\|^2).
$$
Using Lemma \ref{le-shock}, we have
$$
K_2(t)\leq C\d\int_0^t|\dot{\mathbf{X}}(\tau)|\,\|\pa_{_{\xi_1}}^2\phi\|\,\|(v^s)^{-\mathbf{X}}_{_{\xi_1}}\|_{L^2(\mathbb{R})}d\tau
\leq C\d^2\int_0^t(|\dot{\mathbf{X}}(\tau)|^2+\|\sqrt{|p'(v)|}\pa_{_{\xi_1}}^2\phi\|^2)d\tau.
$$
Using \eqref{F}, we have
\begin{align} \nm
K_3(t)&\leq C\d\int_0^t\int_{\mathbb{T}^2}\int_{\mathbb{R}}\Big(\Big|h_1-(h^s_1)^{-\mathbf{X}}-\f{p(v)-p((v^s)^{-\mathbf{X}})}{\s_*}\Big|
+|p(v)-p((v^s)^{-\mathbf{X}})|\\ \nm
&\quad+\big|\pa_{_{\xi_1}}\big(p(v)-p((v^s)^{-\mathbf{X}})\big)\big|+|(v^s)^{-\mathbf{X}}_{_{\xi_1}}||p(v)-p((v^s)^{-\mathbf{X}})|\Big)|\pa_{_{\xi_1}}^2\phi|
|(v^s)^{-\mathbf{X}}_{_\xi}|d\xi_1d\xi'd\tau\\ \nm
&\leq C\d^2\int_0^t\big(\sqrt{G_3(\tau)}+\sqrt{G^s(\tau)}+\sqrt{D(\tau)}\big)\|\pa_{_{\xi_1}}^2\phi\|d\tau\\ \nm
&
\leq C\d^2\int_0^t(G_3(\tau)+G^s(\tau)+D(\tau)+\|\sqrt{|p'(v)|}\pa_{_{\xi_1}}^2\phi\|^2)d\tau.
\end{align}
By using \eqref{F} again, it holds
\be\label{vF}
v F=\s_* (v-(v^s)^{-\mathbf{X}})+u_1-(u^s_1)^{-\mathbf{X}}=\s_* \phi+\psi_1.
\ee
Thus we can get
$$
\ba
K_4(t)&\leq C\int_0^t\int_{\mathbb{T}^2}\int_{\mathbb{R}}|\na_{_\xi}^2\phi|(|\na_{_\xi}^2\phi|+|\na_{_\xi}^2\psi|)|(v^s)^{-\mathbf{X}}_{_{\xi_1}}|d\xi_1d\xi'd\tau\\
&\qquad+C\int_0^t\int_{\mathbb{T}^2}\int_{\mathbb{R}}|\na_{_\xi}^2\phi|(|\na_{_\xi}\phi|+|\na_{_\xi}\psi|)|(v^s)^{-\mathbf{X}}_{_{\xi_1}}|d\xi_1d\xi'd\tau\\
&\leq C\d^2\int_0^t\|\na_{_\xi}^2\phi\|(\|\na_{_\xi}\phi\|_{H^1}+\|\na_{_\xi}\psi\|_{H^1})d\tau \\
&\leq C\d^2\int_0^t(\|\sqrt{|p'(v)|}\na_{_\xi}^2\phi\|^2+D(\tau)+G^s(\tau)+\|\na_{_\xi}\psi\|^2_{H^1})d\tau.
\ea
$$
By using the equation \eqref{dao-perturb-vu}$_1$, we can compute the term $K_5(t)$ as
$$
\ba
K_5(t)&=-\int_0^t\dot{\mathbf{X}}(\tau)\int_{\mathbb{T}^2}\int_{\mathbb{R}}\r \pa_{\xi_1}\psi_1(v^s)^{-\mathbf{X}}_{_{\xi_1\xi_1\xi_1}}d\xi_1d\xi'd\tau
+\int_0^t\int_{\mathbb{T}^2}\int_{\mathbb{R}}F\pa_{\xi_1}\psi_1(v^s)^{-\mathbf{X}}_{_{\xi_1\xi_1\xi_1}}d\xi_1d\xi'd\tau\\
&\quad +\sum_{i=1}^3\int_0^t\int_{\mathbb{T}^2}\int_{\mathbb{R}}\r\pa_{\xi_i}\psi\cdot\big[\nabla_{\xi}\pa_{\xi_i}u\cdot\nabla_\xi\phi+\nabla_{\xi} u\cdot\nabla_\xi\pa_{\xi_i}\phi+\pa_{\xi_i}u\cdot\na_{\xi}(\na_{\xi}\phi)\big]d\xi_1d\xi'd\tau\\
&\quad +\sum_{i=1}^3\int_0^t\int_{\mathbb{T}^2}\int_{\mathbb{R}}\r\pa_{\xi_i}\psi\cdot\big[\nabla_{\xi}\pa_{\xi_i}(v F) (v^s)^{-\mathbf{X}}_{\xi_1}+\na_{_{\xi}}(vF)\pa_{_{\xi_i}}(v^s)^{-\mathbf{X}}_{_{\xi_1}}
+\pa_{_{\xi_i}}(vF)\na_{_{\xi}}(v^s)^{-\mathbf{X}}_{_{\xi_1}}\big]d\xi_1d\xi'd\tau\\
&\quad-\sum_{i=1}^3\int_0^t\int_{\mathbb{T}^2}\int_{\mathbb{R}}\r\pa_{\xi_i}\psi\cdot\big[\na_{\xi}\pa_{\xi_i}v\div_{\xi}\psi+\partial_{\xi_i}v\nabla_\xi\div_\xi \psi+\nabla_\xi v \partial_{\xi_i}\div_\xi \psi\big]d\xi_1d\xi'd\tau\\
&\quad-\sum_{i=1}^3\int_0^t\int_{\mathbb{T}^2}\int_{\mathbb{R}}\pa_{\xi_i}\psi\cdot\na_{\xi}\pa_{\xi_i}\div_{\xi}\psi d\xi_1d\xi'd\tau=:\sum_{i=1}^6K_{5,i}(t).
\ea
$$
Using Lemma \ref{le-shock}, we have
$$
K_{5,1}(t)\leq C\d^2\int_0^t|\dot{\mathbf{X}}(\tau)|^2d\tau+C\d^2\int_0^t\|\na_{_{\xi}}\psi\|^2d\tau.
$$
Using \eqref{F} and \eqref{daoshu}, we have
$$
\ba
K_{5,2}(t)&\leq C \int_0^t\int_{\mathbb{T}^2}\int_{\mathbb{R}}|\na_{_{\xi}}\psi||(v^s)^{-\mathbf{X}}_{_{\xi_1\xi_1\xi_1}}|\Big(
\Big|h_1-(h^s_1)^{-\mathbf{X}}-\f{p(v)-p((v^s)^{-\mathbf{X}})}{\s_*}\Big|\\
&\quad+|p(v)-p((v^s)^{-\mathbf{X}})|+|\pa_{_{\xi_1}}\big(p(v)-p((v^s)^{-\mathbf{X}})\big)|(v^s)^{-\mathbf{X}}_{_{\xi_1}}\Big)d\xi_1d\xi'd\tau\\
&\leq C\d^2\int_0^t(\|\na_{_{\xi}}\psi\|^2+G_3(\tau)+G^s(\tau)+D(\tau))d\tau.
\ea
$$
By assumption \eqref{assumption}, Cauchy inequality and Sobolev inequality, it holds
\be\label{L3}
\ba
&\|\phi\|_{L^3}+\|\psi\|_{L^3}\leq C\|\phi\|^{\f12}\|\phi\|_{L^6}^{\f12}+C\|\psi\|^{\f12}\|\psi\|_{L^6}^{\f12}
\leq C\|\phi\|_{H^1}+ C\|\psi\|_{H^1}\leq C\chi,\\
&\|\na_{_{\xi}}\phi\|_{L^3}+\|\na_{_{\xi}}\psi\|_{L^3}\leq C\|\na_{_{\xi}}\phi\|^{\f12}\|\na_{_{\xi}}\phi\|_{L^6}^{\f12}
+C\|\na_{_{\xi}}\psi\|^{\f12}\|\na_{_{\xi}}\psi\|_{L^6}^{\f12}\\
&\qquad\qquad\qquad\qquad\quad\leq C\|\na_{_{\xi}}\phi\|_{H^1}+ C\|\na_{_{\xi}}\psi\|_{H^1}\leq C\chi.
\ea
\ee
Then we have
$$
\ba
K_{5,3}(t)&\leq C \int_0^t\int_{\mathbb{T}^2}\int_{\mathbb{R}}|\na_{_{\xi}}\psi|\big(|\na^2_{_{\xi}}\psi||\na_{_{\xi}}\phi|+|\na_{_{\xi}}\psi||\na^2_{_{\xi}}\phi|\big)
d\xi_1d\xi'd\tau\\
&\quad+C\int_0^t\int_{\mathbb{T}^2}\int_{\mathbb{R}}|\na_{_{\xi}}\psi|\big(|(u^s_1)^{-\mathbf{X}}_{_{\xi_1\xi_1}}||\pa_{_{\xi_1}}\phi|
+|(u^s_1)^{-\mathbf{X}}_{\xi_1}||\na_{_{\xi}}\pa_{_{\xi_1}}\phi|\big)d\xi_1d\xi'd\tau\\
&\leq C\int_0^t\|\na_{_{\xi}}\psi\|_{L^6}\|\na_{_{\xi}}(\phi,\psi)\|_{L^3}\|\na^2_{_{\xi}}(\phi,\psi)\|d\tau
+C\d^2\int_0^t\|\na_{_{\xi}}\psi\|\big(\|\pa_{_{\xi_1}}\phi\|+\|\na_{_{\xi}}\pa_{_{\xi_1}}\phi\|\big)d\tau\\
&\leq  C\chi\int_0^t(\|\na_{_{\xi}}^2\phi\|^2+\|\na_{_{\xi}}\psi\|_{H^1}^2)d\tau
+C\d^2\int_0^t(\|\na_{_{\xi}}^2\phi\|^2+\|\na_{_{\xi}}\psi\|^2+D(\tau)+G^s(\tau))d\tau.
\ea
$$
We use \eqref{vF} again and Lemma \ref{le-shock}, it holds
$$
\ba
K_{5,4}(t)&\leq C\d^2\int_0^t(\|\na_{_{\xi}}^2(\phi,\psi)\|^2+\|\na_{_{\xi}}(\phi,\psi)\|^2)d\tau\\
&\leq C\d^2\int_0^t(\|\na_{_{\xi}}^2\phi\|^2+\|\na_{_{\xi}}\psi\|_{H^1}^2+D(\tau)+G^s(\tau))d\tau.
\ea
$$
Similarly to $K_{5,3}(t)$, we have
$$
\ba
K_{5,5}(t)&\leq C \int_0^t\int_{\mathbb{T}^2}\int_{\mathbb{R}}|\na_{_{\xi}}\psi|\big(|\na_{_{\xi}}\psi||\na^2_{_{\xi}}\phi|+|\na_{_{\xi}}\phi||\na^2_{_{\xi}}\psi|\big)
d\xi_1d\xi'd\tau\\
&\quad+C\int_0^t\int_{\mathbb{T}^2}\int_{\mathbb{R}}|\na_{_{\xi}}\psi|\big(|(v^s)^{-\mathbf{X}}_{_{\xi_1\xi_1}}||\na_{_{\xi}}\psi|
+|(v^s)^{-\mathbf{X}}_{_{\xi_1}}||\na^2_{\xi}\psi|\big)d\xi_1d\xi'd\tau\\
&\leq C\int_0^t\|\na_{_{\xi}}\psi\|_{L^6}\|\na_{_{\xi}}(\phi,\psi)\|_{L^3}\|\na^2_{_{\xi}}(\phi,\psi)\|d\tau
+C\d^2\int_0^t\|\na_{_{\xi}}\psi\|\big(\|\na_{_{\xi}}\psi\|+\|\na_{_{\xi}}^2\psi\|\big)d\tau\\
&\leq C\chi\int_0^t(\|\na_{_{\xi}}^2\phi\|^2+\|\na_{_{\xi}}\psi\|_{H^1}^2)d\tau
+C\d^2\int_0^t\|\na_{_{\xi}}\psi\|_{H^1}^2d\tau.
\ea
$$
Integration by parts over $\mathbb{R}\times\mathbb{T}^2$ yields
$$
K_{5,6}(t)=\sum_{i=1}^3\int_0^t\int_{\mathbb{T}^2}\int_{\mathbb{R}}(\pa_{_{\xi}}\div_{_{\xi}}\psi)^2d\xi_1d\xi'd\tau
=\int_0^t\|\na_{_{\xi}}\div_{_{\xi}}\psi\|^2d\tau.
$$
Thus, the combination of the above estimates yields
$$
\ba
K_5(t)&\leq  \f{1}{8}\int_0^t\|\sqrt{|p'(v)|}\na_{\xi}^2\phi\|^2d\tau+C\int_0^t\|\na_{\xi}^2\psi\|^2d\tau
+C \d^2 \int_0^t|\dot{\mathbf{X}}(\tau)|^2d\tau\\
&\quad+C\d^2\int_0^t(G_3(\tau)+G^s(\tau)+D(\tau))d\tau+C(\d+\chi)\int_0^t\|\na_{\xi}\psi\|^2d\tau.
\ea
$$
Using Cauchy inequality and \eqref{L3}, it holds
$$
\ba
K_6(t)&\leq C\int_0^t\int_{\mathbb{T}^2}\int_{\mathbb{R}}|\na_{_\xi}^2\phi|\big[|\na_{_\xi}^2\psi||\na_{_\xi}\phi|
+|\na_{_\xi}\psi|^2+|\na_{_\xi}^2\phi||\na_{_\xi}\psi|+|\na_{_\xi}\phi|^2\big]d\xi_1d\xi'd\tau\\
&\quad+C\int_0^t\int_{\mathbb{T}^2}\int_{\mathbb{R}}|\na_{_\xi}^2\phi|\big[|\pa_{_{\xi_1}}\phi||(u^s_1)^{-\mathbf{X}}_{_{\xi_1\xi_1}}|
+|\pa_{_{\xi_1}}\psi||(u^s_1)^{-\mathbf{X}}_{_{\xi_1}}|+|\na_{_\xi}^2\phi||(u^s_1)^{-\mathbf{X}}_{_{\xi_1}}|+|\na_{_{\xi}}\phi||(v^s)^{-\mathbf{X}}_{_{\xi_1}}|\big]d\xi_1d\xi'd\tau\\
&\leq C\int_0^t\|\na_{_\xi}^2\phi\|\big[\|\na_{_\xi}^2\psi\|_{L^6}\|\na_{_\xi}\phi\|_{L^3}
+\|\na_{_\xi}\psi\|_{L^6}\|\na_{_\xi}\psi\|_{L^3}+\|\na_{_\xi}^2\phi\|\|\na_{_\xi}\psi\|_{L^\infty}\big]d\tau\\
&\quad+C\int_0^t\|\na_{_\xi}^2\phi\|\|\na_{_\xi}\phi\|_{L^3}\|\na_{_\xi}\phi\|_{L^6}d\tau
+C\d^2\int_0^t\|\na^2_{_\xi}\phi\| (\|\na_{_{\xi}}(\phi,\psi)\|+\|\na^2_{_\xi}\phi\|)d\tau\\
&\leq C(\d+\chi)\int_0^t(\|\sqrt{|p'(v)|}\na_{_\xi}^2\phi\|^2+D(\tau)+G^s(\tau)+\|\na_{_\xi}\psi\|_{H^2}^2)d\tau,
\ea
$$
where we have use the fact that
$$
\|\na_{_{\xi}}^2\phi\|^2\|\na_{_{\xi}}\psi\|_{L^{\infty}}\leq C\|\na_{_{\xi}}^2\phi\|^2\|\na_{_{\xi}}\psi\|_{H^2}
\leq C\chi\|\na_{_{\xi}}^2\phi\|\|\na_{_{\xi}}\psi\|_{H^2}\leq C\chi(\|\sqrt{|p'(v)|}\na_{_{\xi}}^2\phi\|^2+\|\na_{_{\xi}}\psi\|_{H^2}^2).
$$
Similarly, we have
$$
\ba
K_7(t)&\leq C\int_0^t\int_{\mathbb{T}^2}\int_{\mathbb{R}}(|\na_{_\xi}\phi||\phi|
+|\phi||(v^s)^{-\mathbf{X}}_{_{\xi_1}}|+|\na_{_\xi}\phi|)|\na_{_\xi}^2\phi||(v^s)^{-\mathbf{X}}_{_{\xi_1}}|d\xi_1d\xi'd\tau\\
&\quad+ C\d\int_0^t\int_{\mathbb{T}^2}\int_{\mathbb{R}}|\phi||\pa_{_{\xi_1}}^2\phi||(v^s)^{-\mathbf{X}}_{_{\xi_1}}|d\xi_1d\xi'd\tau \\
&\leq C\d\int_0^t(\|\sqrt{|p'(v)|}\na_{_\xi}^2\phi\|^2+D(\tau)+G^s(\tau))d\tau.
\ea
$$
Using Cauchy inequality and \eqref{L3} again, we have
$$
\ba
K_8(t)&\leq C\int_0^t\int_{\mathbb{T}^2}\int_{\mathbb{R}}|\na_{_\xi}^2\psi||\na_{_\xi}^2\phi|(|\na_{_\xi}\phi|
+|(v^s)^{-\mathbf{X}}_{_{\xi_1}}|)d\xi_1d\xi'd\tau\\
&+C\int_0^t\int_{\mathbb{T}^2}\int_{\mathbb{R}}|\na_{_\xi}^2\phi|(|\na_{_\xi}\psi||\na_{_\xi}^2\phi|+|\na_{_\xi}\psi||(v^s)^{-\mathbf{X}}_{_{\xi_1\xi_1}}|)d\xi_1d\xi'd\tau\\
&\leq C\int_0^t\big[\|\na_{_\xi}\phi\|_{L^3}\|\na_{_\xi}^2\psi\|_{L^6}\|\na_{_\xi}^2\phi\|+\|\na_{_\xi}\psi\|_{L^\infty}\|\na_{_\xi}^2\phi\|^2\big]d\tau
+C\d^2\int_0^t\|\na_{_\xi}^2\phi\|\big[\|\na_{_\xi}^2\psi\|+ \|\na_{_\xi}\psi\|\big]d\tau\\
&\leq C(\d+\chi)\int_0^t\big[\|\sqrt{|p'(v)|}\na_{_\xi}^2\phi\|^2+\|\na_{_\xi}\psi\|_{H^2}^2\big]d\tau.
\ea
$$
Substituting the above estimates into \eqref{v-2nd-3} and using Lemma \ref{le-basic-new}, Lemma \ref{le-1st-v}
and Lemma \ref{le-1st-u}, we can obtain the desired inequality \eqref{es-v-2nd}. The proof of Lemma \ref{le-2nd-v} is completed.

\hfill $\Box$

%
%
%
%

\subsection{Estimates for $\|\na_{_\xi}^2(u-(u^s)^{-\mathbf{X}})\|$}

\begin{lemma}\label{le-2nd-u}
Under the hypotheses of Proposition \ref{priori}, there exists constant $C>0$ independent of $\nu$, $\d$, $\chi$ and $T$,
such that for all $t\in[0,T]$, it holds
\be\label{es-u-2nd}
\|\na_{_\xi}^2(u-(u^s)^{-\mathbf{X}})(t)\|^2+\int_0^t\|\na_{_\xi}^3(u-(u^s)^{-\mathbf{X}})\|^2d\tau
\leq C\big(\|v_0-v^s\|_{H^2}^2+\|u_0-u^s\|_{H^2}^2\big).
\ee
\end{lemma}
\textbf{\emph{Proof}}: Multiplying \eqref{dao-perturb-vu}$_2$ by $-\Delta_{_\xi}\pa_{_{\xi_i}}\psi$,
and summating $i$ from 1 to 3, then integrating the resultant equations over $\mathbb{R}\times\mathbb{T}^2$, it holds
\be\label{u-2nd}
\ba
&\f{d}{dt}\int_{\mathbb{T}^2}\int_{\mathbb{R}}\f{|\na_{_\xi}^2\psi|^2}{2}d\xi_1d\xi'
+\underbrace{\mu\int_{\mathbb{T}^2}\int_{\mathbb{R}}v|\na_{_\xi}\Delta_{_\xi}\psi|^2d\xi_1d\xi'
+(\mu+\lambda)\int_{\mathbb{T}^2}\int_{\mathbb{R}}v|\na_{_\xi}^2\div_{_\xi}\psi|^2d\xi_1d\xi'}_{\mathbf{D}_3(t)}\\
&=:\sum_{i=1}^{8}L_i(t),
\ea
\ee
where
\begin{align}\nm
L_1(t)&=-\sum_{i,j=1}^3\int_{\mathbb{T}^2}\int_{\mathbb{R}}\pa_{\xi_j}u\cdot\na_{_\xi}\pa_{_{\xi_i}}\psi\cdot\pa_{\xi_j}\pa_{_{\xi_i}}\psi d\xi_1d\xi'
+\int_{\mathbb{T}^2}\int_{\mathbb{R}}\div_{_\xi}u\f{|\na^2_{_\xi}\psi|^2}{2}d\xi_1d\xi',\\ \nm
L_2(t)&=\sum_{i=1}^3\int_{\mathbb{T}^2}\int_{\mathbb{R}}vp'(v)\Delta_{_\xi}\pa_{_{\xi_i}}\psi\cdot\na_{_\xi}\pa_{_{\xi_i}}\phi d\xi_1d\xi',\\ \nm
L_3(t)&=\sum_{i=1}^3\int_{\mathbb{T}^2}\int_{\mathbb{R}}\pa_{_{\xi_i}}u\cdot\na_{_\xi}\psi\cdot\Delta\pa_{_{\xi_i}}\psi d\xi_1d\xi'
+\sum_{i=1}^3\int_{\mathbb{T}^2}\int_{\mathbb{R}}\pa_{_{\xi_i}}(vp'(v))\Delta_{_\xi}\pa_{_{\xi_i}}\psi\cdot\na_{_\xi}\phi d\xi_1d\xi',\\ \nm
L_4(t)&=\sum_{i=1}^3\int_{\mathbb{T}^2}\int_{\mathbb{R}}\pa_{_{\xi_i}}\big(v\big(p'(v)-p'((v^s)^{-\mathbf{X}})\big)\big)(v^s)^{-\mathbf{X}}_{_{\xi_1}}\Delta_{_\xi}\pa_{_{\xi_i}}\psi_1
d\xi_1d\xi'\\ \nm
&+\int_{\mathbb{T}^2}\int_{\mathbb{R}}v\big(p'(v)-p'((v^s)^{-\mathbf{X}})\big)(v^s)^{-\mathbf{X}}_{_{\xi_1\xi_1}}\Delta_{_\xi}\pa_{_{\xi_1}}\psi_1
d\xi_1d\xi',\\ \nm
L_5(t)&=-\dot{\mathbf{X}}(t)\int_{\mathbb{T}^2}\int_{\mathbb{R}}\Delta_{_\xi}\pa_{_{\xi_1}}\psi_1(u^s_1)^{-\mathbf{X}}_{_{\xi_1\xi_1}}d\xi_1d\xi',\\\nm
L_6(t)&=\int_{\mathbb{T}^2}\int_{\mathbb{R}}vF(u^s_1)^{-\mathbf{X}}_{_{\xi_1\xi_1}}\Delta_{_\xi}\pa_{_{\xi_1}}\psi_1d\xi_1d\xi',\\ \nm
L_7(t)&=\sum_{i=1}^3\int_{\mathbb{T}^2}\int_{\mathbb{R}}\pa_{_{\xi_i}}(vF)(u^s_1)^{-\mathbf{X}}_{_{\xi_1}}\Delta_{_\xi}\pa_{_{\xi_i}}\psi_1d\xi_1d\xi',\\ \nm
L_8(t)&=(\mu+\lambda)\sum_{i=1}^3\int_{\mathbb{T}^2}\int_{\mathbb{R}}(\Delta_{_\xi}\pa_{_{\xi_i}}\psi-\na_{_\xi}\pa_{_{\xi_i}}\div_{_\xi}\psi)
\cdot\na_{_\xi} v\,\pa_{_{\xi_i}}\div_{_\xi}\psi d\xi_1d\xi',\\ \nm
&-\sum_{i=1}^3\int_{\mathbb{T}^2}\int_{\mathbb{R}}\pa_{_{\xi_i}}v(\mu\Delta_{_\xi}\psi+(\mu+\lambda)\na_{_\xi}\div_{_\xi}\psi)\cdot\Delta\pa_{_{\xi_i}}\psi d\xi_1d\xi'.
\end{align}
Using Cauchy inequality and \eqref{L3}, it holds
$$
\ba
L_1(t)&\leq C\|\na_{_\xi}\psi\|_{L^3}\|\na_{_\xi}^2\psi\|_{L^6}\|\na_{_\xi}^2\psi\|+C\d^2\|\na_{_\xi}^2\psi\|^2
\leq C(\d+\chi)\|\na_{_\xi}^2\psi\|_{H^1}^2\\
&\leq C(\d+\chi)(\mathbf{D}_3(t)+\|\na_{_\xi}^2\psi\|^2).
\ea
$$
Using Cauchy inequality, we have
$$
L_2(t)\leq \f18\mathbf{D}_3(t)+C\|\na_{_\xi}^2\phi\|^2.
$$
It follows from Cauchy inequality and \eqref{L3} that
$$
\ba
L_3(t)&\leq C\big[\|\na_{_\xi}\psi\|_{L^3}\|\na_{_\xi}\psi\|_{L^6}
+\|\na_{_\xi}\phi\|_{L^3}\|\na_{_\xi}\phi\|_{L^6}
+\d^2\|\na_{_{\xi}}(\phi,\psi)\|\big]\sqrt{\mathbf{D}_3(t)}\\
&\leq C(\d+\chi)(\mathbf{D}_3(t)+\|\na_{_\xi}\psi\|^2_{H^1}+\|\na_{_\xi}\phi\|^2_{H^1})\\
&\leq C(\d+\chi) (\mathbf{D}_3(t)+G^s(t)+D(t)+\|\na_{_\xi}\psi\|_{H^1}^2+\|\na_{_\xi}^2\phi\|^2).
\ea
$$
Using Lemma \ref{le-shock}, we have
$$
\ba
L_4(t)&\leq C\int_{\mathbb{T}^2}\int_{\mathbb{R}}\big(|\na_{_\xi}\phi|+|(v^s)^{-\mathbf{X}}_{_{\xi_1}}||\phi|\big)|(v^s)^{-\mathbf{X}}_{_{\xi_1}}||\na^3_{_\xi}\psi_1|d\xi_1d\xi'
+C\d\int_{\mathbb{T}^2}\int_{\mathbb{R}}|\phi||(v^s)^{-\mathbf{X}}_{_{\xi_1}}||\na_{_\xi}^3\psi_1|d\xi_1d\xi'\\
&\leq C\d(\mathbf{D}_3(t)+G^s(t)+\|\na_{_\xi}\phi\|^2)\leq C\d(\mathbf{D}_3(t)+G^s(t)+D(t)),
\ea
$$
and
$$
L_5(t)\leq C\d|\dot{\mathbf{X}}(t)|\sqrt{\mathbf{D}_3(t)}\, \|(v^s)^{-\mathbf{X}}_{_{\xi_1}}\|_{L^2(\mathbb{R})}
\leq C\d^{\f52}|\dot{\mathbf{X}}(t)|\sqrt{\mathbf{D}_3(t)}\leq C\d^{\f52}\mathbf{D}_3(t)+C\d^{\f52}|\dot{\mathbf{X}}(t)|^2.
$$
Using the definition of $F$ \eqref{F}, it holds
$$
\ba
L_6(t)&\leq C\d\int_0^t\int_{\mathbb{T}^2}\int_{\mathbb{R}}\Big(\Big|h_1-(h^s_1)^{-\mathbf{X}}-\f{p(v)-p((v^s)^{-\mathbf{X}})}{\sigma_*}\Big|
+|p(v)-p((v^s)^{-\mathbf{X}})|\\
&\quad+\big|\pa_{_{\xi_1}}\big(p(v)-p((v^s)^{-\mathbf{X}})\big)\big|+|(v^s)^{-\mathbf{X}}_{_{\xi_1}}||p(v)-p((v^s)^{-\mathbf{X}})|\Big)
|(v^s)^{-\mathbf{X}}_{_{\xi_1}}||\Delta_{_\xi}\pa_{_{\xi_1}}\psi_1|d\xi_1d\xi'd\tau\\
&\leq C\d^2(\sqrt{G_3(t)}+\sqrt{G^s(t)}+\sqrt{D(t)})\|\Delta\pa_{_{\xi_1}}\psi_1\|
\leq C\d^2(G_3(t)+G^s(t)+D(t)+\mathbf{D}_3(t)).
\ea
$$
Using \eqref{vF} and  Cauchy inequality, it holds
$$
\ba
L_7(t)\leq C\d^2\|\na_{_\xi}(\phi,\psi_1)\|\|\Delta_{_\xi}\pa_{_{\xi_1}}\psi_1\|
\leq C\d^2(\mathbf{D}_3(t)+G^s(t)+D(t)+\|\na_{_\xi}\psi_1\|^2).
\ea
$$
Using Cauchy inequality and \eqref{L3}, we have
$$
\ba
L_8(t)&\leq C\sqrt{\mathbf{D}_3(t)}\|\na_{_\xi}\phi\|_{L^3}(\|\Delta_{_\xi}\psi\|_{L^6}+\|\na_{_\xi}\div_{_\xi}\psi\|_{L^6})
+C\d^2\sqrt{\mathbf{D}_3(t)}(\|\Delta_{_\xi}\psi\|+\|\na_{_\xi}\div_{_\xi}\psi\|)\\
&\leq C(\d+\chi)\sqrt{\mathbf{D}_3(t)}(\|\Delta_{_\xi}\psi\|_{H^1}+\|\na_{_\xi}\div_{_\xi}\psi\|_{H^1})\\
&\leq C(\d+\chi)(\mathbf{D}_3(t)+\|\Delta_{_\xi}\psi\|_{H^1}^2+\|\na_{_\xi}\div_{_\xi}\psi\|_{H^1}^2).
\ea
$$
Substituting the above estimates into \eqref{u-2nd}, and integrating the resultant equations over $[0,t]$ for any $t\leq T$,
using Lemmas \ref{le-1st-v}, \ref{le-1st-u} and \ref{le-2nd-v}, we can get the desired inequality \eqref{es-u-2nd}. The
proof of Lemma \ref{le-2nd-u} is completed.

\hfill $\Box$

\subsection{Proof of Proposition \ref{priori}}

We use \eqref{daoshu} to have
$$
\|\na_{_\xi}(v-(v^s)^{-\mathbf{X}})\|^2\leq C( D(t)+G^s(t)),
$$
which together with  Lemmas \ref{le-basic-new}-\ref{le-2nd-u} yields \eqref{full-es}.
In addition, using \eqref{X} and \eqref{new-shock-2}$_2$, and together with Lemma \ref{le-shock} and the assumption \eqref{assumption}, we have
$$
|\dot{\mathbf{X}}(t)|\leq \f{C}{\d}(\|\big(p(v)-p((v^s)^{-\mathbf{X}})\big)\|_{L^{\infty}}
+\|v-(v^s)^{-\mathbf{X}}\|_{L^{\infty}})\int_{\mathbb{T}^2}\int_{\mathbb{R}}(v^s)^{-\mathbf{X}}_{\xi_1}d\xi_1d\xi'
\leq C\|v-(v^s)^{-\mathbf{X}}\|_{L^{\infty}},
$$
which implies \eqref{point-X}. The proof of Proposition \ref{priori}
is completed.

\hfill $\Box$

%
%
%
%

%
%
%

%
%
%

\noindent\textbf{Acknowledgment.} The research of T. Wang is partially supported by NSFC grant No. 11971044 and BJNSF grant No. 1202002.
The research of Y. Wang is partially supported by the NNSFC grants No. 12090014 and 11688101.

\end{document}